\documentclass{siamltexmm}

\usepackage{amssymb,bm,color}
\usepackage[leqno]{amsmath}
\usepackage[font=normalsize]{subfig}
\newcommand{\rd}{\mathrm{d}}

\title{A Study in Three-Dimensional Chaotic Dynamics:\\ Granular Flow and Transport in a Bi-Axial Spherical Tumbler\thanks{Received by the editors...}} 

\author{Ivan~C. Christov\footnotemark[2]\ \footnotemark[3]\ \footnotemark[4]\ \footnotemark[7]
\and Richard~M. Lueptow\footnotemark[5]\ \footnotemark[7] 
\and Julio~M. Ottino\footnotemark[6]\ \footnotemark[7]
\and Rob Sturman\footnotemark[8]}

\begin{document}
\maketitle

\renewcommand{\thefootnote}{\fnsymbol{footnote}}

\footnotetext[2]{Department of Engineering Sciences and Applied Mathematics, Northwestern University, Evanston, IL 60208, USA (\email{christov@u.northwestern.edu}). Supported in part by a Walter P.\ Murphy Fellowship from the Robert R.\ McCormick School of Engineering and Applied Science.}
\footnotetext[3]{Department of Mechanical and Aerospace Engineering, Princeton University, Princeton, NJ 08544, USA. Supported in part by NSF Grant DMS-1104047.}
\footnotetext[4]{Present address: Theoretical Division and Center for Nonlinear Studies, Los Alamos National Laboratory, Los Alamos, NM 87545, USA.}
\footnotetext[5]{Department of Mechanical Engineering and The Northwestern Institute on Complex Systems (NICO), Northwestern University, Evanston, IL 60208, USA (\email{r-lueptow@northwestern.edu}).}
\footnotetext[6]{Department of Chemical and Biological Engineering and The Northwestern Institute on Complex Systems (NICO), Northwestern University, Evanston, IL 60208, USA (\email{jm-ottino@northwestern.edu}).}
\footnotetext[7]{The work of these authors was supported in part by NSF Grant CMMI-1000469.}
\footnotetext[8]{Department of Applied Mathematics, University of Leeds, Leeds LS2 9JT, UK  (\email{r.sturman@maths.leeds.ac.uk}).}

\renewcommand{\thefootnote}{\arabic{footnote}}

\newcommand{\slugmaster}{%
\slugger{MMedia}{xxxx}{xx}{x--x}}

\begin{abstract}
We study three-dimensional (3D) chaotic dynamics through an analysis of transport in a granular flow in a half-full spherical tumbler rotated sequentially about two orthogonal axes (a bi-axial ``blinking'' tumbler). The flow is essentially quasi-two-dimensional in any vertical slice of the sphere during rotation about a single axis, and we provide an explicit exact solution to the model in this case. Hence, the cross-sectional flow can be represented by a twist map, allowing us to express the 3D flow as a linked twist map (LTM). We prove that if the rates of rotation about each axis are equal, then (in the absence of stochasticity) particle trajectories are restricted to two-dimensional (2D) surfaces consisting of a portion of a hemispherical shell closed by a ``cap''; if the rotation rates are unequal, then particles can leave the surface they start on and traverse a volume of the tumbler. The period-one structures of the governing LTM are examined in detail: analytical expressions are provided for the location of period-one curves, their extent into the bulk of the granular material, and their dependence on the protocol parameters (rates and durations of rotations). Exploiting the restriction of trajectories to 2D surfaces in the case of equal rotation rates about the axes, a method is proposed for identifying and constructing 3D Kolmogorov--Arnold--Moser (KAM) tubes around the normally elliptic period-one curves. The invariant manifold structure arising from the normally hyperbolic period-one curves is also examined. When the motion is restricted to 2D surfaces, the structure of manifolds of the hyperbolic points in the bulk differs from that corresponding to hyperbolic points in the flowing layer. Each is reminiscent of a template provided by a non-integrable perturbation to a Hamiltonian system, though the governing LTM is not. This highlights the novel 3D chaotic behaviors observed in this model dynamical system.
\end{abstract}



\pagestyle{myheadings}
\thispagestyle{plain}
\markboth{I.~C. CHRISTOV, R.~M. LUEPTOW, J.~M. OTTINO, AND R. STURMAN}{THREE-DIMENSIONAL CHAOS IN GRANULAR FLOW}

\section{Introduction}

Chaotic dynamics in higher dimensions remains relatively unexplored. The great variety of transport and mixing behaviors of three-dimensional (3D) flows, many of which are distinct from those in two-dimensional (2D) systems, have been studied only sparingly \cite{wig-fof} and even fewer studies \cite{emc93} have considered higher-than-three dimensional dynamics. The study presented here makes clear why. The system considered---an idealization of a laboratory granular mixing flow---is, superficially, simple, quite possibly one of the simplest 3D flows. But the wealth of complexities associated with analyzing 3D dynamics, even when restricting ourselves to the simplest member of the family, a half-full case, demonstrates that the need to develop mathematics for such 3D chaotic dynamics can keep researchers busy for a long time to come. We discuss novel 3D chaotic behaviors exhibited in the \emph{blinking spherical tumbler flow}, a 3D chaotic flow of granular materials in a bi-axial spherical tumbler \cite{jfg03,mlo07,smow08,jlosw10,jcol12}. The name comes from an analogy to an early example of 2D chaotic fluid mixing---the blinking vortex flow \cite{a84}---because a 3D spherical tumbler sequentially rotated about two distinct axes ``blinks'' between the two rotations. Such a mode of operation is distinct from, e.g., a \emph{transitory} dynamical system in which the transition between two types of dynamics occurs smoothly in finite time \cite{mm11}; here, the transition is instantaneous.

Early work explored quasi-3D flows (specifically, 2D steady flows with a periodically repeating cross-section in the third dimension), which were easily analyzed using the techniques applicable to 2D time-periodic flows \cite{o89}. This approach has proven fruitful leading to many examples, including but not limited to, the partitioned-pipe mixer \cite{kfo87}, the twisted-pipe mixer \cite{jta89}, the Kenics mixer \cite{gapm03},  mixing tanks with impellers \cite{lam99,hok12}, the staggered herringbone mixer \cite{sdamsw02}, screw extrusion mixers \cite{hkk04}, and the rotated arc mixer \cite{grbgh06}. 

However, even earlier, Arnold \cite{a65} (followed by H\'enon \cite{h66}) showed that 3D chaotic fluid particle trajectories occur in a spatially periodic idealized fluid flow [the so-called in Arnold--Beltrami--Childress (ABC) flow], thus motivating another set of early work on 3D chaotic dynamics \cite{dfghms86,fkp88,betal89,hk91,zklh93}. Later, a more realistic and bounded 3D flow, a poloidal Hill vortex, was considered \cite{bm90,sns91,ks99,bb99} as a model system for mixing within droplets. Mixing in the 3D flow induced by the motion of a droplet along a serpentine channel has also been  examined \cite{sti03,ss05,cvcsa08} more recently. Cartwright {\it et al.}~\cite{cfp96} showed the existence of chaotic fluid particle trajectories in a spherical Couette flow. Meanwhile, Ashwin and King~\cite{ak95} analyzed 3D fluid motion in a Taylor--Couette flow between eccentric cylinders, and chaotic particle trajectories were found in a wavy-vortex Taylor--Couette flow using both simulated \cite{r98} and experimental \cite{al02} velocity fields. Soon, laboratory studies of the \emph{fundamental mixing tank} (an upright cylinder with an angled impeller driving the flow from the top) \cite{fkmo} provided experimental verification of the theoretical results. Cavity flows have also proven attractive for studying 3D mixing. Anderson {\it et al.}~\cite{agpvm99} considered a box whose side-walls are allowed to translate in two coordinate directions, while Speetjens {\it et al.}~\cite{mmsch02,sch04,psc10} studied an upright cylinder whose entire top (and/or bottom) wall is allowed to execute arbitrary motions. Clearly, the notion of 3D chaotic advection has proven quite attractive for the analysis of fluid mixing devices. The 3D dynamical systems framework has even found applications in pollution transport in estuaries \cite{stir03}. 

In contrast, the detailed mathematical study of 3D granular flows has only recently begun. What appears to be one of the first quasi-3D granular flows (specifically, rocking motions of a rotating cylindrical container that give a weak three-dimensionality to the flow) was studied experimentally and numerically (using the discrete element method) in \cite{WightmanSim}. Meanwhile, 3D effects in granular flows in rotating cylindrical and spherical tumblers have been analyzed experimentally \cite{go03}. Through a dynamical systems framework, it has become possible to gain a deeper understanding of how mixing and segregation (demixing) manifest themselves in 3D granular flows \cite{mlo07}. Most research on 3D granular flows, however, has remained on the practical side, exploring mixing of powders in industrial devices \cite{lbcg07,dbc08}, rather than focusing on the underlying mathematics of 3D chaotic dynamics.

Once it had been established that fluid particle trajectories in non-turbulent 3D fluid flows could be chaotic, the natural question arose: What are the possible types of dynamics (regular, chaotic or, perhaps, other) that can be observed in 3D? Whereas a 2D non-integrable Hamiltonian flow can, roughly speaking, be decomposed into Kolmogorov--Arnold--Moser (KAM) islands (regular regions) and tangles of manifolds (chaotic regions), there is no similar obvious classification for a 3D chaotic flow. Nevertheless, regular and chaotic regions still exist \cite{fkp88}, which has suggested (at least) two types of classifications for 3D chaotic flows based on either structure of invariant tori by analogy to the KAM theorem \cite{cfp96} or local flow structure near periodic points \cite{mmsch02}. For example, in certain 3D chaotic flows, it is possible to have chaotic regions restricted to 2D spheroidal surfaces, while KAM islands, which by construction do not admit chaotic trajectories, may have a 3D character manifesting themselves as KAM-like \emph{tubes} \cite{dfghms86,kfo87,fkp88,ko92} in spatially periodic domains and KAM \emph{tori} \cite{fkmo} in bounded domains. 

Many of the dynamical systems considered so far were chosen, in part, because of their simplicity or due to their representation as a perturbation of an integrable Hamiltonian template. Such constructions have been the subject of a number of studies, concurrent with studies of chaotic 3D fluid particle motion, on the generic (mathematical) properties of 3D dynamical systems \cite{fkp88,cs90,mackay,bajer,mw94,vm12}. Specifically, the question of invariants of the motion and their classification has been an important topic \cite{hm98,gm,mjm05}. The notion of \emph{action-angle} variables and, specifically, \emph{action-action-angle} variables as we discuss below, for 3D dynamical systems are also a useful concept \cite{cfp96,vwg08,vwg10,vnm06,smow08} to prove rigorous results and to better characterize invariant structures. In the present work, on the other hand, we explore the 3D chaotic dynamics emerging from a \emph{linked twist map} (LTM) \cite{sow06}. Though the LTM can be constructed from the composition of two fully integrable action--action--angle maps, the combined map is capable of producing 3D chaotic dynamics.

Recently, the question of 3D coherent structures created by chaotic advection has been further explored in some experimentally realizable flows \cite{sch06a,sch06,psc10,zstc12,sc13}. For example, one new (uniquely 3D) phenomenon emerging from the latter studies is \emph{resonance-induced merger} of coherent structures. Put simply, in the parameter regime of the system that is only a slight perturbation of the regime in which all motion is constrained to invariant surfaces, it is possible to find particle trajectories that still do not traverse the whole volume of the domain. They remain restricted to \emph{two} of the previously invariant surfaces, which have now become fused through a ``channel'' surrounding a parabolic periodic point of the flow. 

Naturally, the study of coherent structures in 3D chaotic flows is underpinned by the understanding of fixed and periodic points of the flow map. In 3D, fixed and periodic points can be either isolated or part of continuous segments, thus classifications have been developed for both situations \cite{mmsch02}. While curves of periodic (or fixed) points lead to essentially 2D dynamics on families of surfaces intersecting said curves, isolated periodic points can lead to genuinely 3D dynamics \cite{sch04,sc13}. Analysis of fundamental symmetries \cite{fo92} can provide guidance in locating periodic (or fixed) points but locating all periodic (or fixed) points of a 3D flow map is generally a nontrivial task. The blinking spherical tumbler flow that we consider is defined by the pairwise composition of rotations about orthogonal axes. This fact provides a natural structure to the 3D flow map, which can be exploited in locating period-one points. Specifically, a sufficient (but not necessary) condition that a point be a period-one point of the flow map is that it is a period-one point of each rotation. This is the approach we follow below. However, in the case that the individual rotations have no non-trivial period-one points, an isolated period-one point of the composite flow map is still guaranteed by Brouwer's fixed-point theorem \cite{b11} (see, e.g., the example of the blinking cubical tumbler flow in \cite{mlo07}). We do not focus on isolated period-one points in this study, but seek to analyze in full the simplest periodic structures related to the symmetries in the system.

Here, we explore 3D chaotic dynamics in the physical system of a blinking spherical tumbler half-filled with granular material. Most granular flows encountered in nature and industry are 3D. However, the modeling and mathematical theory of 3D flows, whether granular or fluid, is far from straightforward. What is more, visualization of 3D flows remains a challenge \cite{cyymm11}. Previous work on the kinematics of granular flow \cite{hkgmo99,mclo06,col10} has shown that there is already a wealth of interesting phenomena occurring in quasi-2D granular flows. Thus, it is expected that 3D granular flows will likely reveal even more striking representative features of complex systems far from equilibrium. To this end, motivated by the successful application of dynamical systems theory to granular flow in a blinking spherical tumbler initiated in \cite{mlo07,smow08}, the present work focuses on  providing a more thorough treatment of the 3D features of this flow and, in particular, extending the results from \cite{smow08}.

A difficulty for previous studies of experimentally realizable 3D chaotic flows (e.g., \cite{sch06a,sch06,psc10,zstc12,sc13}) has been the inability to perform any of the calculations required to locate period-one structures or invariant surfaces and determine their dependence on the parameters analytically. However, as shown in \S\ref{sec:radial_disp} and \ref{sec:period_one}, these calculations can be done for the granular flow considered in the present work. Moreover, the blinking spherical tumbler flow is realizable in the laboratory \cite{mlo07,smow08,jlosw10}, which opens the door for future experimental verification of the theoretical results presented herein.

\section{Continuum model of granular flow in a blinking spherical tumbler}
\label{sec:3d_const_ang}

In this section, we summarize the continuum model \cite{jfg03,mlo07} for the kinematics of granular flow in the continuous-flow (rolling) regime for a half-full spherical tumbler rotated alternately about two orthogonal axes [for the present purposes, the $z$- and $x$-axes as shown in Fig.~\ref{fig:biaxsphere}(a)]. A variety of other continuum models of 3D granular flow exist. More details can be found in recent reviews on the topic \cite{mlo07,at09,k11}.

\subsection{Rotation about the $z$-axis}
\label{sec:rot_z_axis}

\begin{figure}
\centerline{\includegraphics[width=0.8\textwidth]{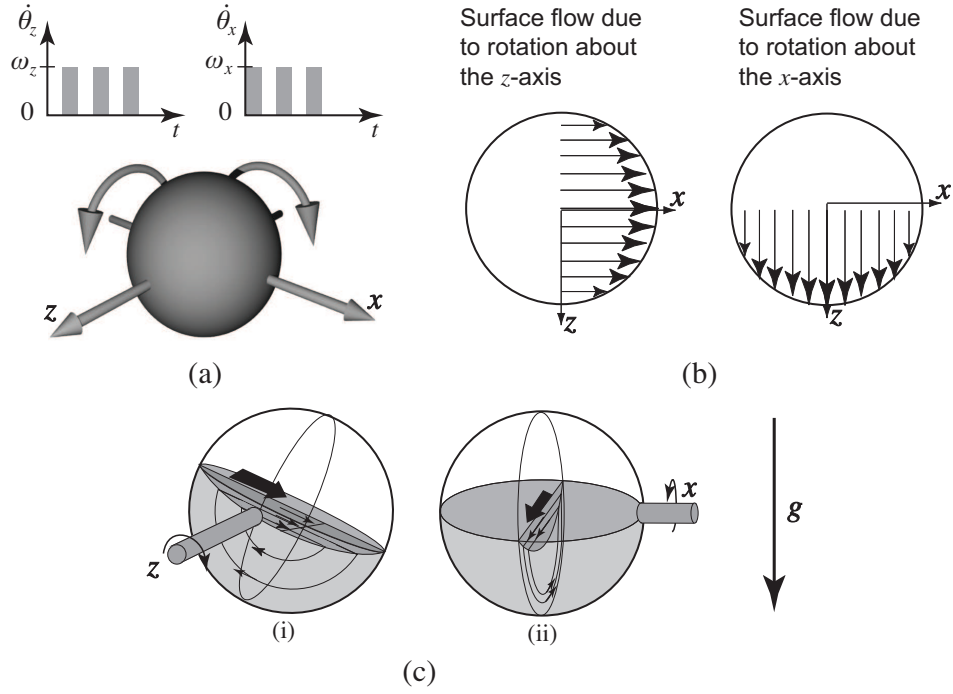}}
\caption{Schematic of the blinking protocol for a spherical mixer (a). The $z$- and $x$-axis are coplanar and orthogonal to each other, intersecting at the origin of the sphere. Top view of the surface velocity profile due to rotation around each axis (b). The flow fields in the fixed cross-sections for clockwise rotation around the $z$- and $x$-axis, respectively, are shown in panels (c)(i) and (c)(ii). Adapted from \cite{smow08} \copyright\ 2008 Cambridge University Press with permission.}
\label{fig:biaxsphere}
\end{figure}

Without loss of generality, we take the first rotation of the half-full spherical tumbler of radius $R$ to be clockwise about the $z$-axis (denoted by the $z$ subscripts) at a rate of $\omega_z(>0)$ for a duration $\tau_z$, typically such that the angle of rotation $\theta_z = \omega_z\tau_z$ is of the order of $\pi$ radians or less. The coordinate system used is Cartesian with the origin at the center of the sphere, with the flow being in the $x$-direction in this case, and the transverse direction being the $z$-direction [see the left portion of Fig.~\ref{fig:biaxsphere}(b)]. Under this version of the continuum model, the velocity in the flow direction, $v_{z,x}$, is assumed to vary linearly with the depth $y$ in the flowing layer, which has thickness $\delta_z(x,z)$ and free-surface half-length $L(z)$, so that the shear rate $\dot{\gamma_z} = \partial v_{z,x}/\partial y$ is constant. The remaining bulk of the granular material is in solid body rotation. In an experiment, the flowing layer is at an angle, referred to as the dynamic angle of repose, with respect to the horizontal. In the model, it is assumed that the coordinate system is rotated backwards by this fixed angle. Also, we assume the transport in the transverse direction is negligible, a reasonable assumption over the rotation angles considered here \cite{zduol12}.

Let us introduce the following dimensionless variables:
\begin{equation}
\{x^\star, y^\star,  z^\star,  L^\star,  \delta^\star\} = \{ x,y,z,L,\delta\}/R,\qquad t^\star = \omega_z t,\qquad \psi^\star = \psi/(R^2\omega_z).
\label{eq:nd_vars_z}
\end{equation}
Then, leaving the $\star$ superscript understood, the equations of pathlines (i.e., for the advection of a passive scalar) in the tumbler \cite{mlo07,smow08} are
\begin{subequations}\label{eq:vel_z}\begin{align}
\frac{\rd x}{\rd t} &= v_{z,x}(x,y,z) = \begin{cases} \displaystyle\frac{1}{\epsilon_z^2} \left[\delta_z(x,z)+y\right], &\quad y>-\delta_z(x,z),\\ y, &\quad\text{otherwise},\end{cases}\label{eq:vel_z_dxdt}\displaybreak[3]\\
\frac{\rd y}{\rd t} &= v_{z,y}(x,y,z) = \begin{cases} \displaystyle\frac{x\, y}{\delta_z(x,z)}, &\quad y>-\delta_z(x,z),\\-x, &\quad\text{otherwise},\end{cases}\label{eq:vel_z_dydt}\displaybreak[3]\\
\frac{\rd z}{\rd t} &= v_{z,z}(x,y,z) = 0,\label{eq:vel_z_dzdt}
\end{align}\end{subequations}
where $\epsilon_z$ is defined below. Here, the assumptions of incompressibility and no transport in the transverse direction allows us to introduce a streamfunction $\psi$  \cite[\S3.12]{o89} in terms of the $x$ and $y$ coordinates alone:
\begin{equation}
\psi(x,y,z) = \begin{cases}\psi_\mathrm{fl} = \displaystyle{\frac{1}{\epsilon_z^2}}\left[\delta_z (x, z) y + \tfrac{1}{2}y^2 \right], &\quad y>-\delta_z(x,z),\\[3mm] \psi_\mathrm{b} = \frac{1}{2}(x^2+y^2), &\quad\text{otherwise}.\end{cases}
\label{eq:streamfunction}
\end{equation}
Clearly, $v_{z,x} = \partial\psi/\partial x$ and $v_{z,y} = - \partial\psi/\partial y$. It can be shown \cite[\S2.1.1]{smow08} that the streamfunction can be made continuous, if necessary, across the interface $y=-\delta_z(x,z)$.

The shape of the flowing layer is given by
\begin{equation}
\delta_z(x,z) = \delta_{0}(z) \sqrt{1-\frac{x^2}{L(z)^2}} = \epsilon_z\sqrt{1 - x^2 - z^2},
\label{eq:delta_z}
\end{equation}
where the second equality follows from the relations $\delta_{0}(z) = L(z) \sqrt{{\omega_z}/{\dot{\gamma_z}}}$ and $L(z) = \sqrt{1-z^2}$ for the \emph{local} (dimensionless) depth and half-length of the flowing layer, respectively. Since $L(0) = 1$, we have $\epsilon_z := \delta_{0}(0)/L(0) = \delta_0(0)$, meaning $\epsilon_z$ is the \emph{maximal} (dimensionless) depth (at $x=z=0$). Of course, $\epsilon_z \equiv \sqrt{{\omega_z}/{\dot{\gamma_z}}}$, and it coincides with the definition of the small parameter in the model from \cite{col10}, except $L$ and $\delta_0$ now depend on $z$ rather than on $t$. Here we note that, physically, a restriction $\epsilon_z \lesssim 0.5$ can be enforced as it does not make sense for the flowing layer to be thicker than (about) half the tumbler's radius ($=1$ in dimensionless units). This assumption is supported by the available experimental data \cite{jol04,ffo07,pakcol12}; in fact, typically $\epsilon_z  \simeq 0.1$.

\subsection{Rotation about the $x$-axis}
\label{sec:rot_x_axis}

The blinking spherical tumbler protocol alternates between rotations about two orthogonal axes [recall Fig.~\ref{fig:biaxsphere}(a)]. Once the rotation about the $z$-axis is complete, a rotation about the $x$-axis at constant rate $\omega_x$ for a duration $\tau_x$ follows. In a laboratory experiment, the tumbler would need to be rotated backward by the static angle of repose to make the free surface horizontal, before changing the axis of rotation, and then rotated forward by the static angle or repose, before rotation and flow about the second axis occurs. Both intermediate solid body rotation steps are performed in a way such that the material remains in static equilibrium and, consequently, there is no flow. As far as the mathematical analysis goes, however, the intermediate solid body rotation steps are inconsequential. Thus, since the rotations are completely independent, we can make $t$ and $\psi$ dimensionless for rotation about the $x$-axis using $\omega_x$. All other quantities are made dimensionless as in Eq.~\eqref{eq:nd_vars_z}. 

\begin{figure}
\centerline{\includegraphics[width=0.6\textwidth]{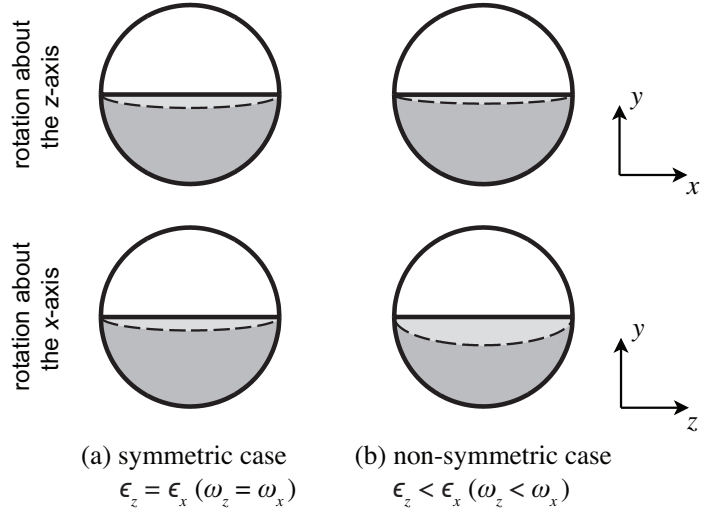}}
\caption{Each row is a schematic representation of a cross-section perpendicular to the $z$- or $x$-axis of rotation of the spherical tumbler. In the symmetric case (a), the flowing layer thickness is the same during each rotation. For the non-symmetric case (b), the flowing layer is thinner for rotation about the $z$-axis than for rotation about the $x$-axis because rotation rates are unequal, specifically $\omega_z < \omega_x$.}
\label{fig:sym_nonsymm}
\end{figure}

Similar to $z$-axis rotation, the interface between the flowing layer and the bulk is given by
\begin{equation}
\delta_x(x,z) = \epsilon_x\sqrt{1-x^2-z^2}\,,
\label{eq:delta_x}
\end{equation}
where $\epsilon_x = \sqrt{\omega_x/\dot\gamma_x}$. Note that if $\epsilon_z = \epsilon_x$, then $\delta_z(x,z) = \delta_x(x,z)$. We shall refer to this as the \emph{symmetric case} in the discussion below. Then, $\epsilon_z \ne \epsilon_x$ can be termed the \emph{non-symmetric case}. It is important to note that different values for $\epsilon_z$ and $\epsilon_x$ arise in physical systems due to different rotation rates $\omega_z$ and $\omega_x$ \cite{jol04,ffo07,pakcol12}. Thus, the symmetric case occurs for identical rotation rates about the two axes of rotation, while the non-symmetric case arises from different rotation rates. This is illustrated in Fig.~\ref{fig:sym_nonsymm}.

It should be clear that the discussion from \S\ref{sec:rot_z_axis} can be repeated here as well (after formally exchanging ``$x$'' and ``$z$'' in all equations and relations) to obtain the velocity field corresponding to rotation about the $x$-axis.

\section{A necessary condition for radial transport}
\label{sec:radial_disp}

In this section, we consider the trajectory of a single tracer particle in the blinking spherical tumbler flow described in the previous section. The results in this section establish the key difference between the symmetric and non-symmetric cases. 

In the symmetric case, i.e., when $\epsilon_z = \epsilon_x$, tracer particle trajectories are constrained to surfaces determined by the tracer's initial position \cite{smow08}. The proof \cite[Appendix~B]{smow08} follows from the fact that if $\epsilon_z = \epsilon_x$, then $\delta_z(x,z) \equiv\delta_x(x,z)$, which endows the dynamical system with a ``symmetry'' and, hence, an invariant structure. On the other hand, for the non-symmetric case ($\epsilon_z \ne \epsilon_x$), we show below that trajectories are not necessarily restricted to invariant surfaces, thereby allowing for radial transport in the flow or, equivalently, transport \emph{across} invariant surfaces (not just \emph{within} them).

Let us violate the conditions of \cite[Theorem~3.2]{smow08} by allowing $\epsilon_z\ne\epsilon_x$. Then, consider a tracer initially located at the point $(x_0,y_0,z_0)^\top$ in the bulk with initial radius $r_i$, i.e., $x_0^2 + y_0^2 + z_0^2 = r_i^2$. During rotation about the $z$-axis, the trajectory is restricted to the plane $z=z_0$ and along the curve given implicitly by 
\begin{equation}
\psi_\mathrm{b} = \frac{1}{2}(x^2+y^2) = const. = \frac{1}{2}(x_0^2 + y_0^2) = \frac{1}{2}(r_i^2-z_0^2).
\end{equation}
Thus, all points $(x,y,z_0)^\top$ on this streamline in the bulk verify $x^2 + y^2 + z_0^2 = r_i^2$. A similar argument holds for rotation about the $x$-axis. Therefore, under this model, a tracer particle in the bulk at a distance $r_i$ from the origin remains a distance $r_i$ from the origin when the axis of rotation is switched while the particle is still in the bulk.

Consider the same tracer particle as above, but for the case when it reaches the flowing layer during the first rotation. The point at which the particle reaches the flowing layer, $(x_1,y_1,z_1=z_0)^\top$, must lie on the intersection of the curves
\begin{subequations}\begin{align}
r_i^2 &= x_1^2+y_1^2 + z_1^2 = x_1^2+y_1^2 + z_0^2, \\
y_1 &= -\delta_z(x_1,z_0) = -\epsilon_z\sqrt{1 - x_1^2 - z_0^2}\,.
\end{align}\end{subequations}
The solution to these equations for $x_1$ and $y_1$, where we must take the solution with $y_1 < 0$, corresponds to the point where this trajectory enters the flowing layer:
\begin{equation}
x_1 = -\sqrt{\frac{r_i^2 - \epsilon_z^2 }{1-\epsilon_z^2} - z_0^2},\qquad
y_1 = -\epsilon_z\sqrt{\frac{1-r_i^2}{1-\epsilon_z^2}}.
\label{eq:z_rot_x1}
\end{equation}
In the flowing layer, the tracer follows a fixed streamline as the tumbler continues to rotate about the $z$-axis. Using Eqs.~\eqref{eq:z_rot_x1}, this streamline can be found to be
\begin{equation}
\psi_\mathrm{fl} = \frac{1}{\epsilon_z^2}\left[\delta_z (x, z) y + \tfrac{1}{2}y^2 \right] = const. = \frac{1}{\epsilon_z^2}\left[\delta_z (x_1, z_0) y_1 + \tfrac{1}{2}y_1^2 \right] = \frac{r_i^2 - 1}{2(1 - \epsilon_z^2)}. 
\label{C_2}
\end{equation}

Now, suppose that rotation about the $z$-axis stops while the tracer is still in the flowing layer at its final position $(x_2,y_2,z_2=z_0)^\top$, which is a function of $\epsilon_z$ and $\theta_z$ (the angle through which the tumbler was rotated about the $z$-axis). Further, assume the depth of the flowing layer corresponding to rotation about the $x$-axis is such that the tracer remains in the flowing layer upon switching axes of rotation. Then, rotating about the $x$-axis, $x=x_2$ remains fixed, and the tracer's motion is along the streamline in the flowing layer given by
\begin{equation}
\psi_\mathrm{fl} = \frac{1}{\epsilon_x^2}\left[\delta_x (x, z) y + \tfrac{1}{2}y^2 \right] = const. = \frac{1}{\epsilon_x^2}\left[\delta_x (x_2, z_0) y_2 + \tfrac{1}{2}y_2^2 \right] =: \kappa.
\label{eq:psi_fl_x_rot}
\end{equation}
Note that, at the moment at which the axes of rotation is switched, we must have \emph{both}
\begin{subequations}\label{eq:psis_switch_pt}\begin{align}
\frac{1}{\epsilon_z^2}\left[\delta_z (x_2, z_0) y_2 + \tfrac{1}{2}y_2^2 \right] &= \frac{r_i^2 - 1}{2(1 - \epsilon_z^2)},\label{eq:psis_switch_pt1}\\
 \frac{1}{\epsilon_x^2}\left[\delta_x (x_2, z_0) y_2 + \tfrac{1}{2}y_2^2 \right] &= \kappa.\label{eq:psis_switch_pt2}
\end{align}\end{subequations}
Subtracting Eq.~\eqref{eq:psis_switch_pt1} from Eq.~\eqref{eq:psis_switch_pt2}, we find that
\begin{equation}
\kappa = \kappa(x_2,y_2,z_0,\epsilon_x,\epsilon_z) = \frac{r_i^2 - 1}{2(1 - \epsilon_z^2)} + \left(\frac{1}{\epsilon_x^2} - \frac{1}{\epsilon_z^2}\right)\tfrac{1}{2}y_2^2 + \left(\frac{1}{\epsilon_x} - \frac{1}{\epsilon_z}\right)y_2\sqrt{1-x_2^2-z_0^2}.
\label{eq:kappa_y2_ex_ez}
\end{equation}
Therefore, in the symmetric case when $\epsilon_z = \epsilon_x$ ($\Rightarrow\delta_z\equiv\delta_x$), we have $\kappa = (r_i^2-1)/[2(1-\epsilon_z^2)]$. In principle, one can find $x_2$ and $y_2$ exactly by tracing the flowing layer trajectory using the solution given in Appendix~\ref{sec:exact_2d_sol}. Then, $\kappa$ can be computed as well.

Now, suppose that the particle follows the streamline defined through Eq.~\eqref{eq:psi_fl_x_rot} under rotation about the $x$-axis until it reaches the boundary between the bulk and the flowing layer at the point $(x_3 = x_2, y_3, z_3)^\top$, where $y_3$ and $z_3$ are the solutions of
\begin{equation}
y_3 = -\delta_x(x_2,z_3),\qquad
\kappa = \frac{1}{\epsilon_x^2}\left[\delta_x (x_2, z_3) y_3 + \tfrac{1}{2}y_3^2 \right].
\label{eq:y3_k}
\end{equation}
Noting that $z,y,\kappa<0$ when a particle is exiting the flowing layer, the system in Eq.~\eqref{eq:y3_k} can be solved to give:
\begin{equation}
y_3 = -\sqrt{-2\epsilon_x^2 \kappa}, \qquad
z_3 = -\sqrt{1 - x_2^2 + 2\kappa}.
\label{eqs:exit}
\end{equation}

Finally, we can compute the final radial distance from the origin at which the particle exits the flowing layer. This is accomplished by using Eqs.~\eqref{eqs:exit}:
\begin{equation}
x_3^2 + y_3^2 + z_3^2 = x_2^2 - 2\epsilon_x^2\kappa + 1 - x_2^2 + 2\kappa = 1 + 2\kappa \left( 1 - \epsilon_x^2 \right) =: r_f^2.
\label{eq:r_final}
\end{equation}
Note that if we were in the symmetric case, then $2\kappa \left( 1 - \epsilon_x^2 \right) = r_i^2 - 1$ from the discussion above, so that $r_f^2 \equiv r_i^2$, i.e., trajectories always remain a constant radial distance from the origin (as measured while in the bulk) in this (the symmetric) case.

\begin{figure}[!ht]
\centerline{\includegraphics[width=0.9\textwidth]{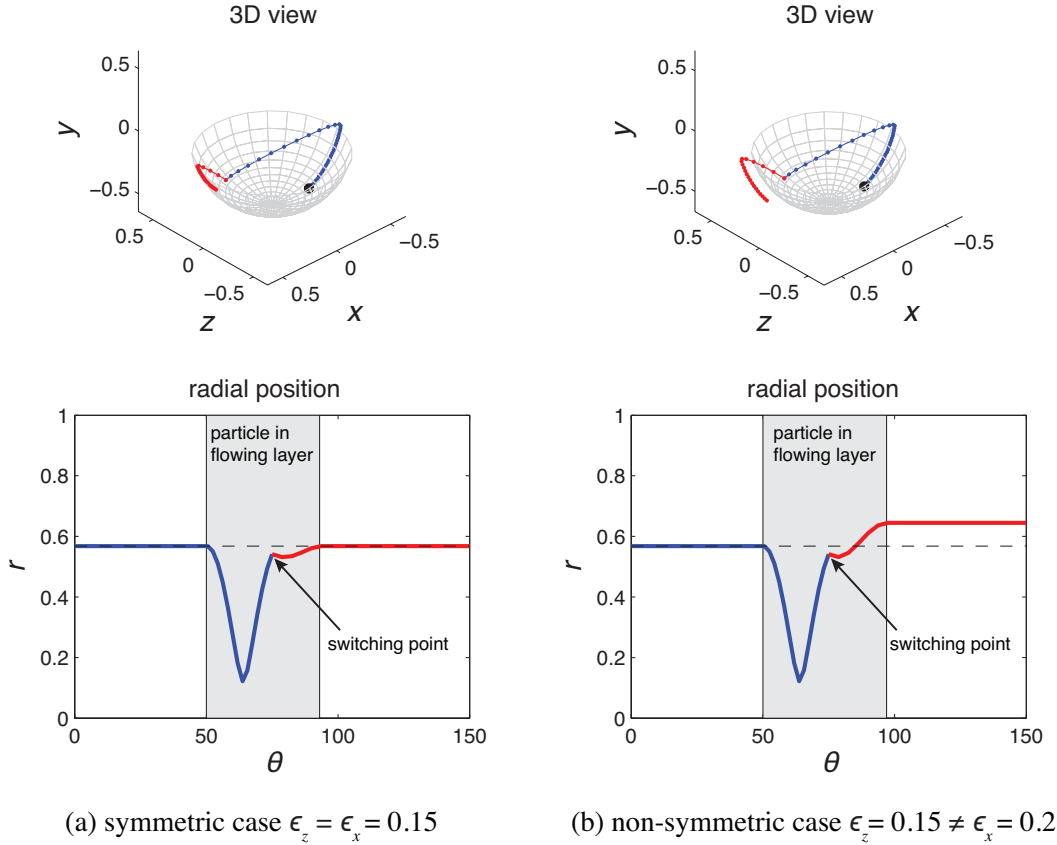}}
\caption{Illustration of the necessary condition for radial transport. The trajectory of the tracer seeded at $(x,y,z)^\top = (-0.25,-0.5,-0.1)^\top$ is visualized in 3D (top row) and through its distance from the origin (bottom row). In the 3D views, a mesh represents part of a hemisphere of the same radius as the initial particle seed location (i.e., part of a streamsurface), which is indicated by a large dot. This makes clear whether the particle remains at the same radius after one iteration of the blinking spherical tumbler protocol with $(\theta_z,\theta_x) = (5\pi/12,5\pi/12)$ in this case. Blue and red curves represent, respectively, the trajectory during the rotation about the $z$ and $x$ axes. In the bottom row, $r = \sqrt{x^2+y^2+z^2}$ is plotted against the \emph{total} angle of rotation $\theta$ (without distinguishing the current axis of rotation) in degrees. The shaded regions correspond to the particle being in the flowing layer. Dashed horizontal lines denote the initial value of $r$.}
\label{fig:radial_change}
\end{figure}

Thus, a particle cannot change the radial distance from the origin that it maintains while in the bulk (and, consequently, there is no \emph{inter}-streamsurface transport)  unless it has a ``switching point'' (i.e., the point on its trajectory at the instant when the axis of rotation is switched) that is in the flowing layer \emph{and} $\epsilon_z\ne\epsilon_x$. Figure~\ref{fig:radial_change} illustrates this result for a specific trajectory and parameter values. This argument provides an alternative proof of \cite[Theorem~3.2]{smow08}. Furthermore, this result is related to the phenomenon of streamline jumping \cite{col10} in the vanishingly thin flowing layer limit. Due to the horizontal and vertical motion of the free surface in a non-circular quasi-2D container, particles can exit the flowing layer onto a solid-body-rotation streamline that is a different distance from the center of rotation than the one they were on prior to entering the flowing layer. However, while streamline jumping occurs continuously in time and is different for all trajectories in a non-circular quasi-2D tumbler, streamsurface switching occurs at discrete time points (when the axis of rotation is switched) in the blinking 3D spherical tumbler. In \S\ref{sec:non-sym_case}, we explore the effect of such switching points in the non-symmetric case on the structure of the dynamics.

\section{Dynamical properties of the flow: period-one invariant structures}
\label{sec:period_one}

In action--action--angle coordinates, the blinking spherical tumbler can be represented as a \emph{linked twist map} \cite{sow06} by the composition of the maps
\begin{subequations}\begin{align}
\mathcal{P}_z(z,I_z,\Omega_{z,0}) &= \big(z,I_z,\Omega_{z,0} + \Omega_z(I_z) \theta_z\big),\\
\mathcal{P}_x(x,I_x,\Omega_{x,0}) &= \big(x,I_x,\Omega_{x,0} + \Omega_x(I_x) \theta_x\big),
\end{align}\end{subequations}
which describe rotation about the $z$-axis and rotation about the $x$-axis, respectively. These are termed \emph{twist} maps since the actions $z$ and $I_z$ (resp.\ $x$ and $I_x$) are preserved by the map, and only the angle, i.e., the third coordinate, changes. Further details of this transformation, including explicit expressions for the actions $I_z$ and $I_x$ and the angles $\Omega_z$ and $\Omega_x$, are given elsewhere \cite[Chapters 2 and 8]{c11}. To formally compose these maps we need a nonlinear change-of-variables transformation $\mathcal{V}$ because the maps are defined in different coordinates \cite{smow08}. Thus, formally, the linked twist map $\mathcal{Q}$ describing the flow in the blinking spherical tumbler is given by 
\begin{equation}
\mathcal{Q} = \mathcal{V}^{-1} \circ \mathcal{P}_z \circ \mathcal{V} \circ \mathcal{P}_x.
\label{eq:q}
\end{equation}
For the purposes of the discussion in the remainder of this work, it is not necessary to derive the specific form of $\mathcal{V}$ as we will not use Eq.~\eqref{eq:q} explicitly. Though an initial assessment of the structure of the map can be found in \cite{smow08}, many details were left for future study. In this section, we present a continuation of the latter examination of the period-one structures of $\mathcal{Q}$ that are also period-one structures of \emph{both} $\mathcal{P}_z$ and $\mathcal{P}_x$. Other period-one structures such as isolated periodic points of $\mathcal{Q}$ \cite{mlo07} are not considered here.

In general, period-one points are significant because they tend to give rise to the largest islands, i.e., regions of material that do not mix with the surroundings. In this 3D setting we see something quite different from the 2D setting; ``islands'' are now volumes (instead of areas) whose ``centers'' are not elliptic periodic points but normally elliptic invariant curves. These  invariant curves are surrounded by ``KAM-like tubes,'' which are barriers to mixing. In fact, the latter structures are distinct from the KAM tubes observed in ABC flows because in our system all trajectories on the ``KAM-like tubes'' in the symmetric case are closed and, hence, they cannot wind around the surface and densely fill it like certain trajectories in ABC flows do \cite{dfghms86}. The normally hyperbolic invariant curves have 2D stable and unstable manifolds \cite{w94}. The nature of trajectories near these invariant curves can be visualized by considering paths of nearby tracers \cite[\S6.4]{mlo07}. More background on the mathematical aspects of these 3D structure can be found elsewhere \cite{mw94}, though few examples of their effects on mixing and transport in realistic flows are available \cite{wig-fof}.

\subsection{Analytical expressions for period-one invariant curves}
\label{sec:location_p1_pts}

The first step in this analysis requires a general expression for the period of a trajectory starting at an arbitrary point in the flowing layer. This is easily generalized from Eqs.~\eqref{eq:x1_from_fl} and \eqref{eq:period} (based on the coordinates shown in Fig.~\ref{fig:circ_tumb} in Appendix~\ref{sec:exact_2d_sol}). To summarize, suppose without loss of generality that the trajectory starts in the flowing layer, $y_0 \in \big(-\delta_z(0,z),0\big)$ for rotation about the $z$-axis, at $x_0=0$, where we recall that $\delta_z(x,z) = \epsilon_z \sqrt{L(z)^2-x^2}$ and $L(z) = \sqrt{1-z^2}$. Then, it can be shown from the expressions in Appendix~\ref{sec:exact_2d_sol} [specifically Eqs.~\eqref{eq:period} and \eqref{eq:x1_from_fl}] that the period $T$ of the trajectory starting in the bulk at $(0,y_0)^\top$ and the coordinates $(x_1,y_1,z)$ of its intersection with the flowing layer are given by
\begin{subequations}\label{eq:T_x1_y1_3d}\begin{align}
T(y_0,z) &= {2}\tan^{-1}\left(\frac{x_1(y_0,z)}{y_1(y_0,z)}\right) + \epsilon_z{\pi},\label{eq:T_y0_z}\\
x_1(y_0,z) &= -\sqrt{L(z)^2 + \frac{2 y_0}{\epsilon_z^2}\left[\epsilon_z L(z) + \tfrac{1}{2}y_0\right]},\\ 
y_1(y_0,z) &= -\delta_z\big(x_1(y_0,z),z\big) = -\sqrt{-2y_0\left[\epsilon_z L(z) + \tfrac{1}{2}y_0\right]},\label{eq:y1_y0_z}
\end{align}\end{subequations}
Note that the square root in Eq.~\eqref{eq:y1_y0_z} is real-valued because $-\delta_z(0,z) = -\epsilon_z L(z)$ and $y_0\in\big(-\delta_z(0,z),0\big)$ by assumption.

To find those points in the flowing layer that return to their initial location after the rotation about the $z$-axis is complete, we set the period equal to the time of rotation about the current axis, i.e., $T(y_0,z) = \theta_z$ and solve for $y_0$ and $z$. To this end, setting $\varphi_z := \tfrac{1}{2}(\theta_z - \epsilon_z {\pi})$, we can re-arrange Eq.~\eqref{eq:T_y0_z} as follows:
\begin{equation}
\frac{L(z)^2}{\epsilon_z^{-2} + \tan^2\varphi_z} = -2y_0\left[\epsilon_z L(z) + \tfrac{1}{2}y_0\right].
\end{equation}
This is a quadratic equation in $y_0$ that we can easily solve:
\begin{equation}
y_{0,\pm} = -\epsilon_z L(z)\left(1 \pm \sqrt{\frac{\epsilon_z^{2}\tan^2\varphi_z}{1 + \epsilon_z^{2}\tan^2\varphi_z}} \, \right),
\label{eq:y0_p1_fl}
\end{equation}
where the negative sign must be taken to ensure that $y_0\in\big(-\delta_{z}(0,z),0\big)$, as assumed. This result is analogous to \cite[Eq.~(C2)]{smow08}.

As discussed in \S\ref{sec:radial_disp} above, trajectories in the flowing layer during rotation about the $z$-axis must satisfy the equation
\begin{equation}
\frac{1}{\epsilon_z^2}y\left[\delta_z(x,z) + \tfrac{1}{2}y\right] = \frac{1}{\epsilon_z^2}y_0\left[\delta_z(0,z) + \tfrac{1}{2}y_0\right].
\end{equation}
Hence,
\begin{equation}
y\left( \epsilon_z \sqrt{L(z)^2 - x^2} + \tfrac{1}{2}y \right) = y_0 \left[ \epsilon_z L(z) + \tfrac{1}{2}y_0 \right],
\label{eq:surf_fl_z}
\end{equation}
where the right-hand side is known thanks to Eq.~\eqref{eq:y0_p1_fl}, and $z$ is treated as a parameter. This equation implicitly defines a 2D surface in the 3D Euclidean space $\mathbb{R}^3$ representing the locus of points initially in the flowing layer that return to their initial location after a rotation by $\theta_z$. Analogously, an implicitly defined surface of period-one points for rotation about the $x$-axis can be found:
\begin{subequations}\label{eq:surf_fl_x}\begin{align}
y\left( \epsilon_x \sqrt{L(x)^2 - z^2} + \tfrac{1}{2}y \right) &= \tilde{y}_0 \left[ \epsilon_x L(x) + \tfrac{1}{2} \tilde{y}_0 \right],\\
\tilde{y}_0 &= -\epsilon_x L(x)\left(1 - \sqrt{\frac{\epsilon_x^2\tan^2\varphi_x}{1 + \epsilon_x^2 \tan^2\varphi_x}} \, \right),
\end{align}\end{subequations}
where $\varphi_x := \tfrac{1}{2}(\theta_x - \epsilon_x {\pi})$.

The intersection of the two surfaces implicitly defined in each of Eqs.~\eqref{eq:y0_p1_fl}, \eqref{eq:surf_fl_z} and Eqs.~\eqref{eq:surf_fl_x} gives rise to curves, which may be lines or hyperbolas when projected into the $(x,z)$-plane, of period-one points in the flowing layer. Generically, these curves can have different (linear) stability types, i.e., elliptic, hyperbolic or parabolic (see, e.g., \cite[\S22]{a65} or \cite[\S5.3]{o89}).

The same procedure leads to expressions for the curves of period-one points in the bulk. Gilchrist \cite{jfg03} was the first to study period-one structures in bulk of the granular flow in a bi-axial spherical tumbler, showing that the points (starting in the bulk) that return to their initial locations after a period of the flow are those that lie at the intersection of two ``bowls'' (parts of prolate spheroids) of period-one points corresponding to each rotation. This approach was further elaborated in \cite{mlo07} using a symmetry analysis applicable to more general container shapes. The equations for the bowls  (derived in Appendices A and C in \cite{smow08}) supplement the derivation of the ``caps'' above. Together, the bowl and cap corresponding to the $z$- or $x$-axis of rotation define the locus of period-one points in the bulk and the flowing layer, respectively.

Thus, we have shown that period-one points of the map $\mathcal{P}_z$ lie on the surface defined by the ``bowl'' $\Sigma_1$ and its ``cap'' $\widehat\Sigma_1$. The bowl, which is a portion of a certain prolate spheroid below the oblate spheroid that separates the flowing layer from the bulk, is given by
\begin{equation}
\Sigma_1 = \left\{ (x,y,z)^\top\in\mathbb{R}^3 \;|\: x^2 + y^2 + \mathfrak{c}_1 z^2 = \mathfrak{c}_1,\; -1 \le y \le -\delta_z(x,z) \right\},
\label{eq:sigma_1}
\end{equation}
where $\delta_z(x,z)$ is defined in Eq.~\eqref{eq:delta_z} and $\mathfrak{c}_1 = \mathfrak{c}(\epsilon_z,\theta_z)$ is a constant with 
\begin{equation}
\mathfrak{c}(\epsilon,\theta) := \frac{\epsilon^2\{1+\tan^2[\varphi(\epsilon,\theta)]\}}{1+ \epsilon^2\tan^2[\varphi(\epsilon,\theta)]},\qquad \varphi(\epsilon,\theta) := \tfrac{1}{2}(\theta - \epsilon\pi).
\label{eq:c_bowls}
\end{equation}
The latter corresponds to \cite[Eqs.~(C4) and (C5)]{smow08}. The cap is implicitly defined as
\begin{equation}
\widehat\Sigma_1 = \left\{ (x,y,z)^\top\in\mathbb{R}^3 \;|\; y\delta_z(x,z) + \tfrac{1}{2}y^2 = \delta_z^2(0,z)\left(\tfrac{1}{2}\mathfrak{d}_1^2 - \mathfrak{d}_1\right),\; -\delta_z(x,z) < y \le 0 \right\},
\end{equation}
where $\mathfrak{d}_1=\mathfrak{d}(\epsilon_z,\theta_z)$ is a constant and
\begin{equation}
\mathfrak{d}(\epsilon,\theta) := 1 - \sqrt{\frac{\epsilon^2\tan^2[\varphi(\epsilon,\theta)]}{1+ \epsilon^2\tan^2[\varphi(\epsilon,\theta)]}}\,,
\end{equation}
with $\varphi(\epsilon,\theta)$ as given in Eq.~\eqref{eq:c_bowls}. 

Similarly, period-one points for the map $\mathcal{P}_x$ lie on a surface consisting of a bowl $\Sigma_2$ with cap $\widehat\Sigma_2$. Now, the bowl is the surface
\begin{equation}
\Sigma_2 = \left\{ (x,y,z)^\top\in\mathbb{R}^3 \;|\: \mathfrak{c}_3 x^2 + y^2 + z^2 = \mathfrak{c}_3,\; -1 \le y \le -\delta_x(x,z) \right\},
\label{eq:sigma_2}
\end{equation}
where $\mathfrak{c}_3 = \mathfrak{c}(\epsilon_x,\theta_x)$.\footnote{We use the notation $\mathfrak{c}_{1}$, $\mathfrak{c}_{3}$ and $\mathfrak{d}_{1}$, $\mathfrak{d}_{3}$ for consistency with the notation in \cite{smow08}, where \emph{dimensional} variables were used leading to the necessity of introducing $\mathfrak{c}_{2}$, $\mathfrak{c}_{4}$ as well. In the present work, $\mathfrak{c}_{2} \equiv \mathfrak{c}_{1}$ and $\mathfrak{c}_{4} \equiv \mathfrak{c}_{3}$.} Its cap is implicitly defined as
\begin{equation}
\widehat\Sigma_2 = \left\{ (x,y,z)^\top\in\mathbb{R}^3 \;|\; y\delta_x(x,z) + \tfrac{1}{2}y^2 = \delta_x^2(x,0)\left(\tfrac{1}{2}\mathfrak{d}_3^2 - \mathfrak{d}_3\right),\; -\delta_x(x,z) < y \le 0 \right\},
\end{equation}
where $\delta_x(x,z)$ is given in Eq.~\eqref{eq:delta_x} and $\mathfrak{d}_3 = \mathfrak{d}(\epsilon_x,\theta_x)$.

\begin{figure}
\centering
\subfloat[$\theta_z=\theta_x=\pi$ and $\epsilon_z = \epsilon_x = 0.15$]{\includegraphics[width=\textwidth]{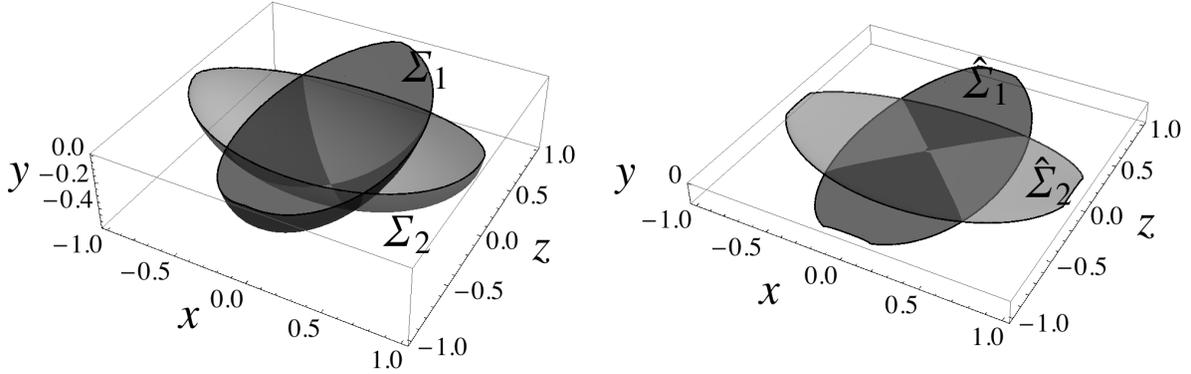}}\\
\subfloat[$\theta_z=12\pi/11$ and $\epsilon_z = 0.2$, while $\theta_x=\pi$ and $\epsilon_x = 0.1$]{\includegraphics[width=\textwidth]{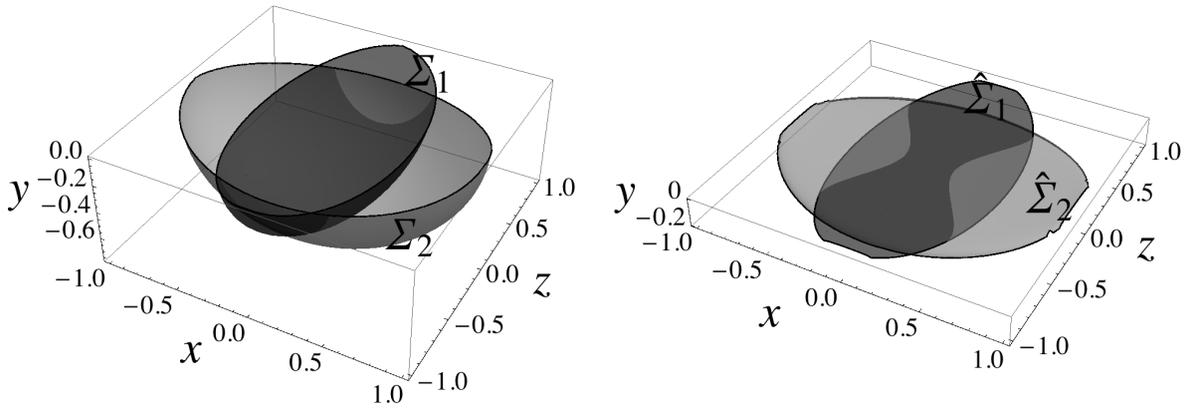}}
\caption{The intersections of the ``bowls'' $\Sigma_{1}$, $\Sigma_{2}$ (left panels) and the ``caps'' $\widehat\Sigma_{1}$, $\widehat\Sigma_{2}$ (right panels) of period-one points of the twist maps $\mathcal{P}_z$ (dark gray surfaces) and $\mathcal{P}_x$ (light gray surfaces), give rise to curves of period-one points for the linked twist map $\mathcal{Q}$.}
\label{fig:spexplain}
\end{figure}

\begin{figure}[!h]
\centering
\subfloat[$\theta_z=\theta_x = \pi$ and $\epsilon_z=\epsilon_x=0.15$]{\includegraphics[width=0.45\linewidth]{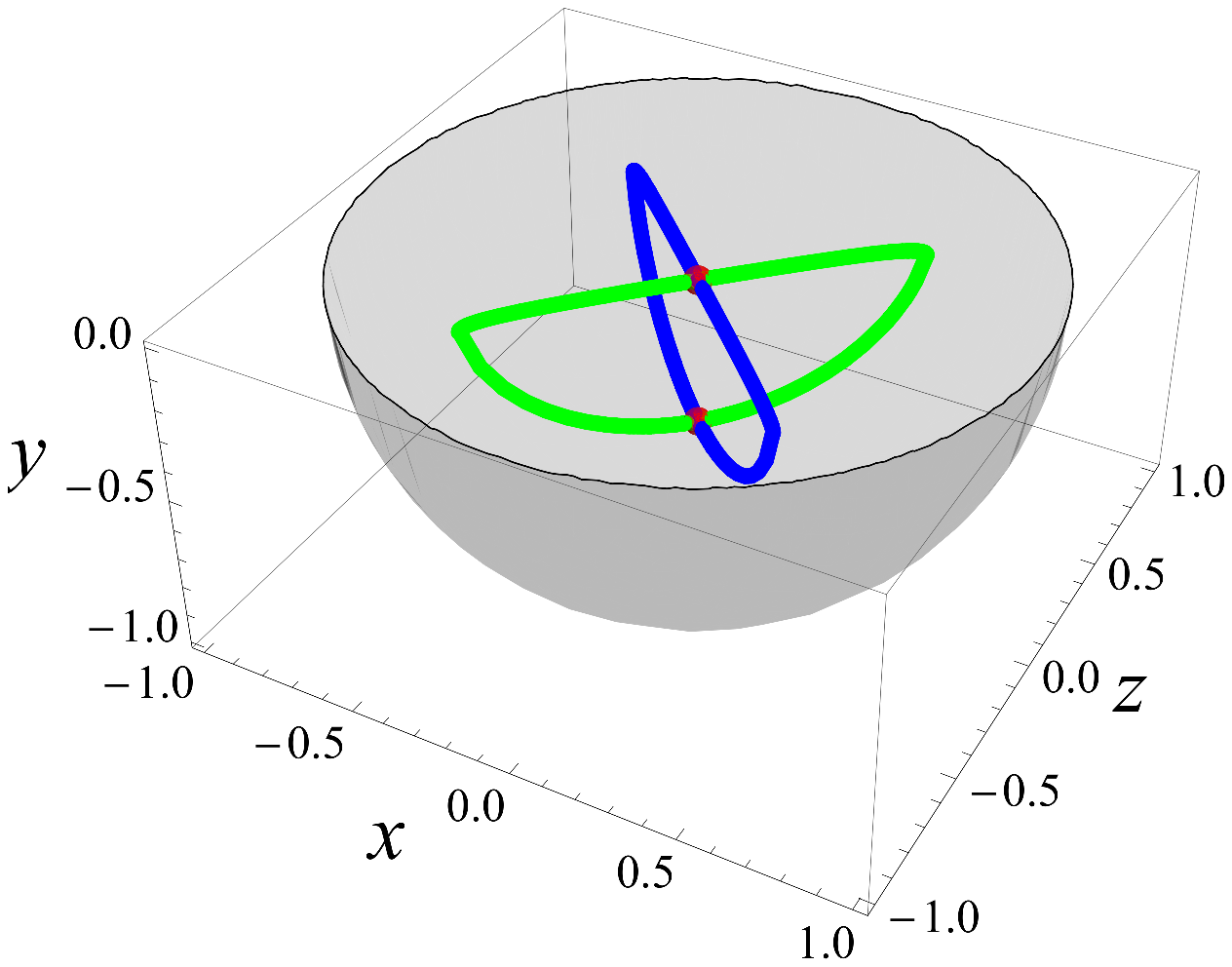}}\hfill
\subfloat[$\theta_z=12\pi/11$ and $\epsilon_z=0.15$; $\theta_x = \pi$ and $\epsilon_x=0.15$]{\includegraphics[width=0.45\linewidth]{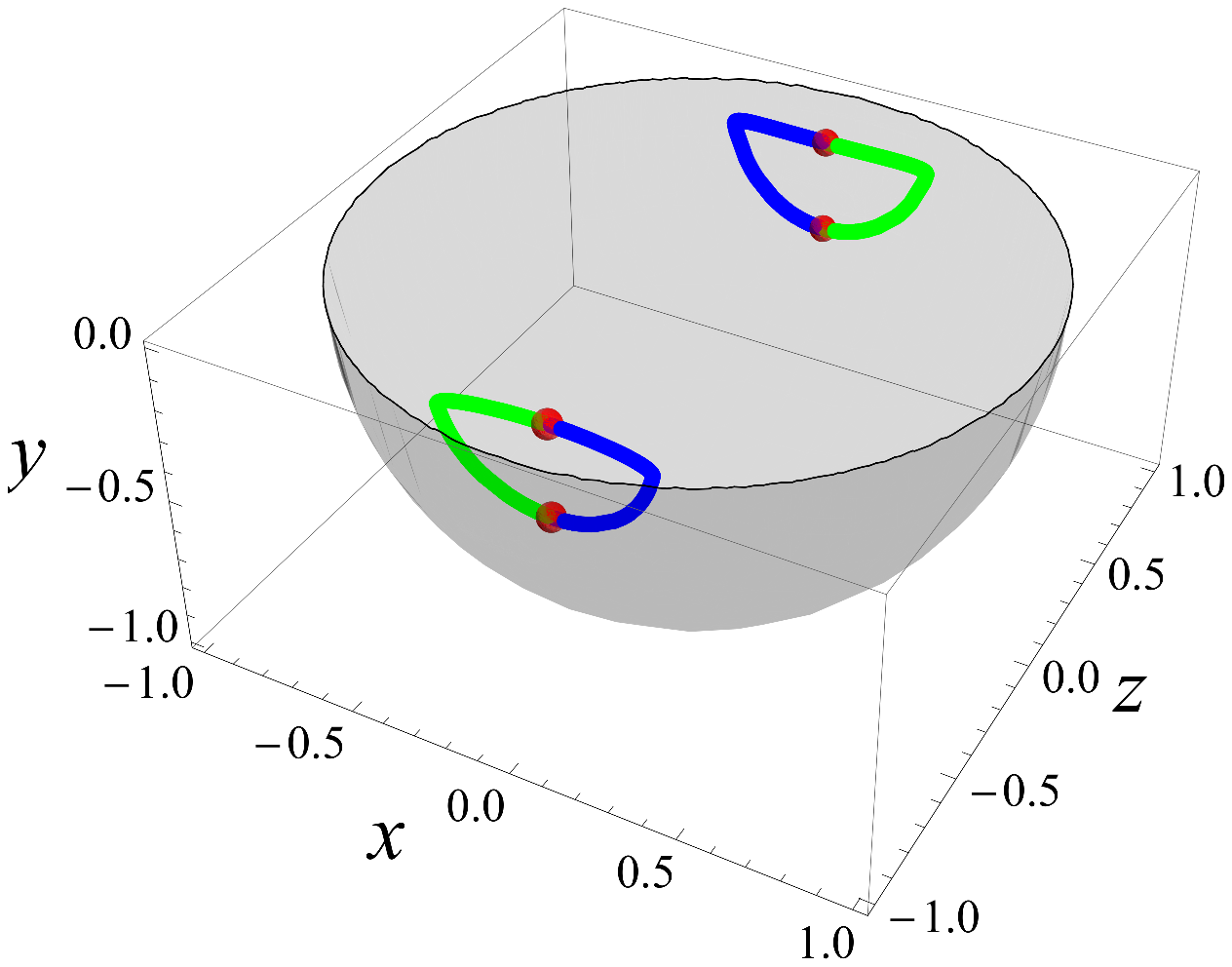}}\\
\subfloat[$\theta_z=12\pi/11$ and $\epsilon_z=0.2$; $\theta_x = \pi$ and $\epsilon_x=0.1$]{\includegraphics[width=0.45\linewidth]{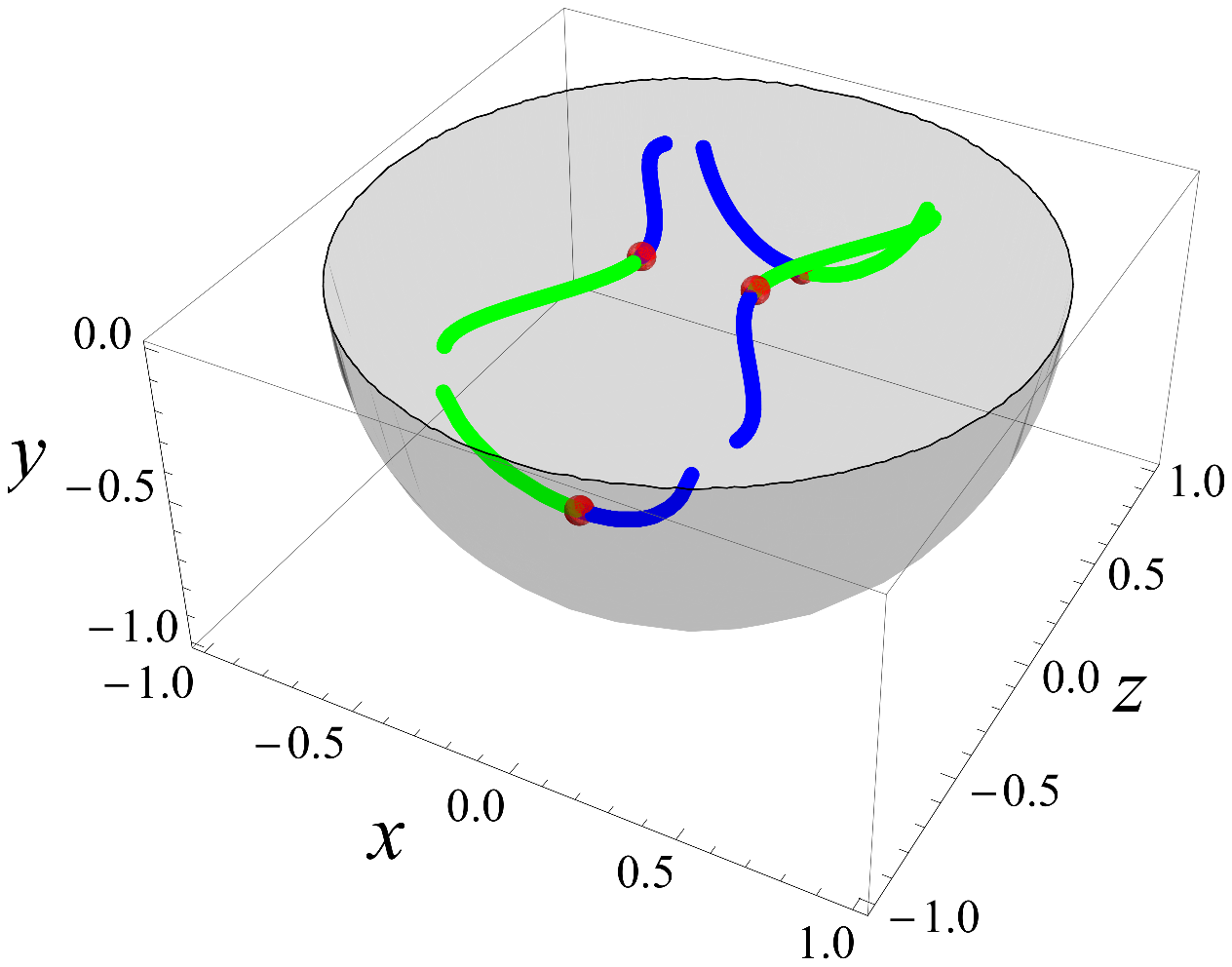}}
\caption{Plots of the period-one invariant curves $\Upsilon(x,y,z)$ of the blinking spherical tumbler flow for various choices of the system parameters. Green represents the normally hyperbolic invariant curves, while blue represents the normally elliptic ones. In the symmetric case with (a) equal angles of rotation $\theta_z = \theta_x$, the curves are interlinked, while for (b) $\theta_z\ne\theta_x$, they separate. In the non-symmetric case (c), there is a gap between the period-one curves in the flowing layer and the bulk. The gray hemisphere represents the outer wall of the spherical tumbler.}
\label{fig:p1_curv_exmpl}
\end{figure}

Schematics of the intersection of the prolate spheroids $\Sigma_1$, $\Sigma_2$ and their caps $\widehat\Sigma_1$, $\widehat\Sigma_2$ are shown separately (for clarity) in Fig.~\ref{fig:spexplain}(a) for the symmetric case with identical rotation angles of $\pi$ about each axis and in Fig.~\ref{fig:spexplain}(b) for the non-symmetric case with different rotation angles ($\theta_z = 12\pi/11$ and $\theta_x=\pi$). Using the software package {\sc Mathematica} (ver.~8.0.4), we are able to  plot implicitly defined 2D surfaces. However, it is difficult to properly resolve the endpoints of the caps, which should be sharp cusps (to match with the bowls at their point of intersection in the symmetric case, for which $\Sigma_1 \cup \widehat\Sigma_1$ and $\Sigma_2 \cup \widehat\Sigma_2$ form a single surface each) instead of the flat, straight edges seen in Fig.~\ref{fig:spexplain} near the plots' framing boxes. Nevertheless, Fig.~\ref{fig:spexplain} represents the geometry of caps and bowls accurately enough for the purposes of visualizing the shapes of these surfaces and their intersections. Note that in Fig.~\ref{fig:spexplain}(b) the union of the bowls $\Sigma_{1,2}$ and the caps $\widehat\Sigma_{1,2}$ do not have to form a single surface each, i.e., there can be a ``gap'' between them, as shown more clearly in Fig.~\ref{fig:p1_curv_exmpl}(b) below.

Period-one points of the linked twist map $\mathcal{Q}$ are thus given by $\Upsilon(x,y,z) = (\Sigma_1 \cap \Sigma_2) \cup (\widehat\Sigma_1 \cap \widehat\Sigma_2)$. In fact, the period-one points comprising $\Upsilon$ lie on \emph{curves}, since $\Upsilon$ is defined by the pairwise intersections of surfaces. The stability of these period-one points can be determined by a symmetry argument \cite[Appendix~E]{smow08}, showing that, in the bulk, $\Upsilon({xz>0})$ is a curve\footnote{The short-hand notation ``$\Upsilon({xz>0})$'' means those parts of the curve $\Upsilon(x,y,z)$ that correspond to $(x,y,z)^\top$ such that the product $xz >0$.} of unstable (hyperbolic) points (i.e., a normally hyperbolic invariant curve), while $\Upsilon({xz<0})$ is a curve of stable (elliptic) points (i.e., a normally elliptic invariant curve). It is evident that this argument applies equally well to the period-one curves in the flowing layer, though a proof is beyond the scope of the present discussion. Furthermore, the point(s) $\Upsilon(xz=0)$ at which the period-one curves change type from normally hyperbolic to normally elliptic is a (are) parabolic (degenerate) point(s). Also, note that we refer to the period-one invariant curves as \emph{normally}-elliptic and \emph{normally}-hyperbolic because the elliptic and hyperbolic behavior occurs in planes normal to the period-one curves. Meanwhile, in the third direction (along the period-one curve), there is no notion of hyperbolicity or ellipticity that applies. This reflects the quasi-2D nature of these period-one curves. Here, the distinction from isolated periodic points, which can possess hyperbolic behavior in all three coordinate directions, is made explicit.

Figure~\ref{fig:p1_curv_exmpl} shows three examples of the shape of the resulting period-one invariant curves, i.e., $\Upsilon(x,y,z)$, for particular example choices of the rotation angles $\theta_z$, $\theta_x$ and flowing layer thicknesses $\epsilon_z$, $\epsilon_x$. Their structure can be quite complicated and depends strongly on the system parameters $\theta_z$, $\theta_x$, $\epsilon_z$ and $\epsilon_x$. If $\theta_z = \theta_x$, then $\mathfrak{c}_1 = \mathfrak{c}_3$, and there is a further symmetry that leads to $\Upsilon$ consisting of two pairs of curves $\Upsilon({xz>0})$ and $\Upsilon({xz<0})$ that intersect at $x=z=0$, once with $y=-\sqrt{\mathfrak{c}_1}$ and again with $y = -\mathfrak{d}_1\epsilon_z$. In the former case, we must discard the positive root because it is outside the filled portion of the tumbler ($y\le 0$). In the latter case, we must discard the root $y=(\mathfrak{d}_1-2)\epsilon_z$ because $\mathfrak{d}_1\in[0,1]$, which means that $(\mathfrak{d}_1-2)\epsilon_z< -\epsilon_z$ is not in the flowing layer. Meanwhile, $y = -\mathfrak{d}_1\epsilon_z$ is always in the flowing layer or on its boundary [$-\epsilon_z\equiv-\delta_z(0,0) \le y \le 0$]. This is easily seen in Fig.~\ref{fig:p1_curv_exmpl}(a), which depicts the period-one curves that arise from the intersections of the bowls and caps shown in Fig.~\ref{fig:spexplain}(a). The curved lower portions of the period-one curves originate from the intersection $\Sigma_1\cap\Sigma_2$ of the bowls shown in the left panel of Fig.~\ref{fig:spexplain}(a), while the nearly straight upper portions of the period-one curves originate from the intersection $\widehat\Sigma_1 \cap \widehat\Sigma_2$ of the caps shown in the right panel of Fig.~\ref{fig:spexplain}(a). The curved lower portions of the period-one curves have been shown for the symmetric case \cite{mlo07}, but the nearly straight upper portions, as well as the period-one curves for the non-symmetric case (discussed next), have not been previously described.

Figure~\ref{fig:p1_curv_exmpl}(b), like  Fig.~\ref{fig:p1_curv_exmpl}(a), corresponds to the symmetric case but with different angles of rotation about each axis. Now, the period one curves separate (i.e., they are not interlinked) because there is no longer a common intersection point at $x=z=0$ (recall the earlier discussion). Figure~\ref{fig:p1_curv_exmpl}(c) shows an example of the non-symmetric case arising from different flowing layer thicknesses (due to unequal rotation rates) about the $z$ and $x$ axes. The parameters chosen correspond to the period-one curves at the intersections of the bowls and caps shown in Fig.~\ref{fig:spexplain}(b). Once again, since $\theta_z\ne\theta_x$, the period-one curves have separated. Furthermore, the period-one curves in the flowing layer and in the bulk are not connected because the flowing layers corresponding to each rotation do not coincide in the non-symmetric case. Understanding the structure of these period-one curves will become important below as we consider barriers to mixing and 3D generalizations of elliptic islands.


\subsection{Period-one points on the bulk portion of invariant surfaces in the symmetric case}
\label{sec:location_p1_pts_surf}

As shown in \S\ref{sec:radial_disp}, when $\epsilon_z=\epsilon_x$ trajectories are restricted to 2D invariant surfaces determined by the initial conditions. Those invariant surfaces that  intersect the period-one invariant curves constructed above can have isolated hyperbolic, elliptic or parabolic fixed points on them. Clearly, it would be useful to know more about the dynamics on the invariant surfaces, specifically surfaces on which period-one points (due to the intersection of the surface with the curves of period-one points constructed in \S\ref{sec:location_p1_pts} above) are present, their stability type, and their location. In this subsection, we answer some of these questions analytically.

From the results in \S\ref{sec:radial_disp}, it follows that trajectories in the symmetric case remain a constant distance from the origin while in the bulk. Therefore, the invariant surface on which they remain restricted has the shape of a portion of a hemispherical ``shell'' in the bulk. Consequently, since the invariant surfaces are nested and cannot intersect each other, we can uniquely parametrize all of them by the radius, $\bar{R}\in[\epsilon_z,1]$ (keeping in mind that only symmetric cases, $\epsilon_z=\epsilon_x$, are considered in this section), of their hemispherical portions in the bulk. However, fixing $\epsilon_z=\epsilon_x$ does not specify the rotation angles, since $\theta_z$ and $\theta_x$ are independent parameters. Thus, we consider two cases: $\theta_z = \theta_x$ and $\theta_z \ne \theta_x$.

\subsubsection{Case 1: $\theta_z = \theta_x$}\label{sec:case1_shells} Here, we assume not only the symmetric case ($\epsilon_z=\epsilon_x$) but also equal angles of rotation about the two axes. From the result in the previous subsection, it immediately follows that $\mathfrak{c}_1 = \mathfrak{c}_3$. Therefore, a period-one point $(x,y,z)^\top$ on the invariant surface with radius $\bar{R}$ in the bulk would have to be a solution of the simultaneous equations
\begin{subequations}\label{eqs:case1}\begin{align}
x^2 + y^2 + \mathfrak{c}_1 z^2 &= \mathfrak{c}_1,\displaybreak[3]\label{eq:case1_e1}\\
\mathfrak{c}_1 x^2 + y^2 + z^2 &= \mathfrak{c}_1,\displaybreak[3]\label{eq:case1_e2}\\
x^2 + y^2 + z^2 &= \bar{R}^2,\label{eq:case1_e3}
\end{align}\end{subequations}
where $\mathfrak{c}_1$ and $\bar{R}$ are known, but $x$, $y$ and $z$ are unknown. Note that the first two equations come from the definition of $\Sigma_1$ and $\Sigma_2$ [i.e., Eqs.~\eqref{eq:sigma_1} and \eqref{eq:sigma_2}] specialized to the case $\theta_z=\theta_x$. From Eq.~\eqref{eq:case1_e1}, $y^2 = \mathfrak{c}_1 - \mathfrak{c}_1 z^2 - x^2$, which we substitute into Eq.~\eqref{eq:case1_e2} to find that $(\mathfrak{c}_1-1)x^2 + (1-\mathfrak{c}_1)z^2 = 0$. Then, Eq.~\eqref{eq:case1_e3} becomes
\begin{equation}
\mathfrak{c}_1 + (1 - \mathfrak{c}_1)z^2 = \bar{R}^2.
\end{equation}
The possible solutions are
\begin{equation}
 x_\pm^* = \pm\sqrt{\frac{\bar{R}^2 - \mathfrak{c}_1}{1 - \mathfrak{c}_1}},\qquad y^* = -\sqrt{\mathfrak{c}_1 - \frac{1+\mathfrak{c}_1}{1-\mathfrak{c}_1}(\bar{R}^2 - \mathfrak{c}_1)},\qquad z_\pm^* = \pm\sqrt{\frac{\bar{R}^2 - \mathfrak{c}_1}{1 - \mathfrak{c}_1}}.
\label{eq:xyz_soln_case1}
\end{equation}
Thus, there are four period-one points in the bulk portion (shell) of this invariant surface:
\begin{multline} 
   (x_1,y_1,z_1)^\top = \left(x_+^*,y^*,z_+^*\right)^\top,\qquad
   (x_2,y_2,z_2)^\top = \left(x_+^*,y^*,z_-^*\right)^\top,\\
   (x_3,y_3,z_3)^\top = \left(x_-^*,y^*,z_+^*\right)^\top,\qquad
   (x_4,y_4,z_4)^\top = \left(x_-^*,y^*,z_-^*\right)^\top.
\label{eq:xyz_list_case1}
\end{multline}

Note that the radius of the shell in the bulk must be greater than the depth of the prolate spheroids corresponding to the flowing layer boundary for there to be any of these period-one points on the shell, i.e., $\bar{R} > \sqrt{\mathfrak{c}_1}$. This guarantees that $x^*_\pm$ and $z^*_\pm$ in Eq.~\eqref{eq:xyz_soln_case1} are real-valued. In addition, we must require that $\bar{R} < \sqrt{2\mathfrak{c}_1/(1+\mathfrak{c}_1)}$ for $y^*$ to be real. (It can be shown that this inequality is always satisfied once the more restrictive condition, in Eq.~\eqref{eq:case1_rb2_ineq} below, is enforced.) As discussed above, the period-one point is hyperbolic if $xz > 0$ and elliptic if $xz < 0$, hence $(x_1,y_1,z_1)^\top$ and $(x_4,y_4,z_4)^\top$ are hyperbolic period-one points, while $(x_2,y_2,z_2)^\top$ and $(x_3,y_3,z_3)^\top$ are elliptic period-one points. 

Furthermore, we cannot ignore the constraint that $y < -\max_{x,z}\big\{\delta_z(x,z),\delta_x(x,z) \big\}$ for all $(x,y,z)^\top$ on the curve because, by construction, the period-one points in Eq.~\eqref{eq:xyz_list_case1} only exist in the bulk. Hence, letting $\epsilon_z = \epsilon_x = \epsilon$, we must verify that
\begin{equation}
-\sqrt{\mathfrak{c}_1 - \frac{1+\mathfrak{c}_1}{1-\mathfrak{c}_1}(\bar{R}^2 - \mathfrak{c}_1)} < -\epsilon\sqrt{1 - \frac{\bar{R}^2 - \mathfrak{c}_1}{1 - \mathfrak{c}_1} - \frac{\bar{R}^2 - \mathfrak{c}_1}{1 - \mathfrak{c}_1}}.
\label{eq:case1_ex_cond}
\end{equation}
After some algebra, the latter inequality can be shown to be equivalent to
\begin{equation}
2\mathfrak{c}_1 - (1+\mathfrak{c}_1)\bar{R}^2 > \epsilon^2\left(1 + \mathfrak{c}_1 - 2\bar{R}^2\right).
\end{equation}
Recalling the restriction that $\bar{R}^2 > \mathfrak{c}_1$ noted after Eq.~\eqref{eq:xyz_list_case1}, we arrive at the following result: period-one points exist on the bulk shell of radius $\bar{R}$ only if
\begin{equation}
\mathfrak{c}_1 < \bar{R}^2 < \frac{\epsilon^2(1 + \mathfrak{c}_1) - 2\mathfrak{c}_1}{2\epsilon^2 - (1+\mathfrak{c}_1)},
\label{eq:case1_rb2_ineq}
\end{equation}
where we note again that this holds only for $\theta_z=\theta_x$ and any allowed $\epsilon(=\epsilon_z=\epsilon_x)$. The left inequality states that for period-one points of this type to exist on bulk portion of the 2D invariant surface, then the surface's radius in the bulk must be larger than the depth of the period-one invariant curves. The right inequality, meanwhile, indicates that beyond a certain critical value of $\bar{R}$, the bulk portion of the invariant surface fully envelops the period-one curves and, hence, does not intersect them. This, again, leads to no period-one points of this type on the bulk portion of said surface.

Here, it is important to note that Eq.~\eqref{eq:case1_rb2_ineq} holds only if $2\epsilon^2 - (1+\mathfrak{c}_1) < 0$. One way to justify this is to recall that, physically, the model does not make sense if the flowing layer is thicker than (about) half the tumbler's radius, i.e., $\epsilon \lesssim 0.5$ or $2\epsilon^2 \lesssim 0.5$. Meanwhile, $\mathfrak{c}_1\ge0\Rightarrow 1+\mathfrak{c}_1\ge 1\Rightarrow2\epsilon^2 - (1+\mathfrak{c}_1) < 0$.

Additionally, it remains to consider the degenerate case of $\bar{R} = \sqrt{\mathfrak{c}_1}$, i.e., the 2D invariant surface on which trajectories are restricted is tangent to (rather than intersects with) the period-one invariant curves. Then, Eq.~\eqref{eq:xyz_soln_case1} becomes
\begin{equation}
x_\pm^* = 0,\qquad y^* = - \sqrt{\mathfrak{c}_1},\qquad z_\pm^* = 0.
\end{equation}
From \cite[Appendix~E]{smow08}, we know this must be a parabolic periodic point. Therefore, picking a shell with radius $\bar{R} = \sqrt{\mathfrak{c}_1}$ guarantees that there is only one period-one point (of the type constructed in \S\ref{sec:location_p1_pts} above) on the surface in the bulk, and it is the parabolic one. This is the case if $\bar{R}$ is chosen so that the invariant surface is tangent to the bottom of the bulk portion of the period-one invariant curves from Fig.~\ref{fig:p1_curv_exmpl}(a), rather than being pierced by them. [Figure~\ref{fig:psecs_dsnbf}(c) below shows such a Poincar\'e section.]

\subsubsection{Case 2: $\theta_z \ne \theta_x$}\label{sec:case2_shells}
For the symmetric case in which the angles of rotation about the axes differ, $\mathfrak{c}_1\ne \mathfrak{c}_3$, Eqs.~\eqref{eqs:case1} are replaced by
\begin{subequations}\label{eqs:case2}\begin{align}
x^2 + y^2 + \mathfrak{c}_1 z^2 &= \mathfrak{c}_1,\displaybreak[3]\label{eq:case2_e1}\\
\mathfrak{c}_3 x^2 + y^2 + z^2 &= \mathfrak{c}_3,\displaybreak[3]\label{eq:case2_e2}\\
x^2 + y^2 + z^2 &= \bar{R}^2,\label{eq:case2_e3}
\end{align}\end{subequations}
where $\mathfrak{c}_1$, $\mathfrak{c}_3$ and $\bar{R}$ are known but $x$, $y$ and $z$ are (again) to be determined. Again, as made explicit by Eq.~\eqref{eq:case2_e3}, the portion of the invariant surface in the bulk is a hemispherical shell. From Eq.~\eqref{eq:case2_e1}, $y^2 = \mathfrak{c}_1 - \mathfrak{c}_1 z^2 - x^2$. Plugging this into Eq.~\eqref{eq:case2_e2}, we find that $(\mathfrak{c}_3-1)x^2 + (1-\mathfrak{c}_1)z^2 = \mathfrak{c}_3 - \mathfrak{c}_1$. Then, Eq.~\eqref{eq:case2_e3} becomes
\begin{equation}
\mathfrak{c}_1 + (1- \mathfrak{c}_1) z^2 = \bar{R}^2.
\end{equation}
Therefore, possible solutions are
\begin{equation}
 x_\pm^* = \pm\sqrt{\frac{\bar{R}^2 - \mathfrak{c}_3}{1 - \mathfrak{c}_3}},\qquad y^* = -\sqrt{\mathfrak{c}_1 - \mathfrak{c}_1 \frac{\bar{R}^2 - \mathfrak{c}_1}{1 - \mathfrak{c}_1} - \frac{\bar{R}^2 - \mathfrak{c}_3}{1 - \mathfrak{c}_3}},\qquad z_\pm^* = \pm\sqrt{\frac{\bar{R}^2 - \mathfrak{c}_1}{1 - \mathfrak{c}_1}}.
\label{eq:xyz_soln_case2}
\end{equation}
It follows that there are four period-one points in the bulk of this shell:
\begin{multline}
   (x_1,y_1,z_1)^\top = \left(x_+^*,y^*,z_+^*\right)^\top,\qquad
   (x_2,y_2,z_2)^\top = \left(x_+^*,y^*,z_-^*\right)^\top,\\
   (x_3,y_3,z_3)^\top = \left(x_-^*,y^*,z_+^*\right)^\top,\qquad
   (x_4,y_4,z_4)^\top = \left(x_-^*,y^*,z_-^*\right)^\top,
\end{multline}
where, again, $(x_1,y_1,z_1)^\top$ and $(x_4,y_4,z_4)^\top$ are hyperbolic, while $(x_2,y_2,z_2)^\top$ and $(x_3,y_3,z_3)^\top$ are elliptic. The resulting dynamics on an invariant surface that contains all four of these fixed points are qualitatively the same as for the case $\theta_z=\theta_x$ discussed in the previous subsection. If $\bar{R}$ is such that the invariant surface intersects the bulk portions of all the period-one curves in Fig.~\ref{fig:p1_curv_exmpl}(b), there will be four period-one points on its bulk portion with $\theta_z\ne\theta_x$.

To guarantee the square roots in Eq.~\eqref{eq:xyz_soln_case2} are real-valued we would need either $\bar{R} > \sqrt{\mathfrak{c}_1}$ or $\bar{R} > \sqrt{\mathfrak{c}_3}$ depending on the relative size of $\mathfrak{c}_1$ and $\mathfrak{c}_3$, which depends on the relative size of $\theta_z$, $\theta_x$, $\epsilon_z$ and $\epsilon_x$ in a {non-trivial} manner. Though it has been suggested that $\theta_z > \theta_x\Rightarrow \mathfrak{c}_1 > \mathfrak{c}_3$ \cite[Appendix~D]{smow08}, this is \emph{not} necessarily true as shown in Appendix~\ref{sec:c1_ineq_c3}.

Now, once again, we have to make sure the constraint $y < -\max_{x,z}\{\delta_z(x,z),\delta_x(x,z)\}$ is satisfied. We are considering the symmetric case, $\epsilon_z = \epsilon_x = \epsilon$, so this constraint becomes
\begin{equation}
-\sqrt{\mathfrak{c}_1 - \mathfrak{c}_1 \frac{\bar{R}^2 - \mathfrak{c}_1}{1 - \mathfrak{c}_1} - \frac{\bar{R}^2 - \mathfrak{c}_3}{1 - \mathfrak{c}_3}} < -\epsilon\sqrt{1 - \frac{\bar{R}^2 - \mathfrak{c}_3}{1 - \mathfrak{c}_3} - \frac{\bar{R}^2 - \mathfrak{c}_1}{1 - \mathfrak{c}_1}}.
\end{equation}
The latter can be manipulated into
\begin{equation}
\left\{(\epsilon^2-\mathfrak{c}_1)(1-\mathfrak{c}_3) + (\epsilon^2-1)(1-\mathfrak{c}_1) \right\}\bar{R}^2
 > (\epsilon^2 - \mathfrak{c}_1)(1-\mathfrak{c}_3) + (\epsilon^2-1)(1-\mathfrak{c}_1)\mathfrak{c}_3.
\label{eq:case2_inequality_intermediate}
\end{equation}
Now, recalling the restriction $\bar{R}^2 > \max\{\mathfrak{c}_1,\mathfrak{c}_3\}$ and noting that $(\epsilon^2-\mathfrak{c}_1)(1-\mathfrak{c}_3) + (\epsilon^2-1)(1-\mathfrak{c}_1)<0$ because $\epsilon < 1$ by assumption while $\epsilon^2 < \epsilon \le \mathfrak{c}_1,\mathfrak{c}_3 \le 1$ by definition (see also Appendix~\ref{sec:c1_ineq_c3}), we can take the last set of inequalities as strict without loss of generality. Applying these to Eq.~\eqref{eq:case2_inequality_intermediate}, we find that period-one points (of the type constructed in \S\ref{sec:location_p1_pts} above) exist on this shell of radius $\bar{R}$ only if
\begin{equation}
\max\{\mathfrak{c}_1,\mathfrak{c}_3\} < \bar{R}^2 < \frac{(\epsilon^2-\mathfrak{c}_1)(1-\mathfrak{c}_3) + (\epsilon^2-1)(1-\mathfrak{c}_1)\mathfrak{c}_3}{(\epsilon^2-\mathfrak{c}_1)(1-\mathfrak{c}_3) + (\epsilon^2-1)(1-\mathfrak{c}_1)}.
\label{eq:case2_inequality_final}
\end{equation}
Unlike the result in \S\ref{sec:case1_shells}, this one holds for $\theta_z\ne\theta_x$ and any allowed $\epsilon$. Also, note that when $\mathfrak{c}_1=\mathfrak{c}_3$ (i.e., $\theta_z=\theta_x$), Eq.~\eqref{eq:case2_inequality_final} reduces to Eq.~\eqref{eq:case1_rb2_ineq}, as required.

Finally, it remains to consider the degenerate case of equality for the first part of Eq.~\eqref{eq:case2_inequality_final}. First, let $\bar{R} = \sqrt{\mathfrak{c}_1}$. Then, Eq.~\eqref{eq:xyz_soln_case2} becomes
\begin{equation}
 x_\pm^* = \pm\sqrt{\frac{\mathfrak{c}_1 - \mathfrak{c}_3}{1 - \mathfrak{c}_3}},\qquad y^* = -\sqrt{\mathfrak{c}_1 - \frac{\mathfrak{c}_1 - \mathfrak{c}_3}{1 - \mathfrak{c}_3}},\qquad z_\pm^* = 0.
\end{equation}
There are now \emph{two} parabolic period-one points on the shell of radius $\bar{R} = \sqrt{\mathfrak{c}_1}$, instead of the single parabolic point present in the $\theta_z=\theta_x$ case. Second, if $\mathfrak{c}_1 > \mathfrak{c}_3$, then let $\bar{R} = \sqrt{\mathfrak{c}_3}$, whence
\begin{equation}
 x_\pm^* = 0,\qquad y^* = -\sqrt{\mathfrak{c}_1 - \mathfrak{c}_1 \frac{\mathfrak{c}_3 - \mathfrak{c}_1}{1 - \mathfrak{c}_1}},\qquad z_\pm^* = \pm\sqrt{\frac{\mathfrak{c}_3 - \mathfrak{c}_1}{1 - \mathfrak{c}_1}}.
\end{equation}
Again, it is clear that we have two parabolic period-one points on the shell, and their locations are mirror images about the $x$-axis. This can occur, for example, if $\bar{R}$ is chosen such that the invariant surface on which trajectories are restricted is precisely tangent to the bulk portions of the period-one curves in Fig.~\ref{fig:p1_curv_exmpl}(b). Note that this tangency would take place at the two ``lowest'' points (i.e., the $y$ coordinate is at a minimum) along the curves.

Finally, it should be noted that, in the above discussion, we neglected the structure of the corresponding period-one points on the cap, which resides in the flowing layer. Since the period-one curves in the flowing layer are given as the intersection of two implicitly defined surfaces, it is unclear whether a tractable parametrization of the analytical expressions, such as Eqs.~\eqref{eq:xyz_soln_case1} and \eqref{eq:xyz_soln_case2} above, can be derived for the location of the period-one points in the flowing layer portion of a given 2D invariant surface.

\subsection{Period-one structures on the invariant surfaces and their persistence in the weakly non-symmetric case}
\label{sec:non-sym_case}

In the preceding subsection, we showed that there can be four, two, one or zero period-one points (depending upon the value of $\bar{R}$) of the type constructed in \S\ref{sec:location_p1_pts} on the bulk part of the invariant surfaces that exists in the symmetric case. In this subsection, we illustrate these results further, and we also address the question of what happens to these surfaces and the period-one structures on them when the flow is perturbed into the non-symmetric case, which allows for radial transport as described in \S\ref{sec:radial_disp}.

\begin{figure}
\centerline{\includegraphics[width=0.95\textwidth]{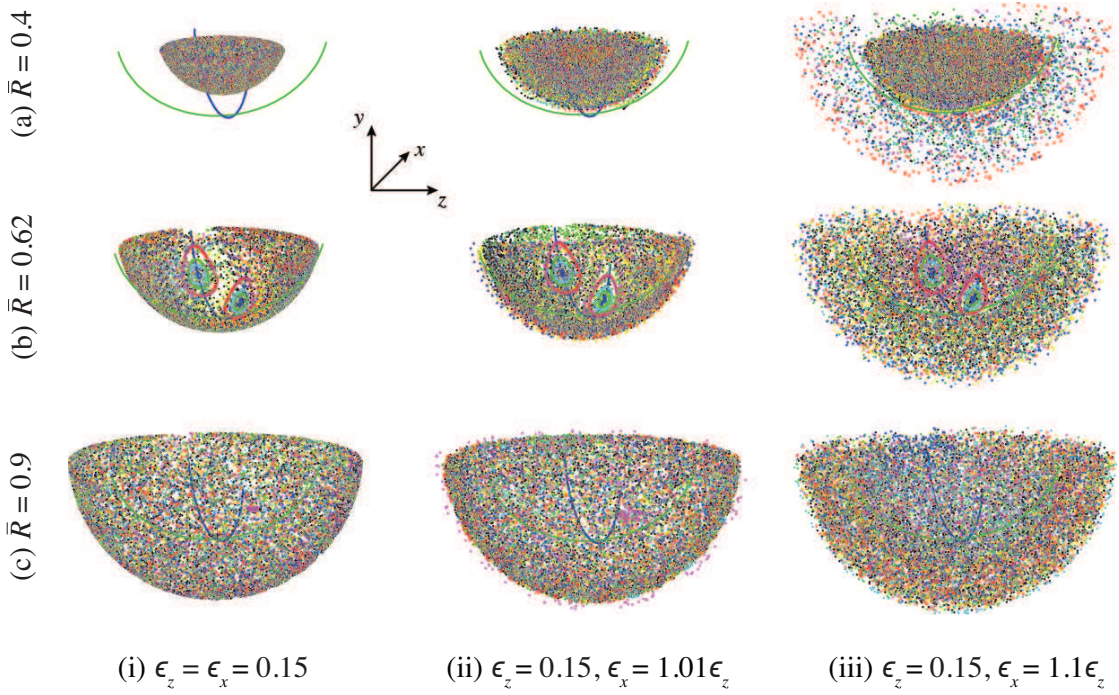}}
\caption{500-period Poincar\'e sections of tracers (initially uniformly distributed along the transect $z=-x$ in the bulk) showing the persistence of invariant surfaces in the weakly non-symmetric case as the flowing layer depth about the $x$-axis is perturbed to be (ii) $101\%$ and (iii) $110\%$ of the flowing layer depth about the $z$-axis. This allows for fully 3D transport, though for weak perturbations (ii), trajectories remain close to the invariant surfaces from the symmetric case (i). For the case of larger perturbation (iii), trajectories no longer remain close to the formerly invariant surfaces. In all plots, the angles through which the tumbler is rotated are $\theta_z = \theta_x = \pi$. For clarity, the period-one curves in the flowing layer are not shown.}
\label{fig:psecs_wns}
\end{figure}

To this end, Fig.~\ref{fig:psecs_wns} shows nine combinations of $\bar{R}$ and $(\epsilon_z,\epsilon_x)$. The symmetric case corresponds to column (i). Tracers are seeded in the bulk at a distance $\bar{R}$ from the center of the tumbler, then tracked over 500 periods of the flow. The plots shown are Poincar\'e sections, i.e., the particle locations are only plotted once per flow period.\footnote{Particle tracking is performed using the straightforward extension of the analytical solution in Appendix~\ref{sec:exact_2d_sol}. This allows for fast and efficient simulations without accumulating error from the numerical integration of the equations of motion, i.e., Eqs.~\eqref{eq:vel_z} and the corresponding ones for rotation about the $x$-axis. In all Poincar\'e sections, different colors correspond to different initial conditions of the tracer particles.} Clearly, in the symmetric case [Fig.~\ref{fig:psecs_wns}(i)], the dynamics are restricted to invariant surfaces and the dynamics on these surfaces are essentially those of a one-and-a-half-degree-of-freedom dynamical system.

Figure~\ref{fig:psecs_wns}(a)(i) shows the symmetric case and an invariant surface that is not pierced by the curves of period-one points. Thus, there are no features indicative of elliptic periodic points on this invariant surface and, as a result, there are also no visible KAM islands; however, hyperbolic period-one points not of the type constructed in \S\ref{sec:location_p1_pts} may still exist. Now, when the system is perturbed slightly into the non-symmetric case, as shown in Fig.~\ref{fig:psecs_wns}(a)(ii), trajectories leave the invariant surface on which they were constrained in the symmetric case, but they still remain ``close'' to it. This can be thought of as an adiabatic change in the system (from the symmetric case to a weakly non-symmetric one), which then allows the application of the \emph{averaging principle} \cite{a89}, at least formally. The latter motivates mathematically the question of why coherent structures such as invariant surfaces should form, and suggests that a weak perturbation cannot destroy them but can make them ``leaky'' \cite{sch06}. Thus, we observe the persistence of the invariant structures from the symmetric case, albeit in a ``fuzzier'' form. However, as the fuzzier version of the invariant surface still does not intersect the curves of period-one points, there are no visible KAM islands in the Poincar\'e section. Then, in Fig.~\ref{fig:psecs_wns}(a)(iii), the flow is perturbed well into the non-symmetric case with $\epsilon_x$ being $110\%$ of $\epsilon_z$. Now, some trajectories travel as far as the outer container walls. 

Figure~\ref{fig:psecs_wns}(b), on the other hand, presents a progression of Poincar\'e sections involving an invariant surface that is pierced by the period-one curves, showing the presence of KAM islands (evident as red closed curves at five and seven o'clock) at precisely the intersection points with the normally elliptic period-one curve [Fig.~\ref{fig:psecs_wns}(b)(i)]. Perturbing this surface into the non-symmetric case, it is evident that the KAM islands persist [Fig.~\ref{fig:psecs_wns}(b)(ii,iii)]. This is an important observation because the lack of KAM theory for the blinking spherical tumbler flow makes it impossible to guarantee the persistence of the KAM islands in the non-symmetric case \emph{a priori}. Since unmixed regions due to elliptic periodic points are present in both the symmetric and non-symmetric cases, they are not a feature unique to the essentially 2D dynamics in the symmetric case. Hence, KAM regions are generic barriers to mixing in this flow even when trajectories depart from invariant surfaces due to non-symmetric rotation rates. Meanwhile, the effect of the hyperbolic period-one points on a given invariant surface is to  create a chaotic region through the portion of the surface that is not a KAM island [see also the bottom view shown in Fig.~\ref{fig:psecs_dsnbf}(a) below]. This is more difficult to visualize in the non-symmetric case [Fig.~\ref{fig:psecs_wns}(b)(ii,iii)] because the trajectories can leave the remnants of the invariant surface.

Figure~\ref{fig:psecs_wns}(c) shows the same progression for an invariant surface that fully encloses the period-one curves. Again we observe the invariant surface becoming ``fuzzier'' as we perturb the flow into the non-symmetric case. Now, however, the changes in the nature of the Poincar\'e section are less dramatic as its size is already quite large (its radius in the bulk is $90\%$ of the container's). Trajectories leaving the surface are thus more likely to traverse the volume it encloses, which is more difficult to visualize. For the symmetric case in Fig.~\ref{fig:psecs_wns}(a), there are no period-one points on the shell because it is not pierced by the curves of period-one points. No KAM regions are evident in the non-symmetric cases [Fig.~\ref{fig:psecs_wns}(c)(ii,iii)], though it is possible for some trajectories to travel towards the center of the filled region, potentially ``spending a long time'' in or around an elliptic region. However, for the same reason that a particle away from the ``remnants'' of a KAM tube from the symmetric case can get there and spend a lot of time in this region, the possibility of 3D transport in the non-symmetric case allows for escape as well.

\begin{figure}
\centerline{\includegraphics[width=0.95\textwidth]{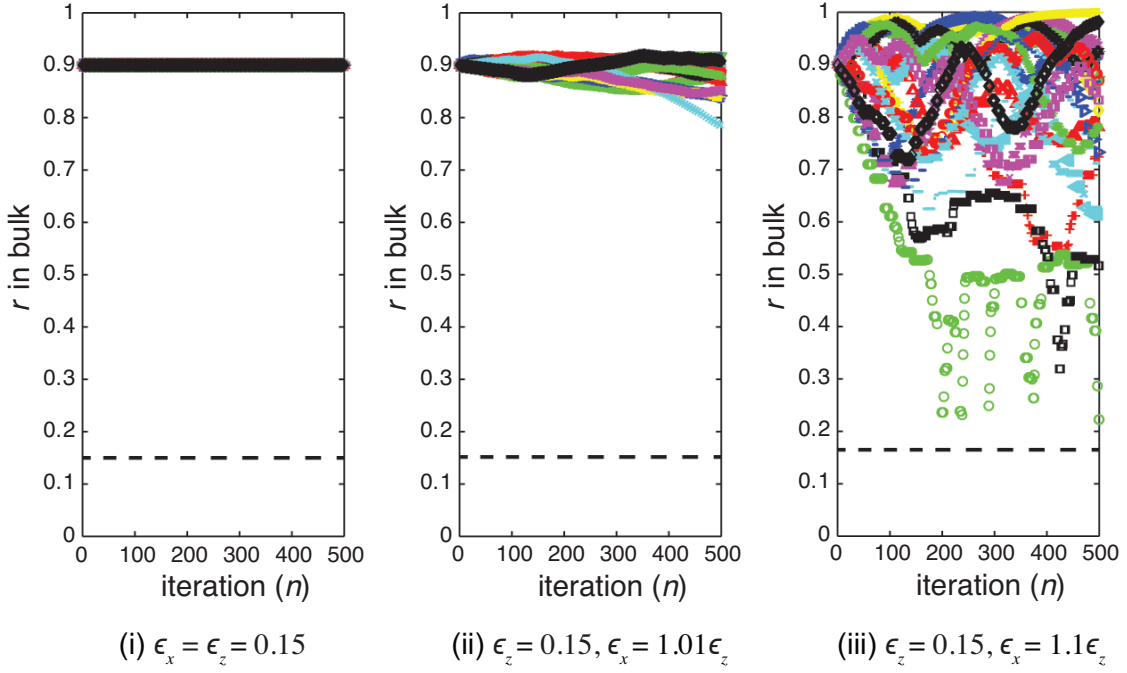}}
\caption{Radius, $r = \sqrt{x^2 + y^2 + z^2}$, in the bulk as a function of protocol iterations (flow periods) $n$ for 20 trajectories initially uniformly distributed along the transect $z=-x$ in the bulk. Different plot shapes/colors correspond to different initial conditions placed on a portion of a hemispherical shell with $\bar{R} = 0.9$ in the bulk [as in Fig.~\ref{fig:psecs_wns}(c)]. As in Fig.~\ref{fig:psecs_wns}, we observe that in the symmetric case (i) trajectories cannot change their radii in the bulk, while under non-symmetric perturbations (ii,iii) they can. On average, the rate of escape from the initial surface is clearly proportional to the strength of the perturbation as evidenced by the widening distributions of slopes. As in Fig.~\ref{fig:psecs_wns}, the angles through which the tumbler is rotated are $\theta_z = \theta_x = \pi$. The thick dashed line near the bottom of each plot denotes $r_\mathrm{min} = \max\{\epsilon_x, \epsilon_z\}$, the minimum radial position a trajectory could have while still being in the bulk.}
\label{fig:escape_rates}
\end{figure}

Finally, to clarify the escape from the invariant surfaces under non-symmetric perturbations, Fig.~\ref{fig:escape_rates} shows the radius in the bulk of a number of representative trajectories as a function of iteration count (periods of the flow). For brevity, we consider only the case of $\bar{R} = 0.9$, which corresponds to the bottom row in Fig.~\ref{fig:psecs_wns}. For the symmetric case in Fig.~\ref{fig:escape_rates}(i), the radius in the bulk of each trajectory remains constant. Meanwhile, for a weak perturbation from the symmetric case in Fig.~\ref{fig:escape_rates}(ii), the radii in the bulk of trajectories change (either increase or decrease) with the number of iterations. Finally, for a stronger perturbation in Fig.~\ref{fig:escape_rates}(iii) the radii are now quite complicated functions of the number of iterations of the protocol, showing that the 3D trajectories of tracers in the non-symmetric case explore nontrivial portions of the tumbler volume. Note that the growth of the radii is typically not smooth, often exhibiting a ``staircase-like'' pattern. As discussed in \S\ref{sec:radial_disp}, even in the nonsymmetric case, tracer particles cannot change their distance from the origin (or switch the streamsurface they are on) unless they have a switching point in the flowing layer. Thus, on average, a trajectory might wander an invariant (stream)surface for a number of flow periods before it encounters such a switching point and moves to another (stream)surface, thereby changing its radial position. Furthermore, note that the typical scale for the slopes in Fig.~\ref{fig:escape_rates}(iii) is $10^{-3}$, which corresponds to a 1\% radial change per 10 protocol iterations. Assuming this rate is sustained over the course of a trajectory, it takes on the order of 100 flow periods for a tracer particle to explore all radial positions. This is a significant gain for the modest $10\%$ increase in $\epsilon_x$, which corresponds to a $21\%$ increase in $\omega_x$ versus $\omega_z$ (for a fixed shear rate $\dot\gamma_x = \dot\gamma_z$).


\subsection{Geometry of KAM tubes in the symmetric case}
\label{sec:sym_islands}

In the previous subsection, we showed that KAM regions in the flow are persistent barriers to mixing. However, the Poincar\'e sections only revealed the local structure of KAM regions near an invariant surface, which we recall consists (in the symmetric case) of a portion of a hemisphere of radius $\bar{R}$ in the bulk with a cap through the flowing layer. In this subsection, we focus on obtaining a better understanding of the global nature of KAM tubes, specifically their shape.

Since the blinking spherical tumbler flow has \emph{curves} of normally elliptic period-one points, it is expected that KAM islands should have a 3D character. In the symmetric case, motion is restricted to 2D invariant surfaces as shown in \S\ref{sec:radial_disp}. Therefore, on each such surface there exists a KAM island. By computing the boundaries of the KAM island on each invariant surface and ``stacking'' these into a 3D volume, we can produce a 3D KAM \emph{tube}. This idea was first proposed in \cite[\S2.3.2]{jfg03} and the existence of such a 3D object was verified by numerical particle tracking using the continuum model \cite{mlo07}. Again, it is important to note that the KAM theory developed for maps in action-angle-angle variables \cite{mw94} has not been extended to linked twist maps such as the one describing the blinking spherical tumbler flow. Therefore, it may be more precise to call these objects ``KAM-like tubes''; we leave this point understood in the discussion.

\begin{figure}
\centering
\includegraphics{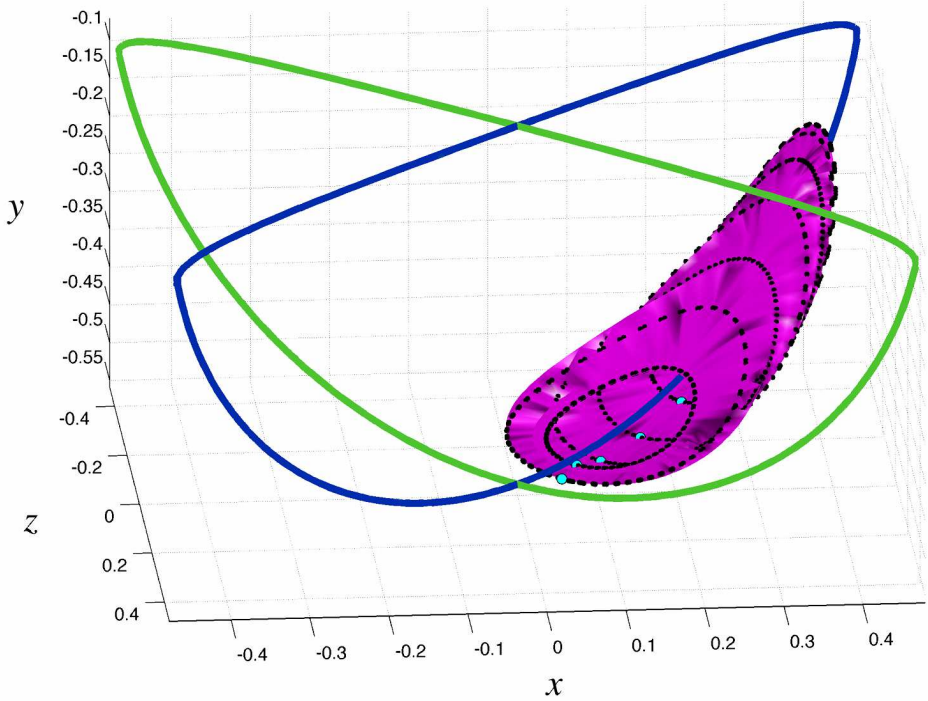}
\caption{A portion of a KAM tube in the symmetric case with $\theta_z=\theta_x = \pi$ and $\epsilon_z=\epsilon_x = 0.15$. Sample locations of tracers along the (approximate) outer boundary of several KAM islands on invariant surfaces parametrized by $\bar{R}$ are shown as black and cyan dots. The period-one normally elliptic (blue) and normally hyperbolic (green) invariant curves are shown, with the KAM tube clearly encircling the normally elliptic one. A second KAM tube exists but it is simply the mirror image across the plane $x=z$ of the one shown.}
\label{fig:kam_tubes}
\end{figure}

Figure~\ref{fig:kam_tubes} shows the result of the proposed construction. A point cloud is generated from Poincar\'e sections on a number of invariant surfaces parametrized by $\bar{R}$ (the radius of the hemi\-spherical-shell portion of these surfaces). On each surface points are seeded until their average radial distance (over 200 flow periods) from the elliptic period-one point on the surface diverges. This gives an acceptable estimate of the outer edge of the KAM region on the invariant surface. A few such islands are shown in Fig.~\ref{fig:kam_tubes}, where a cyan dot indicates the starting point for a tracer particle placed on the outer edge of the KAM region, and black dots indicate the location of this tracer after several periods of the flow. Then, the resulting point cloud of all islands' outer boundaries is triangulated \cite{mycrustopen} into the magenta surface. The resulting KAM tube tightly adheres to the normally elliptic period-one curve and thus has curvature. This is exacerbated as the KAM tube enters the flowing layer and must make a sharp turn with the normally elliptic period-one curve, making the tube difficult to mesh in this region. At the lower end of the tube, as is already evident in Fig.~\ref{fig:kam_tubes}, the tube quickly closes up near the parabolic point at the intersection of the normally elliptic and normally hyperbolic period-one curves. This is also difficult to mesh. Therefore, only a portion of the KAM tube in the bulk is shown. This suffices for the purpose of illustrating the KAM tube's intricate geometry, which is quite distinct from the KAM tori observed in fluid mixing experiments \cite{fkmo,alsm04} for which the dynamics closely resemble those of a  perturbed (non-integrable) Hamiltonian system.

The significance of the KAM tube is that, analogous to a KAM island in a quasi-2D tumbler (see, e.g., \cite{hkgmo99,mclo06,col10}), little mixing can occur within it. This is due to the fact that particle trajectories either cannot leave the KAM tube or spend a long time near it. Outside of the KAM tube, it is expected that the flow is chaotic thanks to the tangle of manifolds of the hyperbolic period-one points as will be described in \S\ref{sec:3d_manifolds}. Therefore, material away from the KAM tube can be effectively mixed, while material in or near the KAM tube cannot.

\begin{figure}
\centerline{\includegraphics[width=0.8\textwidth]{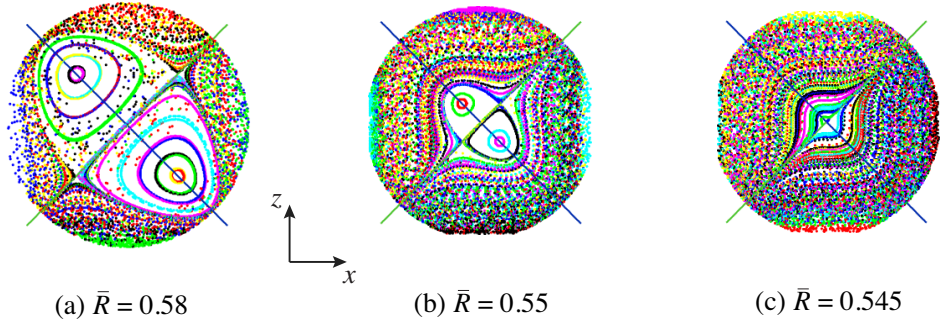}}
\caption{500-period Poincar\'e sections, as seen from below, on invariant surfaces in the symmetric case (with $\theta_z = \theta_x = \pi$ and $\epsilon_z=\epsilon_x=0.15$) as $\bar{R}$ gets progressively closer to the critical value $\sqrt{\mathfrak{c}_1}$($=0.5449\hdots$ for these parameters) from above. For $\bar{R}=\sqrt{\mathfrak{c}_1}$, there is only a single parabolic period-one point on the surface. This progression illustrates the double saddle-node bifurcation occurring on the bulk portion of the invariant surfaces. The normally hyperbolic and normally elliptic period-one invariant curves are shown in green and blue, respectively.}
\label{fig:psecs_dsnbf}
\end{figure}

It is worthwhile exploring the behavior of the KAM tube near the parabolic point further because it appears that such behavior has not been considered previously. Figure~\ref{fig:psecs_dsnbf} shows Poincar\'e sections of the blinking spherical tumbler (in the symmetric case) on different invariant surfaces parameterized once again by $\bar{R}$, as $\bar{R}\to\sqrt{\mathfrak{c}_1}^+$. Recalling the results for \S\ref{sec:case1_shells}, if $\bar{R} = \sqrt{\mathfrak{c}_1}$, the invariant surface has a single period-one point in the bulk, and it is parabolic. From Fig.~\ref{fig:psecs_dsnbf}(a) and (b) it is easy to see that the outer boundary of KAM island on these surfaces shrinks as $\bar{R}\to\sqrt{\mathfrak{c}_1}^+$, eventually becoming too small to easily identify in Fig.~\ref{fig:psecs_dsnbf}(c). This illustrates the ``pinch-off'' behavior of the KAM tube as it approaches the parabolic period-one point. It is also worthwhile pointing out that these Poincar\'e sections show a double saddle-node bifurcation \cite{kuz06,h98} as $\bar{R}\to\sqrt{\mathfrak{c}_1}^+$. That is to say, two elliptic period-one points (nodes) coalesce with two hyperbolic period-one points (saddles) in a symmetric manner as a system parameter, $\bar{R}$, approaches a critical value, $\sqrt{\mathfrak{c}_1}$.

\subsection{Minimizing the spatial extent of period-one barriers to mixing}
\label{sec:sym_opt}

The previous two subsections have shown that KAM tubes in the blinking spherical tumbler flow are generic and potentially significant barriers to mixing. Therefore, it is of interest to determine whether there are parameters such that period-one islands can be eliminated, thereby promoting more complete mixing. In this subsection, we show that it is possible to analytically determine such ``optimal'' parameters that minimize the extent of the period-one islands using the results from \S\ref{sec:location_p1_pts} and \ref{sec:location_p1_pts_surf}.

From Eqs.~\eqref{eq:case1_rb2_ineq} and \eqref{eq:case2_inequality_final}, it is clear that the upper bound on the ``depth,'' $\min_{(x,y,z)\in\Upsilon} y$, of the curves of period-one points of the type constructed in \S\ref{sec:location_p1_pts} is $\max\{\sqrt{\mathfrak{c}_1},\sqrt{\mathfrak{c}_3}\}$. Of course, from Eq.~\eqref{eq:c_bowls}, $\mathfrak{c}_1$ and $\mathfrak{c}_3$ are only functions of $\epsilon_z$, $\theta_z$ and $\epsilon_x$, $\theta_x$, respectively. To minimize the effect of the KAM tubes that surround the curves of normally elliptic period-one curves, one approach would be to make the bowls $\Sigma_1$ and $\Sigma_2$ as shallow as possible, thereby minimizing the volume in the bulk of the spherical tumbler that falls inside the KAM tubes and, hence, the volume that would remain unmixed.
 
Without loss of generality, let us restrict the discussion to the rotation about the axis that gives the deeper bowl of period-one points, noting that the bowls have unequal depths when $\theta_z\ne\theta_x$. Then, supposing that the flowing layer's depth is known, we seek to find the ``optimal'' choice of $\theta$ in the sense described above. Therefore, we seek to minimize $\mathfrak{c}(\epsilon,\theta)$ as a pure function of $\theta$, assuming that the flowing layer's depth $\epsilon$ is prescribed. From Appendix~\ref{sec:c1_ineq_c3}, we immediately find that the result of setting $\partial \mathfrak{c}/\partial \theta = 0$ and solving for $\theta$ as a function of $\epsilon$ is that $\theta_n = (2n+\epsilon)\pi$ minimizes $\mathfrak{c}$, while $\theta_m = (2m+1+\epsilon)\pi$ maximizes $\mathfrak{c}$, where $m,n\in\mathbb{Z}$. Of course, we may take $n=m=0$ without loss of generality. Then, using Eq.~\eqref{eq:c_bowls}, $\mathfrak{c}(\epsilon,\theta_m) = 1$ is the maximum, and $\mathfrak{c}(\epsilon,\theta_n) = \epsilon^2$ is the minimum. This leads to the following result: To make the extent of the period-one curves into the bulk minimal for given $\epsilon_z$ and $\epsilon_x$, choose angles of rotation $\theta_z = \epsilon_z \pi $ and $\theta_x =\epsilon_x \pi $ about each of the two axes, respectively.

\begin{figure}
\centering
\includegraphics[width=0.8\textwidth]{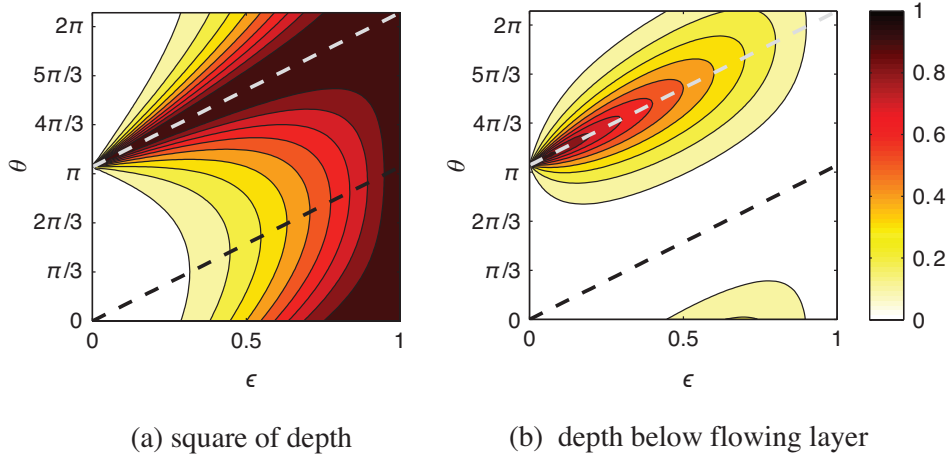}
\caption{Contour plots of the depth of a period-one bowl as a function of $\epsilon$ and $\theta$: (a) square of the depth, ${\mathfrak{c}(\epsilon,\theta)}$, with dark and light dashed lines showing the local minima and maxima, respectively; (b) depth, $\sqrt{\mathfrak{c}(\epsilon,\theta)} - \epsilon$, that the bowl extends below the flowing layer, showing that the period-one curves extend into the flowing layer only for a range of $\theta$ values.}
\label{fig:bowl_depths}
\end{figure}

Figure~\ref{fig:bowl_depths}(a) shows a contour plot of the function $\mathfrak{c}(\epsilon,\theta)$ from Eq.~\eqref{eq:c_bowls}. Clearly, the curves defined by $\theta = \epsilon \pi$ and $\theta = (1+\epsilon) \pi$ in the $(\epsilon,\theta)$-plane are the local minima and maxima, respectively, of $\mathfrak{c}$. The ``full'' 2D optimization problem is simple since $\epsilon$ is constrained to be less than the tumbler radius, which is $1$ in dimensionless units. Meanwhile, the dependence of $\mathfrak{c}$ on $\theta$ is periodic with period $2\pi$. Therefore, knowing the curves of local minima and maxima in the $(\epsilon, \theta)$-plane and the boundaries of the compact set on which the smooth function $\mathfrak{c}$ is defined, say $[0,2\pi]\times[0,1]$, we readily know its global maximum and minimum. From a practical standpoint, the periodicity in $\theta$ simply reflects the fact that, under the assumption of no transport in the direction transverse to the flow, there is no improvement in the transport (or mixing) properties of the flow for rotation angles greater than $2\pi$ radians. In addition, in Fig.~\ref{fig:bowl_depths}(b), the depth of the bowls below the flowing layer [i.e., $\sqrt{\mathfrak{c}(\epsilon,\theta)} - \epsilon$]  is plotted as well. This shows that there is only a limited range of $\theta$ values for any given $\epsilon$ for which the bowls (and consequently curves) of period-one points extend significantly into the bulk.

To summarize: if we assume \emph{a priori} that period-one structures of the type constructed in \S\ref{sec:location_p1_pts} (such as KAM tubes) are the worst barriers to mixing, then the results of this section give the optimal angles for mixing protocols to minimize such KAM tubes in the bulk. Thus, as long as $\epsilon$ and $\theta$ are chosen so that their values fall within the white region in Fig.~\ref{fig:bowl_depths}(b), the curves of period-one points (specifically, the normally elliptic ones) do not extend significantly into the bulk and do not adversely affect the quality of mixing. Furthermore, due to the periodicity in $\theta$, there is no added value (with respect to mixing) in rotating by a large angle (beyond $2\pi$) in the absence of diffusion and/or an axial flow.

\section{Structures associated with hyperbolic period-one points in the symmetric case}
\label{sec:3d_manifolds}

In this section, we present selected numerical results on the structure of manifolds of hyperbolic period-one points in the symmetric case of the model. While in the previous section we considered the most significant barriers to mixing, i.e., KAM tubes, the structure of manifolds also provides  insight into transport in the blinking spherical tumbler flow. We show that the paradigm of the \emph{invariant manifold template} for chaotic transport \cite{blw94} applies. More strikingly, however, we show that not all manifolds form chaotic tangles, leaving some further barriers to mixing in the forms of separatrices. Details of the computational method used are given in Appendix~\ref{app:manifolds}.

\subsection{Geometry of manifolds of hyperbolic period-one points}\label{sec:manifolds_on_shells}

As elaborated in \S\ref{sec:sym_islands} and illustrated in Fig.~\ref{fig:psecs_dsnbf}, the bulk portion of each invariant surface, to which trajectories are restricted in the symmetric case, can have at most two hyperbolic and two elliptic period-one points of the type constructed in \S\ref{sec:location_p1_pts}. The resulting Poincar\'e sections show dynamics (on the invariant surface) that are similar to a perturbed (non-integrable) double saddle-node Hamiltonian system. Therefore, as a stepping stone to understanding the chaotic dynamics of the LTM corresponding to the blinking spherical tumbler flow, it is helpful to consider two classical examples of a non-integrable perturbation to an integrable Hamiltonian double saddle-node system: a simple model of 2D cellular Rayleigh--B\'enard convection \cite{sg88} and the flow generated by an oscillating pair of planar vortices \cite{rklw90}.

\begin{figure}
\begin{center}
\subfloat[streamlines, elliptic (green) and hyperbolic (red) period-one points, and separatrices (blue) of a steady cellular flow]{\includegraphics[width=0.45\textwidth]{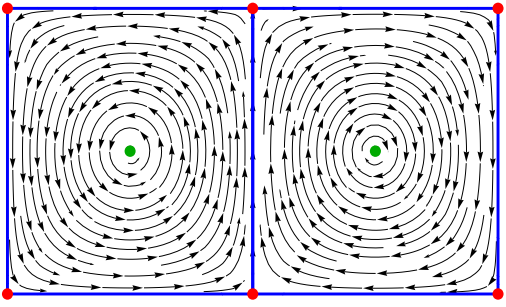}}\hfill
\subfloat[heteroclinic tangle of stable (black) and unstable (red) manifolds formed by the break-up of a separatrix (blue) under a weak perturbation]{\includegraphics[width=0.45\textwidth]{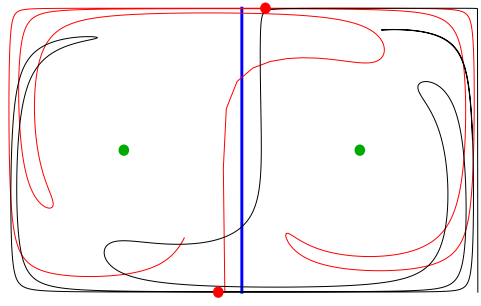}}
\end{center}
\caption{Forcing a time-periodic translation of the hyperbolic points in spatially periodic convection rolls (a), which represent a double saddle-node Hamiltonian dynamical system, leads to the break up of separatrices and the formation of a tangle of transversely intersecting manifolds (b). Lobe dynamics (mapping of area between stable and unstable manifold intersections) allows for transport in and out of each vortex (entrainment and detrainment). All vertical separatrices delineating adjacent vortices break up into tangles under perturbation.}
\label{fig:sepr_perturb}
\end{figure}

\begin{figure}
\centering
\subfloat[streamlines, elliptic (green) and hyperbolic (red) period-one points, and separatrices (blue) of a steadily translating vortex pair]{\includegraphics[width=0.45\textwidth]{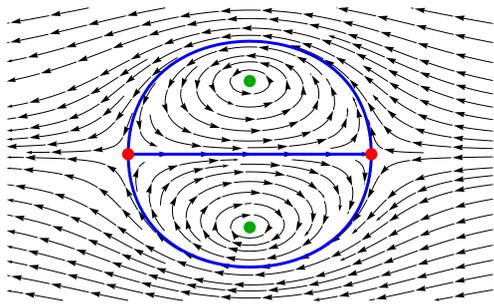}}\hfill
\subfloat[transverse intersection of stable (black) and unstable (red) manifolds of hyperbolic periodic points in the perturbed vortical flow]{\includegraphics[width=0.45\textwidth]{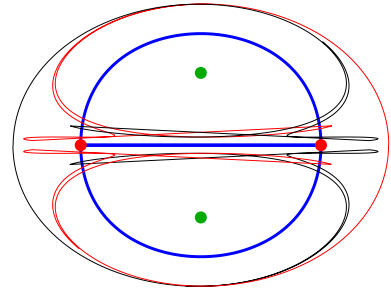}}
\caption{Subjecting a steadily translating vortex pair (a), which is a double saddle-node Hamiltonian dynamical system, using a time-periodic strain-rate field leads to the break up of separatrices and the formation of a manifold tangle (b). Lobes dynamics (mapping of area between stable and unstable manifold intersections) allows for transport in and out of the vortex (entrainment and detrainment). The horizontal separatrix connecting the two hyperbolic points persists under perturbation, preventing transport between the upper and lower half-planes.}
\label{fig:lobe_dynamics}
\end{figure}

When an incompressible thermally conducting fluid between two plates is heated from below, a horizontally periodic array of counter-rotating vortices, known as Rayleigh--B\'enard convection rolls, two cells of which are illustrated in Fig.~\ref{fig:sepr_perturb}(a), forms \cite{sg88}. A heteroclinic trajectory connecting the hyperbolic fixed points on the top and bottom walls acts as a separatrix, preventing  transport between the two convection cells. (Two additional heteroclinic trajectories form the top and bottom boundaries of each convection cell.) Each convection cell has an elliptic fixed point at its center. Beyond a critical temperature difference between the top and bottom plates, an instability causes the flow to become time-dependent but the velocity field remains essentially a perturbation of the steady one. In particular, the time-dependent perturbation causes the convection cells to oscillate horizontally (breaking the left-right reflection symmetry), the separatrix through the middle of Fig.~\ref{fig:sepr_perturb}(a) breaks up into a tangle of stable and unstable manifolds, as shown in Fig.~\ref{fig:sepr_perturb}(b), which exhibits the typical lobe dynamics transport behavior \cite{cw91}. The two separatrices on the top and bottom of each convection cell (and coinciding with the top and bottom walls) remain, bounding the dynamics in the vertical direction. 

Another, topologically equivalent, double saddle-node Hamiltonian system is a pair of steadily translating vortices in an inviscid incompressible fluid, a 2D version of Hill's spherical vortex \cite{h94}. The resulting dynamical system's phase portrait is illustrated in Fig.~\ref{fig:lobe_dynamics}(a). Similarly to the previous example, one separatrix prevents transport between the two vortices, while two other separatrices prevent transport from each vortex into the ambient fluid. Upon introducing a time-periodic strain-rate field that causes the vortices to oscillate, the separatrices demarcating the vortices from the uniform outer flow break up into chaotic tangles \cite{rklw90} as shown Fig.~\ref{fig:lobe_dynamics}(b). However, unlike the previous example, the separatrix the middle that prevents transport between the vortices persists (the system maintains its up-down reflection symmetry), so that transport is only possible between the ambient (unbounded) fluid and each vortex.

Now, we have two possible perturbed-integrable-Hamiltonian-system templates against which to compare the manifold structure of the blinking spherical tumbler flow. A first clue about the structure of the manifold tangles in the blinking spherical tumbler is the lack of a chaotic region between the KAM islands on invariant surfaces with $\bar{R}$ close to $0.54\hdots$ in Fig.~\ref{fig:psecs_dsnbf}(a). It appears that a separatrix persists, connecting the two hyperbolic points in the bulk. Indeed, similarly to the pair of 2D vortices in Fig.~\ref{fig:lobe_dynamics}(a), the blinking spherical tumbler with $\theta_z=\theta_x$ possesses a symmetry of its invariant curves and, hence, its nodes (resp.\ saddles) about $x=z$ (resp.\ $x=-z$). In contrast, Poincar\'e sections on the invariant surfaces on which the period-one points are farther apart [Fig.~\ref{fig:psecs_wns}(b)(i)] appear to have a chaotic region between the KAM islands.

\begin{figure}[!ht]
\centering
\subfloat[manifolds of forward ($x,z>0$)$\qquad\qquad$ hyperbolic point in the flowing layer]{\includegraphics[width=0.5\textwidth]{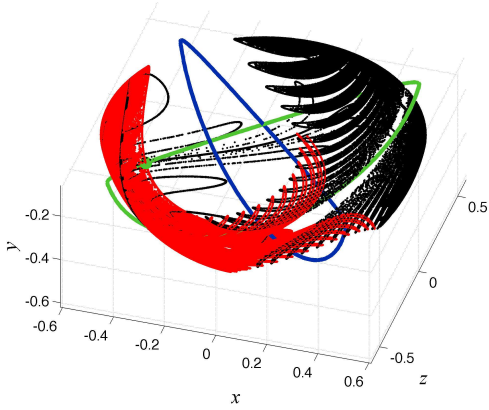}}\hfill
\subfloat[manifolds of rear ($x,z<0$)$\qquad\qquad\qquad$ hyperbolic point in the flowing layer]{\includegraphics[width=0.5\textwidth]{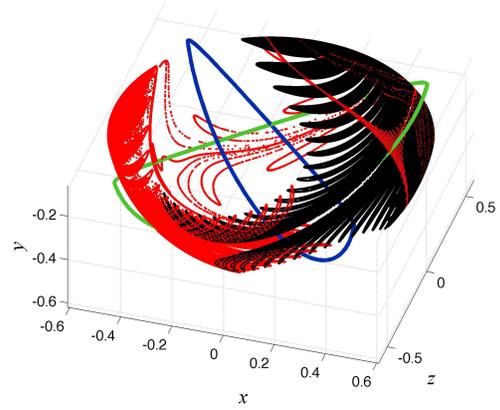}}\\
\subfloat[manifolds of forward ($x,z>0$)$\qquad\qquad\;$ hyperbolic point in the bulk]{\includegraphics[width=0.5\textwidth]{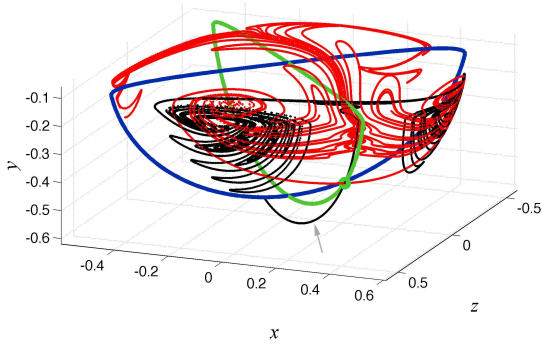}}\hfill
\subfloat[manifolds of rear ($x,z<0$)$\qquad\qquad\qquad$ hyperbolic point in the bulk]{\includegraphics[width=0.5\textwidth]{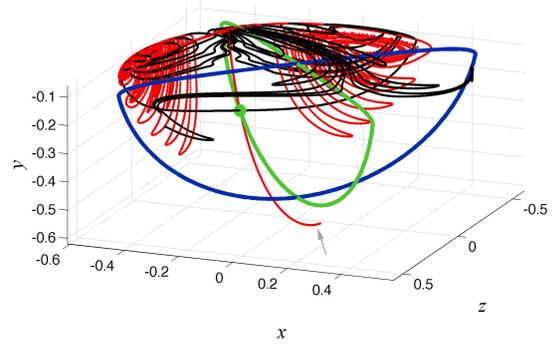}}
\caption{Stable (black) and unstable (red) manifolds of a pair of mirror hyperbolic periodic points in (a,b) the flowing layer and (c,d) the bulk for the invariant surface with $\bar{R}=0.62$. For all panels, $\theta_z=\theta_x = \pi$, $\epsilon_z =\epsilon_x=0.15$. The manifolds are traced out over 27 periods (a,b) or 40 periods (c,d) of the flow. Note the apparent coincidence between a branch of the stable manifold of the forward hyperbolic point in the bulk (c) and a branch of the unstable manifold of the rear hyperbolic point in the bulk (d), indicated by arrows.}
\label{fig:3d_manifolds}
\end{figure}

To reconcile these observations, Fig.~\ref{fig:3d_manifolds} shows traces of the stable and unstable manifolds of the four hyperbolic period-one points on the invariant surface (in the symmetric case) that we studied through Poincar\'e sections in Fig.~\ref{fig:psecs_wns}(b)(i). In Fig.~\ref{fig:3d_manifolds}(a,b), for the two hyperbolic period-one points in the flowing layer, we see the typical meshwork of infinitely nested intersections of stable and unstable manifolds. What is surprising, however, is that this meshwork extends well into the bulk portion of the invariant surface. This suggests that if there is a ``chaotic transport channel'' between the KAM islands, then it is created by a lobe dynamics transport mechanism mediated by the hyperbolic points \emph{in the flowing layer}. Furthermore, as Fig.~\ref{fig:3d_manifolds}(c,d) shows, the structure of the manifolds of the hyperbolic period-one points in the bulk is markedly different from that in the flowing layer [note the orientation of the $z$-axis is reversed between panels (a,b) and (c,d)]. For example, the meshwork of intersections does not extend through the bottom of the bulk part of the surface, rather it remains localized near the tumbler walls close to the flowing layer. Here it is important to emphasize that the tangles of manifolds will eventually fill the entire invariant surface if we trace them over enough periods because the manifolds themselves have infinite length. However, studying a finite length of manifold is enlightening because it gives an organizing framework \cite{blw94} for the observed chaotic dynamics for trajectories passing near the hyperbolic points.

Thus, while the chaotic tangle produced by the manifolds of the hyperbolic period-one points in the flowing layer is similar to the perturbed-integrable-Hamiltonian-system template of 2D cellular Rayleigh--B\'enard convection [Fig.~\ref{fig:lobe_dynamics}(a)], the tangle produced by the manifolds of the hyperbolic period-one points in the bulk retains a symmetry (for $\theta_z=\theta_x$) that makes it more similar to the one of the oscillating vortices [Fig.~\ref{fig:lobe_dynamics}(b)]. Though both the bulk and the flowing layer contain a double saddle-node structure, the dynamics is distinct. This also serves to highlight the linked twist map (LTM) structure, which is not a simple perturbation of an integrable Hamiltonian system, of the blinking spherical tumbler flow. (Here, we should stress that the dynamics on invariant surfaces in the symmetric case is still those of a \emph{volume-preserving} map.) The important point is that the separation of dynamics between the bulk and the flowing layer, with those in the bulk being simply solid body rotation, leads to a transport mechanism mediated primarily by the tangle of manifolds of hyperbolic points in the flowing layer.

\subsection{Heteroclinic and homoclinic connections in the partially chaotic flow}\label{sec:hc_connections}

\begin{figure}
\centering
\subfloat[$\theta_z=\theta_x=\pi$]{\includegraphics[width=0.49\textwidth]{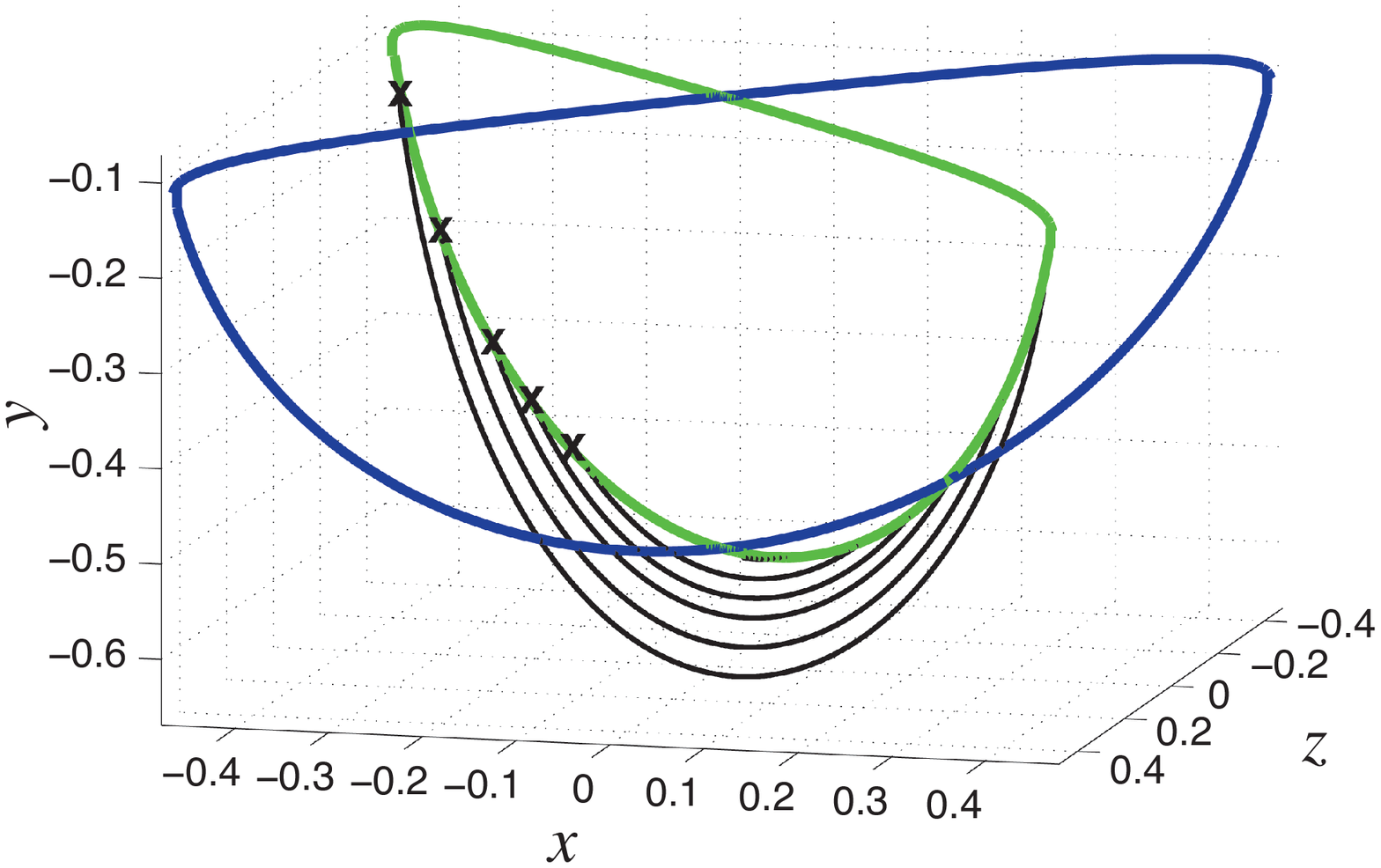}}\hfill
\subfloat[$\theta_z=\pi\ne\theta_x=19\pi/20$]{\includegraphics[width=0.49\textwidth]{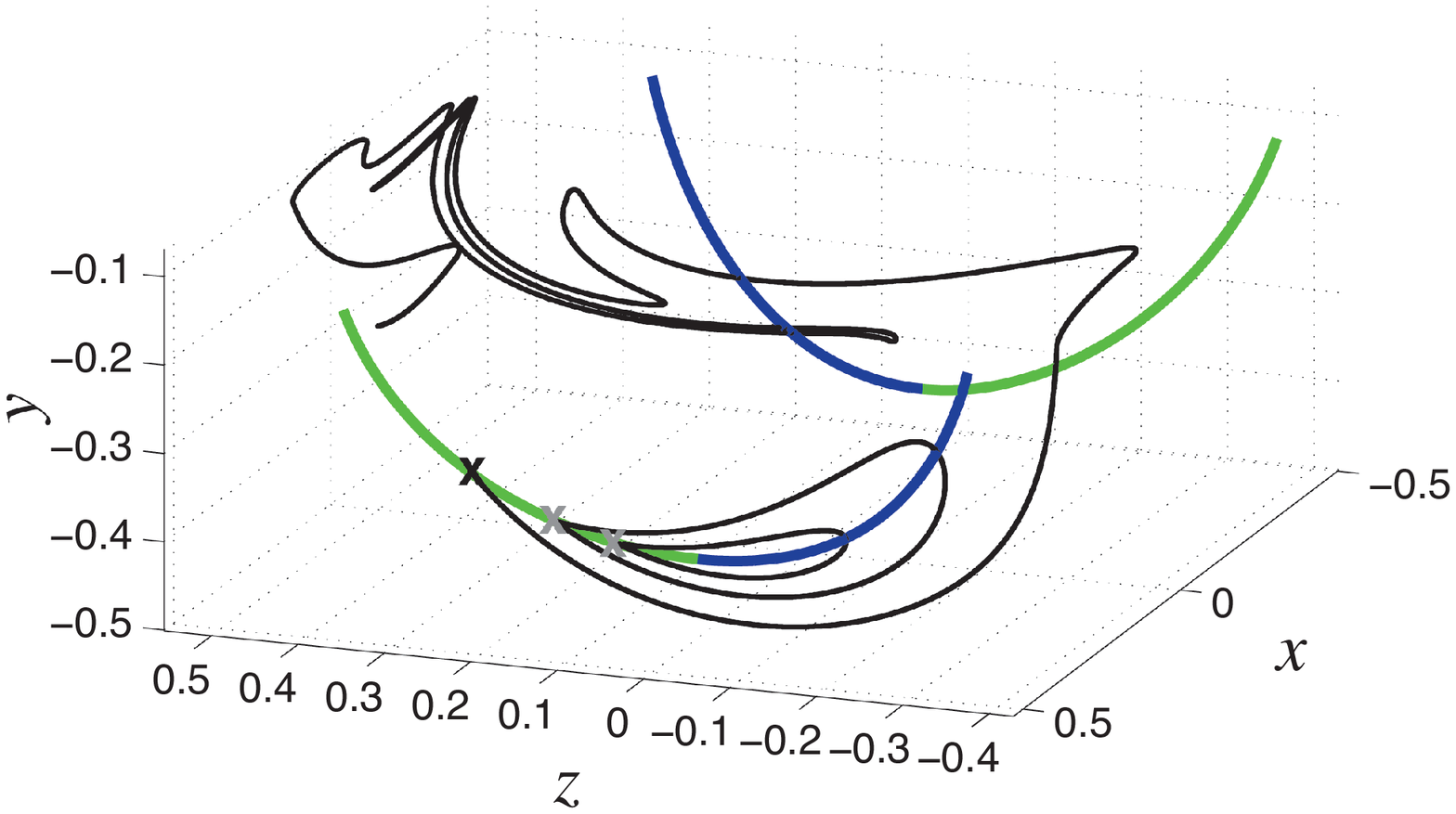}}
\caption{Heteroclinic and homoclinic trajectories (a few shown starting from the hyperbolic points denoted by $\times$s) persist in the bulk of the chaotic blinking spherical tumbler flow, while those in the flowing layer breakup into a tangle of manifolds. The extent of the period-one curves in the flowing layer is not shown in panel (b) out of convenience. In all panels, we are under the conditions for the symmetric case of the model with $\epsilon_z=\epsilon_x = 0.15$.}
\label{fig:barriers}
\end{figure}

The apparent coincidence of a branch of the stable manifold of the forward hyperbolic fixed point in the bulk from Fig.~\ref{fig:3d_manifolds}(c) with the branch of the unstable manifold of the rear hyperbolic fixed point from Fig.~\ref{fig:3d_manifolds}(d) (both indicated by arrows) suggests the persistence of a heteroclinic connection in the chaotic blinking spherical tumbler. To elucidate this observation, in Fig.~\ref{fig:barriers}(a), we trace the stable manifold of the forward hyperbolic fixed (in the Poincar\'e map) point in the bulk until it is within some prescribed distance from the rear hyperbolic fixed point. This allows us to confirm, up to numerical precision, the existence of the heteroclinic connection. Indeed, since we can ``stack'' trajectories across invariant surfaces in this (the symmetric) case, the heteroclinic trajectories form a 2D structure. Note, however, that the heteroclinic trajectories do not lie in the $x=z$ plane of symmetry to which the normally hyperbolic period-one curves are restricted \cite[Appendix~A]{mlo07}. From Fig.~\ref{fig:barriers}(a), it is evident that the shape of the resulting 2D structure is not flat. Rather, it looks more like a ``scoop'' as can be seen from the further forward protrusion of each consecutively lower black curve. Since manifolds are invariant surfaces that material cannot cross, there cannot be trajectories in the Poincar\'e section that pass through these 1D barrier to mixing across the 2D invariant surfaces. Hence, the only way for an orbit of the Poincar\'e map to switch between the two hemispheres of the invariant surface (as divided by a heteroclinic connection) is by passing through the flowing layer.

In \S\ref{sec:location_p1_pts} (specifically, Fig.~\ref{fig:p1_curv_exmpl}), it was shown that when $\theta_z\ne\theta_x$, the period-one curves separate, no longer crossing at some point with $x=z=0$. To see the effect of such change in the period-one curves' topology, Fig.~\ref{fig:barriers}(b) shows some example traces of stable manifolds for $\theta_z\ne\theta_x$. The heteroclinic connections from Fig.~\ref{fig:barriers}(a) are now \emph{homoclinic} connections. That is, the stable manifold does not connect across the tumbler to the opposite hyperbolic fixed point on the invariant surface. Rather, it loops back to the original hyperbolic point, which is denoted by an `$\times$.' This creates local barriers to mixing (separatrices) that are not as pronounced as the hemisphere-delineating ones in Fig.~\ref{fig:barriers}(a). Furthermore, these homoclinic connections do not exist for all period-one points along the normally hyperbolic invariant curve. At some critical point along the curve, this separatrix ceases to exist, and the stable manifold appears to become part of a homoclinic tangle, as is the case for the trace from the leftmost `$\times$' in Fig.~\ref{fig:barriers}(b). It would be of interest to develop theoretical results on the location of the critical point (at which the homoclinic trajectory ceases to exist)  as a function of the system parameters, specifically the angles of rotation $\theta_z$ and $\theta_x$.

\section{Conclusions and outlook}

Canonical examples of 3D chaotic dynamical systems are few, with the more realistic ones being analytically intractable. In this respect, the blinking spherical tumbler offers an interesting duality: it is simple enough to be analytically tractable, yet it retains enough of the kinematics of granular flows to be able to reproduce laboratory experiments \cite{jfg03,go03,mlo07}.

In the present work, an important finding is the necessary condition, formulated in terms of the system's parameters, for radial (i.e., fully 3D) transport in the tumbler. This allows for the identification of a ``symmetric'' versus a ``non-symmetric'' mode of operation of the mixer. Similarly, the derivation of analytical expressions for a class of period-one invariant curves, allows us to understand their structure and dependence on the parameters, which is difficult to do when the periodic points must be found by an extensive numerical search of the domain (as in, e.g., \cite{agpvm99,sch04,psc10}). Furthermore, since the blinking spherical tumbler is a linked twist map (specifically, a composition of two action-action-angle maps), its mathematical structure is not precisely that of a perturbed integrable Hamiltonian system such as the ABC flows first considered as prototypical 3D chaotic systems. Moreover, the computational analysis of the structure of the manifolds of hyperbolic period-one points (in the symmetric case) further highlights the separation of dynamics between the flowing layer and the bulk, each being loosely analogous to a different perturbed-integrable-Hamiltonian-system template. Nevertheless, the analogy to non-integrable Hamiltonian dynamics is helpful in understanding the bifurcations of period-one structures on the invariant surfaces (in the symmetric case) of the flow.

In the future, it would also be interesting to extend the results from \S\ref{sec:3d_manifolds}, specifically to visualize the 2D surfaces of stable and unstable manifolds of the normally hyperbolic invariant curves. When the dynamics are restricted to \emph{planes} in $\mathbb{R}^3$, the ``stacking up'' and meshing procedure that we used to construct a piece of a KAM tube can be applied to short pieces of the 2D stable and unstable manifolds \cite{bk10}. The extension of this approach to the blinking spherical flow does not appear straightforward at this time as the dynamics are restricted to \emph{curved surfaces}, making the meshing step difficult. However, if the latter complication could be overcome, studying the shape and structure of lobes (volumes enclosed by the intersections of the 2D stable and unstable manifolds) would also be of interest, as this could provide new ideas about 3D lobe dynamics of volume preserving maps \cite{lm09,mm11}.

Another avenue of future theoretical work could be to study \emph{resonant periodic trajectories}. That is, since we are able to compute analytically trajectories in the blinking spherical tumbler,  the results in \cite[\S2.3]{smow08} could be extended to the 3D context to find the radial position and angular extent of islands of a given period. Then, it may be possible to arrive at analytical results of how these change with the system's parameters. 

Additionally, the issue of granular segregation has not been addressed here, but it is natural to ask whether some of the dynamics that have been observed affect segregation of particles of different sizes (or densities) in any way, leading to connections such as those made for quasi-2D non-circular tumblers \cite{mclo06}. For example, in the cases where the choice of protocol parameters leads to pathological behaviors, do these exceptional mixing cases correspond to exceptional segregation cases (e.g., possibly no segregation or reverse segregation)? Since KAM regions correlate with lobes in segregation patterns, what are the implications of the results in \S\ref{sec:sym_islands} for ``minimizing'' segregation?

Finally, the more practical aspects of this problem should also be investigated. For example, does changing the fill fraction to less-than-half or more-than-half full affect the dynamics of the flow in a nontrivial way? It is already evident that changing the shape of the container leads to significantly different transport properties, as initial experimental \cite{jfg03} and combined theoretical--experimental \cite{mlo07} studies of blinking cubical tumbler flows have shown. Such cubical tumblers present the intriguing possibility of a system exhibiting two types of vanishing-flowing-layer dynamics known as streamline jumping \cite{col10} due to the non-circular cross-sections of the tumbler and cutting and shuffling \cite{smow08,jlosw10,jcol12} due to the changing axis of rotation.

\appendix
\section{Exact particle trajectory solutions for rotation about a single axis}
\label{sec:exact_2d_sol}
\begin{figure}
\centerline{\includegraphics[width=0.4\textwidth]{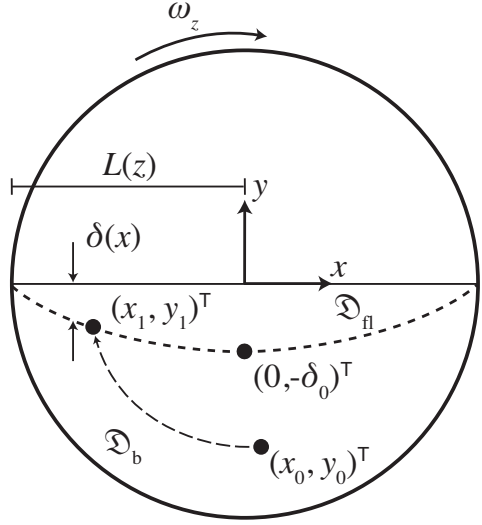}}
\caption{Diagram of a $z=const.$ cross-section of a $50\%$ full spherical tumbler rotating about the $z$-axis (out of the page) that has been rotated backwards by the dynamic angle of repose so the surface of the flowing layer is horizontal, showing the coordinate system and notation. The initial condition of a sample pathline is shown in the bulk (as in \S\ref{sec:traj_case1}) and denoted by $(x_0,y_0)^\top$; its intersection with the flowing layer is denoted by $(x_1,y_1)^\top$.}%
\label{fig:circ_tumb}
\end{figure}

For rotation about a single axis, the ordinary differential equations that define the pathlines of granular flow in a spherical tumbler are an \emph{integrable dynamical system} \cite{smow08,sow06}. Thus, it is expected that the equations can be (formally or directly) integrated to give a closed-form solution for tracer trajectories. Here, we assume the rotation is about the $z$-axis; the case of rotation about the $x$-axis can be obtained by formally exchanging $z$ and $x$ in all expressions. It suffices to consider a cross-section for an arbitrary fixed $z$, as shown in Fig.~\ref{fig:circ_tumb}, and determine $x(t)$ and $y(t)$ (with the $z$-dependence suppressed but understood). Since $z$ is fixed, let $\delta(x) := \delta_z(x,z) = \epsilon \sqrt{L^2 - x^2}$ for simplicity of notation.

We show that the trajectory in the bulk of the material rotating as a solid body [Eqs.~\eqref{eq:xoft_sb} \& \eqref{eq:yoft_sb} below] and the trajectory in the flowing layer [Eqs.~\eqref{eq:xoft_fl} \& \eqref{eq:yoft_fl} below] can be pieced together into a global trajectory in an arbitrary cross-section of the spherical tumbler rotating about a single axis. Composing the latter with the analogous rotation about the $x$-axis is straightforward. Thus, these cross-sectional solutions allow for exact tracking of a 3D particle trajectory in the blinking spherical tumbler.

A pathline that starts at an arbitrary point $(x_0,y_0)^\top$ in the filled portion of a $z=const.$ tumbler cross-section, $\mathfrak{D} = \{(x,y)^\top \,|\, x^2 + y^2 < L^2,\; y < 0\}$, is defined as a solution to Eqs.~\eqref{eq:vel_z_dxdt} and \eqref{eq:vel_z_dydt}. Since the right-hand sides of the latter are defined differently depending on whether $(x,y)^\top$ is in the bulk, $\mathfrak{D}_\mathrm{b} = \{(x,y)^\top \,|\, x^2 + y^2 < L^2,\; y < -\delta(x) \}$, or in the flowing layer, $\mathfrak{D}_\mathrm{fl} = \{(x,y)^\top \,|\, x^2 + y^2 < L^2,\; -\delta(x) < y \}$, we consider each case separately.

\subsection{Case 1: initial condition in the bulk}
\label{sec:traj_case1}
Assume that $y_0 < -\delta(x_0)$, i.e., the trajectory is seeded in the bulk. When it reaches the flowing layer's boundary, the point of intersection being $(x_1,y_1)^\top$, it must satisfy the two simultaneous equations:
\begin{subequations}\label{eq:simul_eq_x1}\begin{align}
x_1^2 + y_1^2 &= x_0^2 + y_0^2,\label{eq:simul_eq_x1_a}\\
y_1 &= -\delta(x_1) = \epsilon\sqrt{L^2-x_1^2},\label{eq:simul_eq_x1_b}
\end{align}\end{subequations}
Equation~\eqref{eq:simul_eq_x1_a} follows from the fact that the value of the streamfunction $\psi_\mathrm{b} = \tfrac{1}{2}(x^2+y^2)$ [second case of Eq.~\eqref{eq:streamfunction}] is constant along a trajectory in the bulk (and thus expresses the fact that trajectories under solid body rotation are circular arcs), while Eq.~\eqref{eq:simul_eq_x1_b} states that the point $(x_1,y_1)^\top$ is on the flowing layer interface. Solving Eqs.~\eqref{eq:simul_eq_x1} for $x_1$, we obtain
\begin{equation}
x_1 = -\sqrt{\frac{\left(x_0^2+y_0^2\right) - \epsilon^2 L^2}{1 - \epsilon^2}},
\label{eq:x1_from_sbr}
\end{equation}
where the negative square root is taken because we seek the point of entry into the flowing layer, which must be at $x_1<0$ by virtue of the fact that we have assumed clockwise rotation of the tumbler.

Next, we construct the portion of this trajectory in the bulk by finding an expression for $y$ as a function of $x$ using the constancy of $\psi_\mathrm{b}$ and the initial condition:
\begin{equation}
\frac{1}{2}\left[x^2 + y(x)^2\right] = \frac{1}{2}\left(x_0^2 + y_0^2\right)\qquad\Rightarrow\qquad y(x) = -\sqrt{x_0^2 + y_0^2 - x^2},
\label{eq:y_of_x}
\end{equation}
where we took the negative square root because $y<0$ by construction. Separating variables in the second case (corresponding to motion in the bulk) of Eq.~\eqref{eq:vel_z_dxdt}, we have $\mathrm{d}t = \mathrm{d}x/y(x)$, which can be integrated for $t(x)$ to yield
\begin{equation}
t_\mathrm{b}(x) = -\sin^{-1}\left(\frac{x}{\sqrt{x_0^2+y_0^2}}\right) + \sin^{-1}\left(\frac{x_0}{\sqrt{x_0^2+y_0^2}}\right).
\label{eq:y_of_x_of_x}
\end{equation}
Obviously $|x|\le \sqrt{x^2 + y^2} = \sqrt{x_0^2+y_0^2}$ while in the bulk, so the $\sin^{-1}$ above is real.
Inverting Eq.~\eqref{eq:y_of_x_of_x}, we obtain
\begin{equation}
x(t) = \sqrt{x_0^2 + y_0^2}\sin\left[-t + \sin^{-1}\left(\frac{x_0}{\sqrt{x_0^2+y_0^2}}\right)\right].
\label{eq:xoft_sb}
\end{equation}
Then, from the second relation in Eq.~\eqref{eq:y_of_x} and the identity $1-\sin^2\zeta = \cos^2\zeta$, we have
\begin{equation}
y(t) = -\sqrt{x_0^2 + y_0^2}\cos\left[-t + \sin^{-1}\left(\frac{x_0}{\sqrt{x_0^2+y_0^2}}\right)\right].
\label{eq:yoft_sb}
\end{equation}
Unsurprisingly, Eqs.~\eqref{eq:xoft_sb} and \eqref{eq:yoft_sb} simply describe solid body rotation.

\subsection{Case 2a: initial condition on the flowing layer interface}
\label{sec:infl}
Now, assume that $y_0 = -\delta(x_0)$, i.e., the trajectory is seeded on the interface between the flowing layer and the bulk. In this case, we first find $y(x)$ using the constancy of $\psi_\mathrm{fl}$ from the first case (corresponding to motion in the flowing layer) of Eq.~\eqref{eq:streamfunction}. The latter can be rewritten as a quadratic equation in $y$:
\begin{equation}
y^2 + 2\delta(x)y - 2 \left[\delta(x_0)y_0 + \tfrac{1}{2}y_0^2\right] = 0,
\end{equation}
which has two solutions:
\begin{equation}
y_\pm(x) = -\delta(x) \pm \sqrt{\delta^2(x) + 2\left[\delta(x_0)y_0 + \tfrac{1}{2}y_0^2\right]}.
\label{eq:yofx_fl}
\end{equation}

We have assumed that $(x_0,y_0)^\top\in\mathfrak{D}_\mathrm{fl}$, so $y(x) \ge -\delta(x)$ $\forall x$ as long as the trajectory remains in the flowing layer. Therefore, $y_+$ is the root we seek. Additionally, we must check that the expression under the square root in Eq.~\eqref{eq:yofx_fl} is always positive. This follows from the fact that the quadratic function of $y_0$ given by $\zeta(y_0) := 2\left[\delta(x_0)y_0 + \tfrac{1}{2}y_0^2\right]$ has an extremum at $y_0 = -\delta(x_0)$, which is a minimum because $\zeta''(y_0) = 2 > 0$ $\forall y_0$. Also, $\zeta(y_0)\le 0$ because $y_0 \ge -\delta(x_0)$ by assumption. Therefore, the least that the quantity under the square root in Eq.~\eqref{eq:yofx_fl} can be is $\delta^2(x)-\delta^2(x_0)$ at $y_0 = -\delta(x_0)$. For instance, if take $(x_0,y_0)^\top$ to lie on the curve $y_0 = -\delta(x_0)$, by symmetry of the flowing layer and the trajectory therein, we have $|x| \le |x_0|$ and $\delta(x) > \delta(x_0)$ $\forall x$ along the pathline, and $y_+$ is always real. Thus, the expression for $y_+$, as written, is always well-defined for pathlines starting at the interface between the flowing layer and the bulk. We return to the case when $(x_0,y_0)^\top$ does \emph{not} lie on the flowing layer boundary [i.e., $y_0 > -\delta(x_0)$] in \S\ref{sec:bookkeeping}.

For convenience, let us define the constant
\begin{equation}
\varkappa = \varkappa(x_0,y_0) := \frac{2}{\epsilon^2}\left[\delta(x_0)y_0 + \tfrac{1}{2}y_0^2\right],
\label{eq:kappa}
\end{equation}
which is a function of only $\epsilon$ and the initial condition, hence it is fixed for each pathline. Then, we eliminate $y$ [using $y_+$ from Eq.~\eqref{eq:yofx_fl}] from the first case (corresponding to motion in the flowing layer) of Eq.~\eqref{eq:vel_z_dxdt} to obtain, after separating variables and simplifying,
\begin{equation}
\rd t = \frac{\rd x}{\epsilon^{-2}\left[\delta(x) -\delta(x) + \sqrt{\delta^2(x) + \epsilon^2\varkappa} \right]} = \frac{\mathrm{d}x}{\epsilon^{-1}\sqrt{L^2 + \varkappa - x^2}}.
\end{equation}
This is easily integrated for $t(x)$, as in \S\ref{sec:traj_case1} above, arriving at
\begin{equation}
t_\mathrm{fl}(x) = \epsilon \sin^{-1}\left(\frac{x}{\sqrt{L^2 + \varkappa}}\right) - \epsilon \sin^{-1}\left(\frac{x_0}{\sqrt{L^2 + \varkappa}}\right),
\label{eq:tf}
\end{equation}
for a trajectory starting at the flowing layer interface.
Based on the discussion regarding $y_\pm$, $\varkappa \ge - \epsilon^{-2}\delta^2(x_0)$ $\Rightarrow$ $\sqrt{L^2 + \varkappa} \ge \sqrt{L^2 - \epsilon^{-2}\delta^2(x_0)} = |x_0|$, thus $\sqrt{L^2 + \varkappa}$ is always real. Also, since $|x| \le |x_0| \le \sqrt{L^2 + \varkappa}$ in the flowing layer, the $\sin^{-1}$ above is also real-valued. Inverting Eq.~\eqref{eq:tf}, we obtain
\begin{equation}
x(t) = \sqrt{L^2 + \varkappa} \sin\left[\frac{t}{\epsilon} + \sin^{-1}\left(\frac{x_0}{\sqrt{L^2 + \varkappa}}\right)\right].
\label{eq:xoft_fl}
\end{equation}

Finally, we find $y(t)$ using the $y_+$ from Eq.~\eqref{eq:yofx_fl}:
\begin{equation}
y(t) =  -\epsilon\sqrt{L^2 - x(t)^2} + \sqrt{\epsilon^2\left[L^2 - x(t)^2\right] + \epsilon^2\varkappa},
\label{eq:yt_fn_xt}
\end{equation}
where the second square root is real provided $y_0=-\delta(x_0)$.
Substituting the expressing for $x(t)$ from Eq.~\eqref{eq:xoft_fl} and using the identity $1-\sin^2\zeta = \cos^2\zeta$, Eq.~\eqref{eq:yt_fn_xt} becomes
\begin{multline}
y(t) = -\epsilon\sqrt{L^2 - (L^2 + \varkappa)\sin^2\left[\frac{t}{\epsilon} + \sin^{-1}\left(\frac{x_0}{\sqrt{L^2 + \varkappa}}\right)\right]}\\
+ \epsilon \sqrt{L^2 + \varkappa}\cos\left[\frac{t}{\epsilon} + \sin^{-1}\left(\frac{x_0}{\sqrt{L^2 + \varkappa}}\right)\right].
\label{eq:yoft_fl}
\end{multline}

\subsection{Case 2b: initial condition in the flowing layer}
\label{sec:bookkeeping}
Now, assume that $y_0 > -\delta(x_0)$, i.e., the trajectory is seeded in the flowing layer. Using the constancy of $\psi_\mathrm{fl}$ [from the first case of Eq.~\eqref{eq:streamfunction}] along a pathline in the flowing layer and $y_1\equiv-\delta(x_1)$ [from the definition that $(x_1,y_1)^\top$ is the point at which the trajectory \emph{enters} the flowing layer from the bulk] yields
\begin{equation}
\frac{1}{\epsilon^{2}}\left[ \delta(x_0)y_0 + \tfrac{1}{2}y_0^2\right] = \frac{1}{\epsilon^2}\left[-\delta^2(x_1) + \tfrac{1}{2}\delta^2(x_1)\right] = \frac{1}{2}\left(x_1^2 - L^2\right).
\end{equation}
Thus, we find what $x_1$ \emph{should be} for this trajectory:
\begin{equation}
x_1 = -\sqrt{L^2 + \frac{2}{\epsilon^{2}}\left[\delta(x_0)y_0 + \tfrac{1}{2}y_0^2\right]},
\label{eq:x1_from_fl}
\end{equation}
where the negative root is taken because our convention is that $x_1<0$. Using the solution in \S\ref{sec:infl}, we can ``trace'' the trajectory from the $(x_1,y_1)^\top$, as computed from Eq.~\eqref{eq:x1_from_fl}, but time ``starts'' at the given $x_0$ ($|x_0| < |x_1|$). Therefore, we need to find by how much we should advance (in time) the evaluation of the solution. To this end, we use Eq.~\eqref{eq:tf} with $(x_1,y_1)^\top$ in place of $(x_0,y_0)^\top$ and $(x_0,y_0)^\top$ as the final point $(x,y)^\top$:
\begin{equation}
\mathfrak{t}_\text{fl so far} = \epsilon \left[\sin^{-1}\left(\frac{x_0}{\sqrt{L^2 + \varkappa(x_1,y_1)}}\right) - \sin^{-1}\left(\frac{x_1}{\sqrt{L^2 + \varkappa(x_1,y_1)}}\right)\right].
\end{equation}
Equation~\eqref{eq:kappa} gives $\varkappa(x_1,y_1) = {2}\epsilon^{-2}\left[-\delta^2(x_1) + \tfrac{1}{2}\delta^2(x_1) \right] = x_1^2 - L^2$, whence
\begin{equation}
\mathfrak{t}_\text{fl so far} = \epsilon \left[\sin^{-1}\left(\frac{x_0}{|x_1|}\right) + \frac{\pi}{2}\right].
\end{equation}
The minus sign in front of $\pi/2$ is gone because $x_1/\sqrt{x_1^2} = -1$, i.e., we must take the negative branch of the square root because $x_1<0$ by definition.
Therefore, if our trajectory starts \emph{in} the flowing layer, we find the point $(x_1,y_1)^\top$ at which it would have entered into the flowing layer, and then we use the solution given in \S\ref{sec:infl} evaluated at $t + \mathfrak{t}_\text{fl so far}$ to account for the fact that we start at a later point along this trajectory.

\subsection{Period of a trajectory starting at an arbitrary point in $\mathfrak{D}$}

We begin by finding the point of intersection of the trajectory and the flowing layer, i.e., $(x_1,y_1)^\top$, from either Eq.~\eqref{eq:x1_from_sbr} or Eq.~\eqref{eq:x1_from_fl}. This now works for trajectories starting at an arbitrary point both in \emph{either} the bulk $\mathfrak{D}_\mathrm{b}$ or the flowing layer $\mathfrak{D}_\mathrm{fl}$. The period $T$ is the time it takes to go around any closed streamline (in the \emph{clockwise} direction), hence it is equal to the sum of the time spent in the bulk [Eq.~\eqref{eq:y_of_x_of_x} with $(-x_1,y_1)^\top$ as starting point and $(x_1,y_1)^\top$ as the endpoint] and the time spent in the flowing layer [Eq.~\eqref{eq:tf} with $(x_1,y_1)^\top$ as starting point and $(-x_1,y_1)^\top$ as the endpoint], which gives
\begin{equation}
T = -2\sin^{-1}\left(\frac{x_1}{\sqrt{x_1^2 + y_1^2}}\right) + \epsilon\pi = -2\sin^{-1}\left(\frac{x_1}{\sqrt{(1-\epsilon^2)x_1^2 + \epsilon^2 L^2}}\right) + \epsilon\pi,
\label{eq:period}
\end{equation}
where $x_1$ is given by either Eq.~\eqref{eq:x1_from_sbr} [for $(x_0,y_0)^\top\in\mathfrak{D}_\mathrm{b}$] or Eq.~\eqref{eq:x1_from_fl} [for $(x_0,y_0)^\top\in\mathfrak{D}_\mathrm{fl}$].
Using the identity $\sin^{-1}\left({x_1}\Big/{\sqrt{x_1^2+y_1^2}}\right) = -\tan^{-1}\left({x_1}/{y_1}\right)$ for $x_1,y_1<0$ in the first term of Eq.~\eqref{eq:period} makes this expression similar to the special case of $x_0=0$ in \cite[Appendix~A]{smow08}. Notice that the time spent in the flowing layer is $\epsilon\pi$ for \emph{all} streamlines, in agreement with experiments \cite{zduol12}.

\section{Dependence of the bowl depths $\mathfrak{c}_1$, $\mathfrak{c}_3$ on the rotation angles $\theta_z$, $\theta_x$}
\label{sec:c1_ineq_c3}

Here, the goal is to determine conditions on the parameters $\theta_z$, $\theta_x$, $\epsilon_z$ and $\epsilon_x$ that guarantee $\mathfrak{c}_1<\mathfrak{c}_3$. The latter inequality is relevant to the results on the existence of period-one points in the bulk of invariant surfaces in the symmetric case (\S\ref{sec:case2_shells}). To this end, we wish to verify the inequality
\begin{equation}
\mathfrak{c}_1 < \mathfrak{c}_3 \qquad\Leftrightarrow\qquad \frac{\epsilon_z^2(1+\tan^2\varphi_z)}{1+ \epsilon_z^2\tan^2\varphi_z} < \frac{\epsilon_x^2(1+\tan^2\varphi_x)}{1+\epsilon_x^2\tan^2\varphi_x},
\label{eq:c2c4_inequal}
\end{equation}
where we recall that $\mathfrak{c}$ and $\varphi$ are defined in Eq.~\eqref{eq:c_bowls}. To make any claims about the latter inequality we must know something about the extrema of the function involved. First, note that we are considering the symmetric case so that $\epsilon_z = \epsilon_x=\epsilon$ is given and fixed and also suppose that $\theta_z=\theta_x = \theta$ ($\Rightarrow\varphi_z = \varphi_x = \varphi$). Then, to extremize the function in the inequality in Eq.~\eqref{eq:c2c4_inequal} we take its derivative and set it equal to zero:
\begin{multline}
\frac{\partial}{\partial \theta} \left( \frac{\epsilon^2\{1+\tan^2[( \theta - \epsilon\pi)/2]\}}{1+\epsilon^2\tan^2[( \theta - \epsilon\pi)/2]} \right) 
= \frac{\epsilon^2\tan[( \theta - \epsilon\pi)/2]\sec^2[( \theta - \epsilon\pi)/2]}{1+\epsilon^2\tan^2[( \theta - \epsilon\pi)/2]}\\ - \frac{\epsilon^4\{1+\tan^2[( \theta - \epsilon\pi)/2]\}\tan[( \theta - \epsilon\pi)/2]\sec^2[( \theta - \epsilon\pi)/2]}{\{1+\epsilon^2\tan^2[( \theta - \epsilon\pi)/2]\}^2}=0.
\label{eq:to_minimize}
\end{multline}
First, the easy case is that of the numerator vanishing:
\begin{multline}
\epsilon^2\tan[( \theta - \epsilon\pi)/2]\sec^2[( \theta - \epsilon\pi)/2]\{1 + \epsilon^2\tan^2[( \theta - \epsilon\pi)/2]\}\\ - \epsilon^4\{1+\tan^2[( \theta - \epsilon\pi)/2]\}\tan[( \theta - \epsilon\pi)/2]\sec^2[( \theta - \epsilon\pi)/2] = 0.
\end{multline}
This results in two possibilities, either
\begin{equation}
\tan[( \theta - \epsilon\pi)/2] = 0 \quad\Rightarrow\quad (\theta - \epsilon\pi)/2 = n\pi,\quad n\in\mathbb{Z},
\end{equation}
or
\begin{equation}
\{1 + \epsilon^2\tan^2[( \theta - \epsilon\pi)/2]\} - \epsilon^2\{1+\tan^2[( \theta - \epsilon\pi)/2]\} = 0 \quad\Rightarrow\quad \epsilon = \pm 1.
\end{equation}
Since $\epsilon \ge 0$, and under the above-mentioned physical restriction $\epsilon \le 0.5$, it follows that $\epsilon=\pm1$ are extraneous (unphysical) solutions. Second, we must not neglect the case when the quantity in Eq.~\eqref{eq:to_minimize} goes to zero due to the denominator going to $\infty$, i.e., when
\begin{equation}
\tan[(\theta-\epsilon\pi)/2] \to \infty \quad\Rightarrow\quad (\theta-\epsilon\pi)/2 \to \left(m+\tfrac{1}{2}\right)\pi,\quad m\in\mathbb{Z}.
\end{equation}

Since all quantities are periodic in $\theta$ we can consider just the case $n=m=0$ without loss of generality. Then, the expression on each side of the inequality in Eq.~\eqref{eq:c2c4_inequal} above attains its minimum at $\varphi = 0$ and its maximum at $\varphi = \pi/2$. We can restate these conditions in terms of $\theta$ as: minimum at $\theta = \epsilon\pi$ and maximum at $\theta = (1+\epsilon)\pi$. It follows that the minimum values are $\mathfrak{c}_{1} = \mathfrak{c}_{3} = \epsilon^2$ (at $\theta = \epsilon\pi$) and the maximum values are $\mathfrak{c}_{1} = \mathfrak{c}_{3} = 1$ (at $\theta = (1+\epsilon)\pi$).

We have shown that the quantity $\mathfrak{c}(\theta,\epsilon)$ in the inequality $\mathfrak{c}(\theta_z,\epsilon_z) \equiv \mathfrak{c}_1 < \mathfrak{c}_3 \equiv \mathfrak{c}(\theta_x,\epsilon_x)$ that we seek to establish is decreasing for $\theta\in\big[0,\epsilon\pi \big)$, increasing for $\theta\in \big(\epsilon\pi,(1+\epsilon)\pi \big)$ and decreasing again for $\theta\in \big((1+\epsilon)\pi,2\pi \big]$. Hence, the desired inequality is easy to verify in the following special cases:
\begin{equation}
\mathfrak{c}_1 \le \mathfrak{c}_3 \quad\text{for}\quad  \begin{cases} 0 \le \theta_x \le \theta_z \le \epsilon\pi,\\ \epsilon\pi \le \theta_z \le \theta_x \le (1+\epsilon)\pi,\\ (1+\epsilon)\pi \le \theta_x \le \theta_z \le 2\pi, \end{cases}
\end{equation}
and
\begin{equation}
\mathfrak{c}_1 > \mathfrak{c}_3 \quad\text{for}\quad  \begin{cases} 0 \le \theta_z < \theta_x \le \epsilon\pi,\\ \epsilon\pi \le \theta_x < \theta_z \le (1+\epsilon)\pi,\\ (1+\epsilon)\pi \le \theta_z < \theta_x \le 2\pi. \end{cases}
\end{equation}
When $\theta_z$ and $\theta_x$ ``straddle'' the extrema, the values of $\mathfrak{c}_1$ and $\mathfrak{c}_3$ will have to be found and compared to each other in order to establish an inequality. Also, it would appear we cannot say much about the inequality in the non-symmetric case when $\epsilon_z \ne \epsilon_x$.

\section{Computing manifolds by iterating a fundamental domain}
\label{app:manifolds}

To begin, we denote by $\bm{\Phi}$ the Poincar\'e map\footnote{Sometimes referred to as the \emph{Liouvillian map} for such 3D flows \cite{cfp96}.} of the blinking spherical tumbler flow, i.e., $\bm{x} = \bm{\Phi}(\bm{X})$ is the location of a material point $\bm{X}$ (from the filled portion of the tumbler) after a rotation about each of the two axes is complete (one iteration of linked twist map).

The fixed points $\{\bm{X}^*\}$ of $\bm{\Phi}$ are the period-one points discussed in detail in \S\ref{sec:location_p1_pts} and \ref{sec:location_p1_pts_surf}. In the symmetric case, when motion is restricted to 2D surfaces in $\mathbb{R}^3$, the normally hyperbolic points have 1D stable and unstable manifolds on said 2D surfaces. We can construct approximations to finite-length segments of these manifolds by iterating an infinitesimal piece of manifold, termed a \emph{fundamental domain}, through repeated applications of the map $\bm{\Phi}$ (for the case of the unstable manifold) or the map $\bm{\Phi}^{-1}$ (for the case of the stable manifold) \cite[Chapter~6]{pc89}.

To construct a fundamental domain suppose we know the Jacobian $\mathrm{D}\bm{\Phi}$ of the Poincar\'e map and its eigenvalues $\{\lambda_1,\lambda_2,\lambda_2\}$ and eigenvectors $\{\bm{E}_1,\bm{E}_2,\bm{E}_3\}$ at the fixed point, then the line segment connecting $\bm{X}^*$ to $\overline{\bm{X}}$ is a fundamental domain of the manifold whenever $\overline{\bm{X}}$ is a small perturbation of $\bm{X}^*$ in one of the eigendirections:
\begin{equation}
\overline{\bm{X}} = \bm{X}^* + \alpha\bm{E}_i.
\end{equation}
Here, $\alpha$ is chosen to ensure that the manifold's curvature over the given distance from the fixed point is sufficiently small so that the line segment connecting $\bm{X}^*$ and $\overline{\bm{X}}$ is indeed close to the true manifold.

We approximate the linear eigenspace of $\bm{\Phi}$ at $\bm{X}^*$ using finite differences by constructing a uniform Cartesian grid $\{\bm{X}_{i,j,k}\}$ with $\{\bm{x}_{i,j,k}\}$ being its image after $\bm{\Phi}$ is applied. The Jacobian of $\bm{\Phi}$ at $\bm{X}^*$ can be approximated from this grid via $\mathrm{D}\bm{\Phi} = \mathrm{D}\bm{\Phi}_h + \mathcal{O}(h)$ (see, e.g, \cite{lr10}), where
\begin{equation}
\mathrm{D}\bm{\Phi}_h = \frac{1}{2h}\begin{pmatrix} {x_{i+1,j,k}-x_{i-1,j,k}} & {x_{i,j+1,k}-x_{i,j-1,k}} & {x_{i,j,k+1}-x_{i,j,k-1}}\\ {y_{i+1,j,k}-y_{i-1,j,k}} & {y_{i,j+1,k}-y_{i,j-1,k}} & {y_{i,j,k+1}-y_{i,j,k-1}}\\ {z_{i+1,j,k}-z_{i-1,j,k}} & {z_{i,j+1,k}-z_{i,j-1,k}} & {z_{i,j,k+1}-z_{i,j,k-1}} \end{pmatrix}
\end{equation}
for the case of a grid constructed by uniform displacements from $\bm{X}^*$: i.e., $\bm{X}_{i\pm1,j,k} = \bm{X}^* \pm h \hat{\bm{\imath}}$, $\bm{X}_{i,j\pm1,k} = \bm{X}^* \pm h \hat{\bm{\jmath}}$ and $\bm{X}_{i,j,k\pm1} = \bm{X}^* \pm h \hat{\bm{k}}$ with $h\ll 1$ ($h = 10^{-6}$ was used in our calculations). 

A check on the accuracy of the numerical calculation of $\mathrm{D}\bm{\Phi}_h$ is to verify that its eigenvalues satisfy $\lambda_1\lambda_2\lambda_3 \approx 1$ as required by the fact that the flow is divergence-free and thus the Poincar\'e map is volume preserving. In addition, we must also have $\lambda_1 =1$ with the eigenvector $\bm{E}_1$ being the unit tangent vector to the invariant curve of period-one (fixed) points. For those points along the normally hyperbolic invariant curve, $\lambda_{2,3}$ are both real and such that $\lambda_2 = 1/\lambda_3$. Along the normally elliptic invariant curve, $\lambda_{2,3}$ are complex conjugates of unit magnitude. 

An attempt was also made to implement an adaptive algorithm for computing a finite-length portion of  manifolds \cite[Chapter~6]{pc89}, however, the algorithm failed at the interface between the flowing layer and the bulk, where the manifolds are continuous but not differentiable.

\subsection*{Acknowledgments} I.C.C.\ thanks Michel Speetjens for clarifying remarks and helpful discussions about his related research and Paul Park for a careful reading of the manuscript. The authors acknowledge discussions with Stephen Wiggins regarding this work and, specifically, his suggestion of the argument in \S\ref{sec:radial_disp}. I.C.C.\ and R.S.\ would also like to acknowledge the hospitable environment of the 2011 and 2013 SIAM Conferences on Applications of Dynamical Systems in Snowbird, Utah, which provided a catalyzer for much of this collaboration. During the completion of the revision of this manuscript, I.C.C.\ was supported by LANL/LDRD Program through a Feynman Distinguished Fellowship at Los Alamos National Laboratory, which is operated by Los Alamos National Security, L.L.C.\ for the National Nuclear Security Administration of the U.S.\ Department of Energy under Contract No.\ DE-AC52-06NA25396.


\bibliographystyle{siam}
\bibliography{SIADS_093407R}

\begin{thebibliography}{10}

\bibitem{al02}
{\sc A.~Akonur and R.~M. Lueptow}, {\em Chaotic mixing and transport in wavy
  {Taylor--Couette} flow}, Physica D, 167 (2002), pp.~183--196.

\bibitem{agpvm99}
{\sc P.~D. Anderson, O.~S. Galaktionov, G.~W.~M. Peters, F.~N. {van de Vosse},
  and H.~E.~H. Meijer}, {\em Analysis of mixing in three-dimensional
  time-periodic cavity flows}, J. Fluid Mech., 386 (1999), pp.~149--166.

\bibitem{at09}
{\sc I.~S. Aranson and L.~S. Tsimring}, {\em Granular Patterns}, Oxford
  University Press, New York, 2009.

\bibitem{a84}
{\sc H.~Aref}, {\em Stirring by chaotic advection}, J. Fluid Mech., 143 (1984),
  pp.~1--21.

\bibitem{a65}
{\sc V.~I. Arnold}, {\em Sur la topologie des \'{e}coulements stationnaires des
  fluides parfaits}, C. R. Acad. Sci. Paris A, 261 (1965), pp.~17--20.

\bibitem{a89}
\leavevmode\vrule height 2pt depth -1.6pt width 23pt, {\em Mathematical Methods
  of Classical Mechanics}, vol.~60 of Graduate Texts in Mathematics,
  Springer-Verlag, New York, 2nd~ed., 1989.

\bibitem{alsm04}
{\sc P.~E. Arratia, J.~P. Lacombe, T.~Shinbrot, and F.~J. Muzzio}, {\em
  Segregated regions in continuous laminar stirred tank reactors}, Chem. Eng.
  Sci., 59 (2004), pp.~1481--1490.

\bibitem{ak95}
{\sc P.~Ashwin and G.~P. King}, {\em Streamline topology in eccentric {Taylor}
  vortex flow}, J. Fluid Mech., 285 (1995), pp.~215--247.

\bibitem{bajer}
{\sc K.~Bajer}, {\em {Hamiltonian} formulation of the equations of streamlines
  in three-dimensional steady flows}, Chaos Solitons Fractals, 4 (1994),
  pp.~859--911.

\bibitem{bm90}
{\sc K.~Bajer and H.~K. Moffatt}, {\em On a class of steady confined {Stokes}
  flows with chaotic streamlines}, J. Fluid Mech., 212 (1990), pp.~337--363.

\bibitem{blw94}
{\sc D.~Beigie, A.~Leonard, and S.~Wiggins}, {\em Invariant manifold templates
  for chaotic advection}, Chaos Solitons Fractals, 4 (1994), pp.~749--868.

\bibitem{betal89}
{\sc V.~V. Beloshapkin, A.~A. Chernikov, M.~{\relax Ya}. Natenzon, B.~A.
  Petrovichev, R.~Z. Sagdeev, and G.~M. Zaslavsky}, {\em Chaotic streamlines in
  pre-turbulent states}, Nature, 337 (1989), pp.~133--137.

\bibitem{bk10}
{\sc M.~Branicki and A.~D. {Kirwan Jr.}}, {\em Stirring: {The Eckart} paradigm
  revisited}, Int. J. Eng. Sci., 48 (2010), pp.~1027--1042.

\bibitem{b11}
{\sc L.~E. Brouwer}, {\em \"{U}ber abbildung von {M}annigfaltigkeiten}, Math.
  Annalen, 71 (1911), pp.~97--115.

\bibitem{bb99}
{\sc M.~D. Bryden and H.~Brenner}, {\em Mass-transfer enhancement via chaotic
  laminar flow within a droplet}, J. Fluid Mech., 379 (1999), pp.~319--331.

\bibitem{cw91}
{\sc R.~Camassa and S.~Wiggins}, {\em Chaotic advection in a
  {Rayleigh--B\'enard} flow}, Phys. Rev. A, 43 (1991), pp.~774--797.

\bibitem{cfp96}
{\sc J.~H.~E. Cartwright, M.~Feingold, and O.~Piro}, {\em Chaotic advection in
  three-dimensional unsteady incompressible laminar flow}, J. Fluid Mech., 316
  (1996), pp.~259--284.

\bibitem{cvcsa08}
{\sc R.~Chabreyrie, D.~Vainchtein, C.~Chandre, P.~Singh, and N.~Aubry}, {\em
  Tailored mixing inside a translating droplet}, Phys. Rev. E, 77 (2008),
  pp.~036314--1--4.

\bibitem{cyymm11}
{\sc C.-K. Chen, S.~Yan, H.~Yu, N.~Max, and K.-L. Ma}, {\em An illustrative
  visualization framework for {3D} vector fields}, Comput. Graph. Forum, 30
  (2011), pp.~1941--1951.

\bibitem{cs90}
{\sc C.-Q. Cheng and Y.-S. Sun}, {\em Existence of invariant tori in
  three-dimensional measure-preserving mappings}, J. Nonlinear Sci., 47 (1990),
  pp.~275--292.

\bibitem{c11}
{\sc I.~C. Christov}, {\em From Streamline Jumping to Strange Eigenmodes and
  Three-Dimensional Chaos: A Tour of the Mathematical Aspects of Granular
  Mixing in Rotating Tumblers}, PhD thesis, Northwestern University, Evanston,
  Illinois, June 2011.

\bibitem{col10}
{\sc I.~C. Christov, J.~M. Ottino, and R.~M. Lueptow}, {\em Streamline jumping:
  A mixing mechanism}, Phys. Rev. E, 81 (2010), p.~046307.

\bibitem{dfghms86}
{\sc T.~Dombre, U.~Frisch, J.~M. Greene, M.~H\'enon, A.~Mehr, and A.~M.
  Soward}, {\em Chaotic streamlines in the {ABC} flows}, J. Fluid Mech., 167
  (1986), pp.~353--391.

\bibitem{dbc08}
{\sc J.~Doucet, F.~Betrand, and J.~Chaouki}, {\em Experimental characterization
  of the chaotic dynamics of cohesionless particles: application to a
  {V}-blender}, Granular Matter, 10 (2008), pp.~133--138.

\bibitem{emc93}
{\sc R.~W. Easton, J.~D. Meiss, and S.~Carver}, {\em Exit times and transport
  for symplectic twist maps}, Chaos, 3 (1993), pp.~154--165.

\bibitem{fkp88}
{\sc M.~Feingold, L.~P. Kadanoff, and O.~Piro}, {\em Passive scalars,
  three-dimensional volume-preserving maps, and chaos}, J. Stat. Phys., 50
  (1988), pp.~529--565.

\bibitem{ffo07}
{\sc G.~F\'{e}lix, V.~Falk, and U.~D'Ortona}, {\em Granular flows in a rotating
  drum: the scaling law between velocity and thickness of the flow}, Eur. Phys.
  J. E, 22 (2007), pp.~25--31.

\bibitem{fkmo}
{\sc G.~O. Fountain, D.~V. Khakhar, I.~Mezi\'c, and J.~M. Ottino}, {\em Chaotic
  mixing in a bounded three-dimensional flow}, J. Fluid Mech., 417 (2000),
  pp.~265--301.

\bibitem{fo92}
{\sc J.~G. Franjione and J.~M. Ottino}, {\em Symmetry concepts for the
  geometric analysis of mixing flows}, Phil. Trans. R. Soc. Lond. A, 338
  (1992), pp.~301--323.

\bibitem{gapm03}
{\sc O.~S. Galaktionov, P.~D. Anderson, G.~W.~M. Peters, and H.~E.~H. Meijer},
  {\em Analysis and optimization of {Kenics} static mixers}, Int. Polym.
  Process, 18 (2003), pp.~138--150.

\bibitem{jfg03}
{\sc J.~F. Gilchrist}, {\em Geometric Aspects of Mixing and Segregation in
  Granular Tumblers}, PhD thesis, Northwestern University, Evanston, Illinois,
  June 2003.

\bibitem{go03}
{\sc J.~F. Gilchrist and J.~M. Ottino}, {\em Competition between chaos and
  order: Mixing and segregation in a spherical tumbler}, Phys. Rev. E, 68
  (2003), p.~061303.

\bibitem{gm}
{\sc A.~G\'{o}mez and J.~D. Meiss}, {\em Volume-preserving maps with an
  invariant}, Chaos, 12 (2002), pp.~289--299.

\bibitem{hm98}
{\sc H.~Haller and I.~Mezi\'c}, {\em Reduction of three-dimensional,
  volume-preserving flows with symmetry}, Nonlinearity, 11 (1996),
  pp.~319--339.

\bibitem{h98}
{\sc H.~Han{\ss}mann}, {\em The reversible umbilic bifurcation}, Physica D, 112
  (1998), pp.~81--94.

\bibitem{hok12}
{\sc S.~Hashimoto, R.~Osaka, and M.~Kawamata}, {\em Analysis on laminar chaotic
  mixing based on configuration of streak lobes in an impeller-agitated
  vessel}, Ind. Eng. Chem. Res., 51 (2012), pp.~6939--6947.

\bibitem{h66}
{\sc M.~H\'enon}, {\em Sur la topologie des lignes de courant dans un cas
  particulier}, C. R. Acad. Sci. Paris A, 262 (1966), pp.~312--314.

\bibitem{hkgmo99}
{\sc K.~M. Hill, D.~V. Khakhar, J.~F. Gilchrist, J.~J. McCarthy, and J.~M.
  Ottino}, {\em Segregation-driven organization in chaotic granular flows},
  Proc. Natl Acad. Sci. USA, 96 (1999), pp.~11701--11706.

\bibitem{h94}
{\sc M.~J.~M. Hill}, {\em On a spherical vortex}, Phil. Trans. R. Soc. Lond. A,
  185 (1894), pp.~213--245.

\bibitem{hk91}
{\sc D.~D. Holm and Y.~Kimura}, {\em Zero-helicity {Lagrangian} kinematics of
  three-dimensional advection}, Phys. Fluids A, 3 (1991), pp.~1033--1038.

\bibitem{hkk04}
{\sc W.~R. Hwang, K.~W. Kang, and T.~H. Kwon}, {\em Dynamical systems in pin
  mixers of single-screw extruders}, {AIChE} J., 50 (2004), pp.~1372--1385.

\bibitem{jol04}
{\sc N.~Jain, J.~M. Ottino, and R.~M. Lueptow}, {\em Effect of interstitial
  fluid on a granular flowing layer}, J. Fluid Mech., 508 (2004), pp.~23--44.

\bibitem{jta89}
{\sc S.~W. Jones, O.~M. Thomas, and H.~Aref}, {\em Chaotic advection by laminar
  flow in a twisted pipe}, J. Fluid Mech., 209 (1989), pp.~335--357.

\bibitem{jcol12}
{\sc G.~Juarez, I.~C. Christov, J.~M. Ottino, and R.~M. Lueptow}, {\em Mixing
  by cutting and shuffling {3D} granular flow in spherical tumblers}, Chem.
  Eng. Sci., 73 (2012), pp.~195--207.

\bibitem{jlosw10}
{\sc G.~Juarez, R.~M. Lueptow, J.~M. Ottino, R.~Sturman, and S.~Wiggins}, {\em
  Mixing by cutting and shuffling}, EPL, 91 (2010), p.~20003.

\bibitem{k11}
{\sc D.~V. Khakhar}, {\em Rheology and mixing of granular materials}, Macromol.
  Mater. Eng., 296 (2011), pp.~278--289.

\bibitem{kfo87}
{\sc D.~V. Khakhar, J.~G. Franjione, and J.~M. Ottino}, {\em A case study of
  chaotic mixing in deterministic flows: {The} partitioned-pipe mixer}, Chem.
  Eng. Sci., 42 (1987), pp.~2909--2926.

\bibitem{ks99}
{\sc D.~Kroujiline and H.~A. Stone}, {\em Chaotic streamlines in steady bounded
  three-dimensional {Stokes} flows}, Physica D, 130 (1999), pp.~105--132.

\bibitem{ko92}
{\sc H.~A. Kusch and J.~M. Ottino}, {\em Experiments on mixing in continuous
  chaotic flows}, J. Fluid Mech., 236 (1992), pp.~319--348.

\bibitem{kuz06}
{\sc Yu.~A. Kuznetsov}, {\em Saddle-node bifurcation}, Scholarpedia, 1 (2006),
  p.~1859.

\bibitem{lam99}
{\sc D.~J. Lamberto, M.~M. Alvarez, and F.~J. Muzzio}, {\em Experimental and
  computational investigation of the laminar flow structure in a stirred tank},
  Chem. Eng. Sci., 54 (1999), pp.~919--942.

\bibitem{lr10}
{\sc F.~Lekien and S.~D. Ross}, {\em The computation of finite-time {Lyapunov}
  exponents on unstructured meshes and for non-{Euclidean} manifolds}, Chaos,
  20 (2010), p.~017505.

\bibitem{lbcg07}
{\sc M.~Lemieux, F.~Bertrand, J.~Chaouki, and P.~Gosselin}, {\em Comparative
  study of the mixing of free-flowing particles in a {V}-blender and a
  bin-blender}, Chem. Eng. Sci., 62 (2007), pp.~1783--1802.

\bibitem{lm09}
{\sc H.~E. Lomel\'{i} and J.~D. Meiss}, {\em Resonance zones and lobe volumes
  for exact volume-preserving maps}, Nonlinearity, 22 (2009), pp.~1761--1789.

\bibitem{mycrustopen}
{\sc G.~Luigi}, {\em {MyCrustOpen}: Surface reconstruction from scattered
  points cloud}.
\newblock MATLAB File Exchange, 2009.

\bibitem{mackay}
{\sc R.~S. MacKay}, {\em Transport in {3D} volume-preserving flows}, J.
  Nonlinear Sci., 1 (1994), pp.~329--354.

\bibitem{mmsch02}
{\sc V.~S. Malyuga, V.~V. Meleshko, M.~F.~M. Speetjens, H.~J.~H. Clercx, and
  G.~J.~F. van Heijst}, {\em Mixing in the {Stokes} flow in a cylindrical
  container}, Proc. R. Soc. Lond. A, 458 (2002), pp.~1867--1885.

\bibitem{mclo06}
{\sc S.~W. Meier, S.~E. Cisar, R.~M. Lueptow, and J.~M. Ottino}, {\em Capturing
  patterns and symmetries in chaotic granular flow}, Phys. Rev. E, 74 (2006),
  p.~031310.

\bibitem{mlo07}
{\sc S.~W. Meier, R.~M. Lueptow, and J.~M. Ottino}, {\em A dynamical systems
  approach to mixing and segregation of granular materials in tumblers}, Adv.
  Phys., 56 (2007), pp.~757--827.

\bibitem{grbgh06}
{\sc G.~Metcalfe, M.~Rudman, A.~Brydon, L.~J.~W. Graham, and R.~Hamilton}, {\em
  Composing chaos: An experimental and numerical study of an open duct mixing
  flow}, AIChE J., 52 (2006), pp.~9--28.

\bibitem{mw94}
{\sc I.~Mezi\'c and S.~Wiggins}, {\em On the integrability and perturbation of
  three-dimensional fluid flows with symmetry}, J. Nonlinear Sci., 4 (1994),
  pp.~157--194.

\bibitem{mm11}
{\sc B.~A. Mosovsky and J.~D. Meiss}, {\em Transport in transitory dynamical
  systems}, SIAM J. Appl. Dyn. Syst., 10 (2011), pp.~35--65.

\bibitem{mjm05}
{\sc P.~Mullowney, K.~Julien, and J.~D. Meiss}, {\em Blinking rolls: Chaotic
  advection in a three-dimensional flow with an invariant}, SIAM J. Appl. Dyn.
  Syst., 4 (2005), pp.~159--186.

\bibitem{o89}
{\sc J.~M. Ottino}, {\em The Kinematics of Mixing: Stretching, Chaos, and
  Transport}, vol.~3 of Cambridge Texts in Applied Mathematics, Cambridge
  University Press, Cambridge, 1989.

\bibitem{pc89}
{\sc T.~S. Parker and L.~O. Chua}, {\em Practical Numerical Algorithms for
  Chaotic Systems}, Springer-Verlag, New York, 1989.

\bibitem{pakcol12}
{\sc F.~Pignatel, C.~Asselin, L.~Krieger, I.~C. Christov, J.~M. Ottino, and
  R.~M. Lueptow}, {\em Parameters and scalings for dry and immersed granular
  flowing layers in rotating tumblers}, Phys. Rev. E, 86 (2012), p.~011304.

\bibitem{psc10}
{\sc Z.~Pouransari, M.~F.~M. Speetjens, and H.~J.~H. Clercx}, {\em Formation of
  coherent structures by fluid inertia in three-dimensional laminar flows}, J.
  Fluid Mech., 654 (2010), pp.~5--34.

\bibitem{rklw90}
{\sc V.~Rom-Kedar, A.~Leonard, and S.~Wiggins}, {\em An analytical study of
  transport, mixing and chaos in an unsteady vortical flow}, J. Fluid Mech.,
  214 (1990), pp.~347--394.

\bibitem{r98}
{\sc M.~Rudman}, {\em Mixing and particle dispersion in the wavy vortex regime
  of {Taylor--Couette} flow}, AIChE J., 44 (1998), pp.~1015--1026.

\bibitem{sg88}
{\sc T.~H. Solomon and J.~P. Gollub}, {\em Chaotic particle transport in
  time-dependent {Rayleigh--B\'enard} convection}, Phys. Rev. A, 38 (1988),
  pp.~6280--6286.

\bibitem{sti03}
{\sc H.~Song, J.~D. Tice, and R.~F. Ismagilov}, {\em A microfluidic system for
  controlling reaction networks in time}, Angew. Chem. Int. Ed., 42 (2003),
  pp.~767--772.

\bibitem{sc13}
{\sc M.~F.~M. Speetjens and H.~J.~H. Clercx}, {\em Formation of coherent
  structures in a class of realistic {3D} unsteady flows}, in Fluid Dynamics in
  Physics, Engineering and Environmental Applications, J.~Klapp, A.~Medina,
  A.~Cros, and C.~A. Vargas, eds., Springer-Verlag, Berlin/Heidelberg, 2013,
  pp.~139--157.

\bibitem{sch04}
{\sc M.~F.~M. Speetjens, H.~J.~H. Clercx, and G.~J.~F. {van Heijst}}, {\em A
  numerical and experimental study on advection in three-dimensional {Stokes}
  flows}, J. Fluid Mech., 514 (2004), pp.~77--105.

\bibitem{sch06a}
\leavevmode\vrule height 2pt depth -1.6pt width 23pt, {\em Inertia-induced
  coherent structures in a time-periodic viscous mixing flow}, Phys. Fluids, 18
  (2006), pp.~083603--1--15.

\bibitem{sch06}
\leavevmode\vrule height 2pt depth -1.6pt width 23pt, {\em Merger of coherent
  structures in time-periodic viscous flows}, Chaos, 16 (2006),
  pp.~043104--1--8.

\bibitem{stir03}
{\sc J.~R. Stirling}, {\em Chaotic advection, transport and patchiness in
  clouds of pollution in an estuarine flow}, Discret. Contin. Dyn. Syst. B, 3
  (2003), pp.~263--284.

\bibitem{sns91}
{\sc H.~A. Stone, A.~Nadim, and S.~H. Strogatz}, {\em Chaotic streamlines
  inside drops immersed in steady {Stokes} flows}, J. Fluid Mech., 231 (1991),
  pp.~149--166.

\bibitem{ss05}
{\sc Z.~B. Stone and H.~A. Stone}, {\em Imaging and quantifying mixing in a
  model droplet micromixer}, Phys. Fluids, 17 (2005), pp.~063103--1--11.

\bibitem{sdamsw02}
{\sc A.~D. Stroock, S.~K.~W. Dertinger, A.~Ajdari, I.~Mezi\'c, H.~A. Stone, and
  G.~M. Whitesides}, {\em Chaotic mixer for microchannels}, Science, 295
  (2002), pp.~647--651.

\bibitem{smow08}
{\sc R.~Sturman, S.~W. Meier, J.~M. Ottino, and S.~Wiggins}, {\em Linked twist
  map formalism in two and three dimensions applied to mixing in tumbled
  granular flows}, J. Fluid Mech., 602 (2008), pp.~129--174.

\bibitem{sow06}
{\sc R.~Sturman, J.~M. Ottino, and S.~Wiggins}, {\em The Mathematical
  Foundations of Mixing}, vol.~22 of Cambridge Monographs on Applied and
  Computational Mathematics, Cambridge University Press, Cambridge, 2006.

\bibitem{vm12}
{\sc U.~Vaidya and I.~Mezi\'c}, {\em Existence of invariant tori in three
  dimensional maps with degeneracy}, Physica D, 241 (2012), pp.~1136--1145.

\bibitem{vnm06}
{\sc D.~L. Vainchtein, A.~I. Neishtadt, and I.~Mezi\'c}, {\em On passage
  through resonances in volume-preserving systems}, Chaos, 16 (2006),
  p.~043123.

\bibitem{vwg08}
{\sc D.~L. Vainchtein, J.~Widloski, and R.~O. Grigoriev}, {\em Resonant mixing
  in perturbed action-action-angle flow}, Phys. Rev. E, 78 (2008),
  pp.~026302--1--11.

\bibitem{vwg10}
\leavevmode\vrule height 2pt depth -1.6pt width 23pt, {\em Erratum: {Resonant}
  mixing in perturbed action-action-angle flow {[Phys. Rev. E 78, 026302
  (2008)]}}, Phys. Rev. E, 81 (2010), pp.~049904--1.

\bibitem{w94}
{\sc S.~Wiggins}, {\em Normally Hyperbolic Invariant Manifolds in Dynamical
  Systems}, vol.~105 of Applied Mathematical Sciences, Springer-Verlag, New
  York, 1994.

\bibitem{wig-fof}
\leavevmode\vrule height 2pt depth -1.6pt width 23pt, {\em Coherent structures
  and chaotic advection in three dimensions}, J. Fluid Mech., 654 (2010),
  pp.~1--4.

\bibitem{WightmanSim}
{\sc C.~Wightman, M.~Moakher, F.~J. Muzzio, and O.~Walton}, {\em Simulation of
  flow and mixing of particles in a rotating and rocking cylinder}, {AIChE} J.,
  44 (1998), pp.~1266--1276.

\bibitem{zduol12}
{\sc Z.~Zaman, U.~{D'Ortona}, P.~Umbanhowar, J.~M. Ottino, and R.~M. Lueptow},
  {\em Slow axial drift in three-dimensional granular tumbler flow}, Phys. Rev.
  E, 88 (2013), p.~012208.

\bibitem{zklh93}
{\sc {X.-H.} Zhao, {K.-H.} Kwek, {J.-B.} Li, and {K.-L.} Huang}, {\em Chaotic
  and resonant streamlines in the {ABC} flow}, SIAM J. Appl. Math., 53 (1993),
  pp.~71--77.

\bibitem{zstc12}
{\sc J.~Znaien, M.~F.~M. Speetjens, R.~R. Trieling, and H.~J.~H. Clercx}, {\em
  Observability of periodic lines in three-dimensional lid-driven cylindrical
  cavity flows}, Phys. Rev. E, 85 (2012), p.~066320.

\end{thebibliography}

\end{document}